\documentclass[11pt]{article}
\usepackage[a4paper,margin= 2cm]{geometry}
\usepackage{booktabs,longtable,array,multirow}
\usepackage[T1]{fontenc}
\usepackage{hyperref}
\usepackage{enumitem}
\usepackage{amssymb}

\usepackage{listings}

\def\st{\mathord{*}}
\def\-{\raisebox{.75pt}{-}}

\title{Non-Isomorphic Optimal Cayley Graphs }
\author{Francesc Comellas\\
 {\tt\small  francesc.comellas@upc.edu}}
\date{\today}

\begin{document}
\maketitle

\section*{Introduction}

We describe optimal degree-diameter Cayley graphs that are non-isomorphic to those currently reported. The graphs were identified using SageMath 10.8 programs generated within seconds by \texttt{Perplexity AI} and made operational after only minor modifications. One program uses groups from  the GAP SmallGroups Library \cite{BeEiOB05} and constructs Cayley graphs of prescribed order and diameter by examining generating sets either exhaustively (for groups of small order) or through randomized searches. Because the program often explores a large number of generating sets, some computations are extremely time-consuming and may require several days to complete. The resulting Cayley graphs have the same order as their underlying groups.

A second program searches the Cayley graph databases available at \url{https://graphsym.net/} and identifies graphs with minimum diameter for a given order and degree.

These computations revealed numerous graphs that are non-isomorphic to those listed as optimal by Marston Conder in a forthcoming paper \cite{Co26}, and previously described on the Combinatorics Wiki website \cite{Co11}. Several of the newly identified graphs exhibit distinct structural properties, including differences in girth, algebraic connectivity, automorphism groups, cycle counts, and distance distributions. As alternative optimal graphs, they may provide useful examples for future investigations of extremal problems involving degree, diameter, and girth.

We also present four new Cayley graphs that, although they have not yet been proven to be optimal, improve the current records in the table of the largest known Cayley graphs for the undirected degree-diameter problem \cite{Co11,Co26}. These correspond to the parameter pairs (9,3), (12,3), (13,2), and (13,3).

The collection of non-isomorphic graphs reported here is not exhaustive, since the search procedure relies partly on randomized methods, especially for graphs of larger order.

\medskip

\noindent\textbf{Notation.} Throughout the paper, the following abbreviations are used:
\begin{itemize}[leftmargin=1.5cm]
\item MC denotes optimal Cayley graphs listed by Marston Conder on the Combinatorics Wiki:
\url{http://combinatoricswiki.org/wiki/Description_of_optimal_Cayley_graphs_found_by_Marston_Conder}.
\item HGxxxx  denotes a graph from the House of Graphs database:
\url{https://houseofgraphs.org/}.
\item SG denotes GAP's \texttt{SmallGroup} notation. For example, SG24x10x1 represents graph 1 (of several) coming from \texttt{SmallGroup(24,10)}. 
\item P denotes cubic and tetravalent Cayley graphs from the compilation of Poto\v{c}nik, Spiga, and Verret \cite{Po13}, available at \url{https://graphsym.net/}. For example, P84x142 denotes Cayley graph number 142 among the 426 tetravalent Cayley graphs of order 84.
\end{itemize}
The included SageMath links execute code online directly from the sparse6 representation of the graphs embedded in the script. No local SageMath or Python installation is required.

%\tableofcontents
\newpage

%%%%%%%%%%%%%%%%%%%%%%%%%%%%%%%%%%%%%%%%%%%%%%%%%%%%%%%%%%%%
\section{Degree 3}
%%%%%%%%%%%%%%%%%%%%%%%%%%%%%%%%%%%%%%%%%%%%%%%%%%%%%%%%%%%%

We note that all cubic Cayley graphs with order up to 5000 (and  up to isomorphism) have been determined by Poto\v{c}nik, Spiga and Verret \cite{Po13} and can be found at \url{https://graphsym.net/}. Here be calculate the main properties of the optimal graphs that are not isomorphic to one another for each degree-diameter value and give also where can they be found in  the SmallGroups Library.

%%%%%%%%%%%%%%%%%%%%%%%%%%%%%%%
\subsection{\large  (3,3) = 14}
Reference: \url{https://sagecell.sagemath.org/?q=pjltke}
\bigskip

\noindent All 14 optimal Cayley graphs found are isomorphic to the Heawood graph / MC14. \\
\texttt{Heawood/MC14/HG1154/P14x2/SG(14,1)->14}

%%%%%%%%%%%%%%%%%%%%%%%%%%%%%%
\subsection{\large  (3,4)=24}
Reference: \url{https://sagecell.sagemath.org/?q=mzqyth}
\bigskip 

\noindent{\bf 3 non-isomorphic optimal Cayley graphs }\\
\texttt{MC24/Nauru/HG1234/P24x10/SG(24,5)->12/SG(24,8)->12/SG(24,10)->8}. All 32 graphs are isomorphic to one another.\\
\texttt{SG(24,12)->40,SG(24,13)->24,SG(24,14)->24}. Of these 88 graphs, two are not isomorphic to one another or to MC24. These are isomorphic to HG36335 or HG36333. All others are isomorphic to one of these three. 

\begin{center}
\begin{tabular}{lcccccc}
\toprule
Graph         & Avg.dist & Girth & Alg.conn.& Dom.num. & Aut.group  & Edge-trans. \\
\midrule
MC24/Nauru/HG1234/P24x10 & 2.69565 & 6 & 1.000000 &  6 & 144 & True \\
HG36335/P24x9/SG24x12x1   & 2.60870 & 6 & 1.000000 &  6 & 48 & False \\
HG36333/P24x7/SG24x12x2   & 2.86957 & 4 & 0.657077 &  7 & 24 & False \\
\bottomrule
\end{tabular}
\end{center}

\begin{center}
\begin{tabular}{lc}
\toprule
Graph & Distance Distribution \\
\midrule
MC24/Nauru/HG1234/P24x10 & [1,3,6,9,5] \\
HG36335/P24x9/SG24x12x1 &  [1,3,6,11,3] \\
HG36333/P24x7/SG24x12x2 &  [1,3,5,7,8] \\
\bottomrule
\end{tabular}
\end{center}
\bigskip

\begin{center}
{Number of $k$-cycles ($k=3$ to $8$)}\\
\begin{tabular}{lrrrrrr}
\toprule
Graph & $C_3$ & $C_4$ & $C_5$ & $C_6$ & $C_7$ & $C_8$ \\
\midrule
MC24/Nauru/HG1234/P24x10  & 0 & 0 & 0 & 12 & 0 & 54 \\
HG36335/P24x9/SG24x12x1   & 0 & 0 & 0 & 4 & 24 & 30 \\
HG36333/P24x7/SG24x12x2   & 0 & 6 & 0 & 4 & 0 & 15 \\
\bottomrule
\end{tabular}
\end{center}

%%%%%%%%%%%%%%%%%%%%%%%%%%
\subsection{\large  (3,5)=60}
Reference: \url{https://sagecell.sagemath.org/?q=xagjou}
\bigskip 

\noindent{\bf 3 non-isomorphic optimal Cayley graphs }\\
\texttt{MC60/Foster60/HG36682/P60x25/SG(60,5)->80}. 20 of the graphs are isomorphic to MC60 and the other 60 to HG36684/P60x23.\\
\texttt{HG36698/P60x27/SG(60,6)->40/$Z_{4}\rtimes Z_{(7,13)}  Z_{15}$ }. All 40 graphs are isomorphic to one another and to HG36698/P60x27, but not to MC60 or to HG36684/P60x23.

\begin{center}
\begin{tabular}{lcccccc}
\toprule
Graph & Avg.dist & Girth & Alg.conn. & Dom.num. & Aut.group & Edge-trans. \\
\midrule
MC60/Foster60/HG36682/             &  3.67797 & 9 & 0.697224 & 16 &  360 & True \\
$\;\;\;$  P60x25/SG60x5x1                   &          &   &          &    &      &      \\
HG36684/P60x23/SG60x5x2           & 3.81356 & 6 & 0.585786 & 17 & 120 & False \\
HG36698/P60x27/                   & 3.67797 & 9 & 0.650573 & 15 & 60 & False \\
$\;\;\;$   SG60x6x1/$Z_{4}\rtimes_{(7,13)}  Z_{15}$  &          &   &          &    &      &      \\
\bottomrule
\end{tabular}
\end{center}

\begin{center}
\begin{tabular}{lc}
\toprule
Graph & Distance Distribution \\
\midrule
MC60/Foster60/HG36682/P60x25/SG60x5x1                    & [1, 3, 6, 12, 24, 14]\\
HG36684/P60x23/SG60x5x2                                  &  [1, 3, 6, 11, 18, 21] \\
HG36698/P60x27/SG60x6x1/$Z_{4}\rtimes_{(7,13)}  Z_{15}$  &  [1, 3, 6, 12, 24, 14] \\
\bottomrule
\end{tabular}
\end{center}

\begin{center}
{Number of $k$-cycles ($k=3$ to $10$)}
\begin{tabular}{lrrrrrrrr}
\toprule
Graph & $C_3$ & $C_4$ & $C_5$ & $C_6$ & $C_7$ & $C_8$ & $C_9$ & $C_{10}$ \\
\midrule
MC60/Foster60/HG36682/P60x25/SG60x5x1                     & 0 & 0 & 0 & 0 & 0 & 0 & 60 & 108\\
HG36684/P60x23/SG60x5x2                                   & 0 & 0 & 0 & 10 & 0 & 15 & 20 & 12 \\
HG36698/P60x27/SG60x6x1/$Z_{4}\rtimes_{(7,13)}  Z_{15}$   & 0 & 0 & 0 & 0 & 0 & 0 & 60 & 150 \\
\bottomrule
\end{tabular}
\end{center}

%%%%%%%%%%%%%%%%%%%%%%%%%%%%%%%
\subsection{\large  (3,6)=72}
Reference: \url{https://sagecell.sagemath.org/?q=lbcoqt}
\bigskip 

\noindent{\bf 2 non-isomorphic optimal Cayley graphs }\\
\texttt{ MC72/HG36946/P72x34/SG(72,39) -> 72}. All isomorphic\\
\texttt{ HG36950/P72x18/SG(72,40) -> 144}. All isomorphic to one another, but nor to MC72. 

\begin{center}
\begin{tabular}{lcccccc}
\toprule
Graph & Avg.dist & Girth & Alg.conn. & Dom.num. & Aut.group & Edge-trans. \\
\midrule
MC72/HG36946/P72x34/SG72x39  & 4.00000 & 8 & 0.552836 & 20 & 72 & False \\
HG36950/P72x18/SG72x40       & 4.04225 & 8 & 0.438447 & 20 & 72 & False \\
\bottomrule
\end{tabular}
\end{center}

\begin{center}
\begin{tabular}{lc}
\toprule
Graph & Distance Distribution \\
\midrule
MC72/HG36946/P72x34/SG72x39  & [1, 3, 6, 12, 23, 21, 6] \\
HG36950/P72x18/SG72x40       &  [1, 3, 6, 12, 21, 22, 7] \\
\bottomrule
\end{tabular}
\end{center}

\begin{center}
{Number of $k$-cycles ($k=3$ to $10$)}
\begin{tabular}{lrrrrrrrr}
\toprule
Graph & $C_3$ & $C_4$ & $C_5$ & $C_6$ & $C_7$ & $C_8$ & $C_9$ & $C_{10}$ \\
\midrule
MC72/HG36946/P72x34/SG72x39  & 0 & 0 & 0 & 0 & 0 & 9 & 72 & 36\\
HG36950/P72x18/SG72x40       & 0 & 0 & 0 & 0 & 0 & 27 & 72 & 36 \\
\bottomrule
\end{tabular}
\end{center}

%%%%%%%%%%%%%%%%%%%%%%%%%%%%%%%
\subsection{\large  (3,7)=168}
Reference: \url{https://sagecell.sagemath.org/?q=hzexbs}
\bigskip

\noindent All 504 optimal Cayley graphs found are isomorphic to MC168. \\
\texttt{ MC168/HG39743/P168x51/SG(168,11)->168/SG(168,46)-> 336 }  

%%%%%%%%%%%%%%%%%%%%%%%%%%%%%%%
\subsection{\large  (3,8)=300}
Reference: \url{https://sagecell.sagemath.org/?q=bscore}
\bigskip

\noindent All 1200 optimal Cayley graphs found are isomorphic to one another. \\
\texttt{P300x110/SG(300,24) -> 1200 }  

%%%%%%%%%%%%%%%%%%%%%%%%%%%%%%%
\subsection{\large  (3,9)=506}
Reference: \url{https://sagecell.sagemath.org/?q=nybvnq}
\bigskip

\noindent All 506 optimal Cayley graphs found are isomorphic to one another. \\
\texttt{ P506x57/SG(506,1) -> 506 }  

%%%%%%%%%%%%%%%%%%%%%%%%%%%%%%
\subsection{\large  (3,10)=882}
Reference: \url{https://sagecell.sagemath.org/?q=lexghs}
\bigskip

\noindent P882x134 is the graph 134 from the list for order 882 of all cubic graphs  up to 50000 vertices from  Poto\v{c}nik, Spiga and Verret \cite{Po13}. \\
\texttt{  }

%%%%%%%%%%%%%%%%%%%%%%%%%%%%%%%%%%%%%%%%%%%%%%%%%%%%%%%%%%%%
\section{Degree 4}
%%%%%%%%%%%%%%%%%%%%%%%%%%%%%%%%%%%%%%%%%%%%%%%%%%%%%%%%%%%%
We note that all tetravalent Cayley graphs with order up to 1025 (and  up to isomorphism) have been determined by R.J. Evans and P. Poto\v{c}nik and can be found at \url{https://graphsym.net/}. Here be calculate the main properties of optimal graphs that are not isomorphic to one another for each degree-diameter value and give also where can they be found in  the SmallGroups Library.

%%%%%%%%%%%%%%%%%%%%%%%%%%%%%%%
\subsection{\large  (4,2)=13}
Reference: \url{https://sagecell.sagemath.org/?q=lddkwi}
\bigskip

\noindent All 2 optimal Cayley graphs are isomorphic to MC13.\\
\texttt{MC13/P13x3/SG(13,1)-> 2}

%%%%%%%%%%%%%%%%%%%%%%%%%%%%%%
\subsection{\large  (4,3)=30}
Reference: \url{https://sagecell.sagemath.org/?q=wnwooc}
\bigskip

\noindent All 12 optimal Cayley graphs found are isomorphic to MC30.\\
\texttt{MC30/P30x32/HG52497/SG(30,1)-> 12} \\
Note the existence of two non-Cayley but vertex and edge-transitive graphs.

\begin{center}
\begin{tabular}{lccccc}
\toprule
Graph & Avg.dist & Girth & Alg.conn. & Dom.num & Aut.group \\
\midrule
MC30/P30x32/HG52497/SG30x1  & 2.31034 & 5 & 1.76393 & 6 & 60 \\
HG50708 (no Cayley, vtx/edg trans. ) & 2.31034 & 5 & 2.00000 & 6 & 240 \\
HG50709 (no Cayley. vtx/edg trans.) & 2.31034 & 5 & 2.00000 & 6 & 120 \\
\bottomrule
\end{tabular}
\end{center}

%%%%%%%%%%%%%%%%%%%%%%%%%%%%%%
\subsection{\large  (4,4)=84}
Reference: \url{https://sagecell.sagemath.org/?q=esjsnn}
\bigskip 

\noindent{\bf 2 non-isomorphic optimal Cayley graphs}\\
\texttt{MC84, Loz, P84x170, SG(84,2)->84}.   All are isomorphic. \\    \texttt{P84x182, SG(84,10)->72}.  All 72 are isomorphic to one another but  not to MC84.

\begin{center}
\begin{tabular}{lcccccc}
\toprule
Graph & Avg.dist & Girth & Alg.conn. & Dom.num  & Aut.group & Edge-trans. \\
\midrule
MC84/Loz/P84x170/SG84x2 & 3.14458 & 6 & 1.26795 & 20 &  84 & False \\
P84x182/SG84x10         & 3.13253 & 7 & 1.35817 & 20 &  168 & False \\
\bottomrule
\end{tabular}
\end{center}

\begin{center}
\begin{tabular}{lc}
\toprule
Graph & Distance Distribution \\
\midrule
MC84/Loz/P84x170/SG84x2 & [1, 4, 12, 35, 32] \\
P84x182/SG84x10  &  [1, 4, 12, 36, 31] \\
\bottomrule
\end{tabular}
\end{center}

\begin{center}
{Number of $k$-cycles ($k=3$ to $8$)}\\
\begin{tabular}{lrrrrrr}
\toprule
Graph & $C_3$ & $C_4$ & $C_5$ & $C_6$ & $C_7$ & $C_8$ \\
\midrule
MC84/Loz/P84x170/SG84x2 & 0 & 0 & 0 & 14 & 168 & 756 \\
P84x182/SG84x10  & 0 & 0 & 0 & 0 & 168 & 882 \\
\bottomrule
\end{tabular}
\end{center}

%%%%%%%%%%%%%%%%%%%%%%%%%%%%%%
\subsection{\large  (4,5)=216}
Reference: \url{https://sagecell.sagemath.org/?q=hbfzxi}`
\bigskip 

\noindent{\bf 5 non-isomorphic optimal Cayley graphs }\\
\texttt{MC216/P216x3155/SG(216,87)->432}. All  are isomorphic to one another but not to all other cases.. \\    
\texttt{P216x4719/SG(216,90)->216}. All are isomorphic to one another but not to all other cases.\\
\texttt{P216x2202/SG(216,96)->648}.  All are isomorphic to one another but not to all other cases.\\
\texttt{P216x3422/SG(216,157)->288}.  All are isomorphic to one another but not to all other cases.\\
\texttt{P216x5128/SG(216,158)->432}.All are isomorphic to one another but not to all other cases.

\begin{center}
\begin{tabular}{lccccc}
\toprule
Graph & Avg.dist & Girth & Alg.conn. & Aut.group & Edge-trans. \\
\midrule
MC216 / P216x3155/SG(216,87) & 4.01860 & 8 & 1.00000 & 216 & False \\
P216x4719 / SG(216,90)        & 4.01395 & 7 & 1.00000 & 216 & False \\
P216x2202 / SG(216,96)        & 4.01395 & 7 & 1.00000 & 216 & False \\
P216x3422 / SG(216,157)       & 3.97674 & 8 & 1.00000 & 216 & False \\
P216x5128 / SG(216,158)       & 3.98605 & 8 & 1.00000 & 216 & False \\
\bottomrule
\end{tabular}
\end{center}

\begin{center}
\begin{tabular}{lc}
\toprule
Graph & Distance Distribution \\
\midrule
MC216 / P216x3155/SG(216,87)  &   [1, 4, 12, 36, 87, 76] \\
P216x4719 / SG(216,90)         &   [1, 4, 12, 36, 88, 75] \\
P216x2202 / SG(216,96)         &   [1, 4, 12, 36, 88, 75] \\
P216x3422 / SG(216,157)        &   [1, 4, 12, 36, 96, 67] \\
P216x5128 / SG(216,158)        &   [1, 4, 12, 36, 94, 69] \\
\bottomrule
\end{tabular}
\end{center}
\bigskip

\begin{center}
{Number of $k$-cycles ($k=3$ to $9$)}\\
\begin{tabular}{lrrrrrrr}
\toprule
Graph & $C_3$ & $C_4$ & $C_5$ & $C_6$ & $C_7$ & $C_8$  & $C_9$ \\
\midrule
MC216 / P216x3155/SG(216,87)  & 0 & 0 & 0 & 0 & 0   & 567   & 1368  \\
P216x4719 / SG(216,90)         & 0 & 0 & 0 & 0 & 216 & 216   & 1104  \\
P216x2202 / SG(216,96)         & 0 & 0 & 0 & 0 & 216 & 216   & 1104  \\
P216x3422 / SG(216,157)        & 0 & 0 & 0 & 0 & 0   & 324   & 1584  \\
P216x5128 / SG(216,158)        & 0 & 0 & 0 & 0 & 0   & 378   & 1584  \\
\bottomrule
\end{tabular}
\end{center}

%%%%%%%%%%%%%%%%%%%%%%%%%%%%%%%
\subsection{\large  (4,6)=513}
Reference: \url{https://sagecell.sagemath.org/?q=vmydir}
\bigskip

\noindent All 1026  optimal Cayley graphs found are isomorphic to the Loz graph\cite{Lo09} .\\
\texttt{Loz,SG(513,3)-> 1026}

%%%%%%%%%%%%%%%%%%%%%%%%%%%%%%%%%%%%%%%%%%%%%%%%%%%%%%%%%%%%
\section{Degree 5}
%%%%%%%%%%%%%%%%%%%%%%%%%%%%%%%%%%%%%%%%%%%%%%%%%%%%%%%%%%%%

%%%%%%%%%%%%%%%%%%%%%%%%
\subsection{\large   (5,2)=18}
Reference: \url{https://sagecell.sagemath.org/?q=ayofms}
\bigskip

\noindent All 78 optimal Cayley graphs are isomorphic to MC18.\\
\texttt{ MC18/SG(18,3)->6/SG(18,4)->72}
\bigskip

%%%%%%%%%%%%%%%%%%%%%%%%
\subsection{\large  (5,3)=60}
Reference: \url{https://sagecell.sagemath.org/?q=yusqms}
\bigskip

\noindent All 60 optimal Cayley graphs are isomorphic to MC60.\\
\texttt{ MC60/  SG(60,5) -> 60 }
\bigskip

%%%%%%%%%%%%%%%%%%%%%%%%
\subsection{\large  (5,4)=210}
Reference: \url{https://sagecell.sagemath.org/?q=zegrns}
\bigskip

\noindent All 840 Cayley graphs are isomorphic to the graphs coming from semidirect products found by FC.
\texttt{FC$Z_{6}\rtimes_{(4,9)}Z_{35}$/SG(210,2) -> 840}
\bigskip

Note that Sampels \cite{Sa97} gave a $(5,4) = 210$ graph as the index 14 subgroup 
 in $Z_{42}\rtimes_{39)} Z_{70}$   with generators $[28,14]\-[14,56];[7,18]\-[35,48];[21,14]$.  I have not verified if this graph is isomorphic to the others\\
%\texttt{https://link.springer.com/content/pdf/10.1007/BFb0024505?pdf=chapter+toc}

\begin{center}
\begin{tabular}{lcccccc}
\toprule
Graph & Avg.dist & Girth & Alg.conn. & Dom.num  & Aut.group & Edge-trans. \\
\midrule
FC$Z_{6}\rtimes_{(4,9)}  Z_{35}$ / SG(210,2) & 3.39713 & 6 & 1.58579 & 35 & 210 & False \\
\bottomrule
\end{tabular}
\end{center}

\begin{center}
\begin{tabular}{lc}
\toprule
Graph & Distance Distribution \\
\midrule
FC$Z_{6}\rtimes_{(4,9)}  Z_{35}$ / SG(210,2) & [1, 5, 20, 71, 113]\\
\bottomrule
\end{tabular}
\end{center}
\bigskip

\begin{center}
{Number of $k$-cycles ($k=3$ to $7$)}\\
\begin{tabular}{lrrrrr}
\toprule
Graph & $C_3$ & $C_4$ & $C_5$ & $C_6$ & $C_7$  \\
\midrule
FC$Z_{6}\rtimes_{(4,9)}  Z_{35}$ / SG(210,2) & 0 & 0 & 0 & 350 & 1680  \\
\bottomrule
\end{tabular}
\end{center}
\bigskip

%%%%%%%%%%%%%%%%%%%%%%%%
\subsection{\large  (5,5)=648}
Reference: \url{https://sagecell.sagemath.org/?q=ppwekq}
\bigskip

Marston Conder new best Cayley and  undirected degree-diameter \cite{Co26}.
\bigskip

\begin{center}
\begin{tabular}{lccccc}
\toprule
Graph & Avg.dist & Girth & Alg.conn. & Aut.group & Edge-trans. \\
\midrule
MC648 & 4.23957 & 7 & 1.26795 & 648 & False \\
\bottomrule
\end{tabular}
\end{center}

\begin{center}
\begin{tabular}{lc}
\toprule
Graph & Distance Distribution \\
\midrule
MC648   & [1, 5, 20, 80, 252, 290]\\
\bottomrule
\end{tabular}
\end{center}
\bigskip

\begin{center}
{Number of $k$-cycles ($k=3$ to $8$)}\\
\begin{tabular}{lrrrrrr}
\toprule
Graph & $C_3$ & $C_4$ & $C_5$ & $C_6$ & $C_7$ & $C_8$  \\
\midrule
 MC648   & 0 & 0 & 0 & 0 &  1296 & 3564 \\
\bottomrule
\end{tabular}
\end{center}
\bigskip

%%%%%%%%%%%%%%%%%%%%%%%%%%%%%%%%%%%%%%%%%%%%%%%%%%%%%%%%%%%%
\section{Degree 6}
%%%%%%%%%%%%%%%%%%%%%%%%%%%%%%%%%%%%%%%%%%%%%%%%%%%%%%%%%%%%

%%%%%%%%%%%%%%%%%%%%%%%
\subsection{\large  (6,2)=32}
Reference: \url{https://sagecell.sagemath.org/?q=osibem}
\bigskip

\noindent All 80 optimal Cayley graphs are isomorphic to Wegner/MC32.\\
\texttt{ MC32/Wegner/SG(32,8)->32/SG(32,11)->16/SG(32,50)->32}
\bigskip

C. Delorme (1992), in the context of the general degree-diameter problem,  noticed the existence of this  graph with 32 vertices, degree 6 and diameter 2  known as the Wegner graph and provided its edge list.

Paul Hafner pointed out (July 1998) that the Wegner graph is isomorphic (with automorphism group of order 1920) to the graph described by Lutz Twele \cite{Tw97}
(the graph is also listed in Sampels' PhD thesis). 
Hafner provided also a 'prettier' representation for this group and graph, all generators being involutions, $<f>$ is the centre:\\
{\small $< a,b,c,d,e,f | a^2, b^2, c^2, d^2, e^2, f^2, (af)^2, (bf)^2, (cf)^2, (df)^2, (ef)^2, ababf, cbade, dbaec, ebacd, acacf, ecadb >$}\\
For some more information refer to:  M.J.Dinneen, Group-Theoretic Methods for Designing Networks (http://www.cs.auckland.ac.nz/CDMTCS/researchreports/082nznews.pdf).

Guillermo Pineda informed (May 8, 2008) on the existence of a paper in LNCS v.4123 pp. 853-857 (2006) by Sergey G. Molodsov \cite{Mo06} 
reporting a computer search which  generated all largest diameter 2 degree 6 graphs, resulting in 6 non isomorphic graphs with order 32 (one of them is the non Cayley {$K_4\st X_8$}  and another the Wegner graph), and thus all these graphs are optimal for the general degree.diameter problem, but onlythe Wegner graph is optimalin our context.\\
\texttt{\small https://web.mat.upc.edu/francesc.comellas/delta-d/desc\_g/desc\_g2.html\#62}

%%%%%%%%%%%%%%%%%%%%%%%
\subsection{\large  (6,3)=108}
Reference: \url{https://sagecell.sagemath.org/?q=ioyixy}
\bigskip

\noindent All 324 optimal Cayley graphs found are isomorphic to Wohlmuth / MC108. \\
\texttt{ MC108 / Wohlmuth \cite{Wo96}/  Sampels \cite{Sa97} /  SG(108,21)->324 }  
\bigskip

\bigskip

%%%%%%%%%%%%%%%%%%%%%%%
\subsection{\large  (6,4)=384}
Reference: \url{https://sagecell.sagemath.org/?q=mlsjpl}
\bigskip

Marston Conder new  best  Cayley degree-diameter graph \cite{Co26}.
\bigskip 

\begin{center}
\begin{tabular}{lccccc}
\toprule
Graph & Avg.dist & Girth & Alg.conn. & Aut.group & Edge-trans. \\
\midrule
MC384 & 3.44125  & 6 & 2.22800 & 82944 & True \\
\bottomrule
\end{tabular}
\end{center}

\begin{center}
\begin{tabular}{lc}
\toprule
Graph & Distance Distribution \\
\midrule
MC384 & [1, 6, 30, 136, 211]\\
\bottomrule
\end{tabular}
\end{center}
\bigskip

\begin{center}
{Number of $k$-cycles ($k=3$ to $8$)}\\
\begin{tabular}{lrrrrrr}
\toprule
Graph & $C_3$ & $C_4$ & $C_5$ & $C_6$ & $C_7$  & $C_8$  \\
\midrule
MC384  & 0 & 0 & 0 &  1152&  0 & 46368   \\
\bottomrule
\end{tabular}
\end{center}
\bigskip

%%%%%%%%%%%%%%%%%%%%%%%%%%%%%%%%%%%%%%%%%%%%%%%%%%%%%%%%%%%%
\section{Degree 7}
%%%%%%%%%%%%%%%%%%%%%%%%%%%%%%%%%%%%%%%%%%%%%%%%%%%%%%%%%%%%

%%%%%%%%%%%%%%%%%%%%%%%
\subsection{\large  (7,2)=36}
Reference: \url{https://sagecell.sagemath.org/?q=axwlkt}
\bigskip

\noindent All 144 optimal Cayley graphs found are isomorphic to MC36. \\
\texttt{ MC36/ SG(36,9) -> 144 }  
\bigskip

%%%%%%%%%%%%%%%%%%%%%%%
\subsection{\large  (7,3)=168}
Reference: \url{https://sagecell.sagemath.org/?q=fdkifl}
\bigskip 

\noindent{\bf 5 non-isomorphic optimal Cayley graphs }\\
\texttt{ SG(168,43)->588}. MC168/SG168x43x3  are isomorphic \\     
SG168x43x1/SG168x43x2 SG168x43x4/ SG168x43x5  are non-isomorphic to one another and to  MC168 

Note that these five graphs, together with a non-Cayley graph due to G. Exoo \cite{Com25}, also represent the current best-known examples for the undirected degree-diameter problem.

\bigskip

\begin{center}
\begin{tabular}{lccccc}
\toprule
Graph & Avg.dist & Girth & Alg.conn. & Aut.group & Edge-trans. \\
\midrule
 MC168 / SG168x43x3  & 2.67665 & 4 & 3.00000 & 336 & False \\
SG168x43x1 & 2.67665 & 4 & 3.00000 & 336 & False \\
SG168x43x2 & 2.67665 & 4 & 3.00000 & 336 & False \\
SG168x43x4 & 2.67665 & 3 & 3.29453 & 168 & False \\
SG168x43x5 & 2.66467 & 5 & 3.00000 & 168 & False \\
\bottomrule
\end{tabular}
\end{center}

\begin{center}
\begin{tabular}{lc}
\toprule
Graph & Distance Distribution \\
\midrule
MC168/SG168x43x1 /SG168x43x4 / SG168x43x5  & [1, 7, 40, 120] \\
SG168x43x2 & [1, 7, 42, 118] \\
\bottomrule
\end{tabular}
\end{center}
\bigskip

\begin{center}
{Number of $k$-cycles ($k=3$ to $6$)}\\
\begin{tabular}{lrrrrrr}
\toprule
Graph & $C_3$ & $C_4$ & $C_5$ & $C_6$  \\
\midrule
MC168/SG168x43x3      &  0 & 84 &1176 & 2604   \\
SG168x43x1 &  0 & 84 & 840 & 3108   \\
SG168x43x2 &  0 & 84 &1176 & 2884  \\
SG168x43x4 & 56 &  0 & 504 & 3528  \\
SG168x43x5 &  0 &  0 & 336 & 4536  \\
\bottomrule
\end{tabular}
\end{center}

%%%%%%%%%%%%%%%%%%%%%%%%%%%%%%%%%%%%%%%%%%%%%%%%%%%%%%%%%%%%
\section{Degree 8}
%%%%%%%%%%%%%%%%%%%%%%%%%%%%%%%%%%%%%%%%%%%%%%%%%%%%%%%%%%%%

%%%%%%%%%%%%%%%%%%%%%%%
\subsection{\large  (8,2)=48}
Reference: \url{https://sagecell.sagemath.org/?q=upksyu}
\bigskip 

\noindent{\bf 2 non-isomorphic optimal Cayley graphs }\\
\texttt{ SG(48,5) -> 48, SG(48,10) -> 48,  SG(48,15) -> 16,  SG(48,16)  -> 96, SG(48,17) -> 48,\\
SG(48,25) -> 64, SG(48,39) -> 64,  SG(48,47) -> 48}\\
432 graphs, 280 are isomorphic to one another and to MC48 and the other 152 are also isomorphic to one another but non-isomorphic to MC48

\bigskip 

\begin{center}
\begin{tabular}{lcccccc}
\toprule
Graph & Avg.dist & Girth & Alg.conn. & Dom.num  & Aut.group & Edge-trans. \\
\midrule
MC48/SG48x5x1 & 1.82979 & 4 & 4.00000 &  8  &  288 & False \\
SG48x5x13     & 1.82979 & 4 & 4.00000 &  8  &   96 & False \\
\bottomrule
\end{tabular}
\end{center}

\begin{center}
\begin{tabular}{lc}
\toprule
Graph & Distance Distribution \\
\midrule
MC48/SG48x5x1   & [1, 8, 39] \\
SG48x5x13  & [1, 8, 39] \\
\bottomrule
\end{tabular}
\end{center}
\bigskip

\begin{center}
{Number of $k$-cycles ($k=3$ to $6$)}\\
\begin{tabular}{lrrrrrr}
\toprule
Graph & $C_3$ & $C_4$ & $C_5$ & $C_6$  \\
\midrule
MC48/SG48x5x1       &  0 & 228 & 2352 & 10760   \\
SG48x5x13  &  0 & 228 & 2352 & 10712  \\
\bottomrule
\end{tabular}
\end{center}

%%%%%%%%%%%%%%%%%%%%%%%%%%%%%%%%%%%%%%%%%%%%%%%%%%%%%%%%%%%%
\section{Degree 9}
%%%%%%%%%%%%%%%%%%%%%%%%%%%%%%%%%%%%%%%%%%%%%%%%%%%%%%%%%%%%

\subsection{\large  (9,2)=60}
Reference: \url{https://sagecell.sagemath.org/?q=lbfzdc}
\bigskip 

\noindent{\bf 2 non-isomorphic optimal Cayley graphs }\\
\texttt{ SG(60,9) -> 60}   Twelve are isomorphic to MC60 and the others 48 are isomorphic to one another  but not to MC60. 

\begin{center}
\begin{tabular}{lcccccc}
\toprule
Graph & Avg.dist & Girth & Alg.conn. & Dom.num & Aut.group & Edge-trans. \\
\midrule
MC60 /SG60x9x11 & 1.84746 & 3 & 5.00000 & 9 & 480 & False \\
SG60x9x1        & 1.84746 & 4 & 5.34582 & 9 & 120 & False \\
\bottomrule
\end{tabular}
\end{center}

\begin{center}
\begin{tabular}{lc}
\toprule
Graph & Distance Distribution \\
\midrule
MC60 / SG60x9x1   & [1, 9, 50] \\
\bottomrule
\end{tabular}
\end{center}
\bigskip

\begin{center}
{Number of $k$-cycles ($k=3$ to $6$)}\\
\begin{tabular}{lrrrrrr}
\toprule
Graph & $C_3$ & $C_4$ & $C_5$ & $C_6$  \\
\midrule
MC60 / SG60x9x11      & 40 & 330 & 4080 & 23560   \\
SG60x9x1   & 0  & 330 & 4560 & 23800  \\
\bottomrule
\end{tabular}
\end{center}

%%%%%%%%%%%%%%%%%%%%%%%
\subsection{\large  (9,3)=300}
Reference: \url{https://sagecell.sagemath.org/?q=ssbcbz}
\bigskip

Two non-isomorphic graphs to one another that improve the best current value of 294 in the table  (not listed as optimal).  Obtained from SG(300,22) 

\begin{center}
\begin{tabular}{lccccc}
\toprule
Graph      & Avg.dist & Girth & Alg.conn. &  Aut.group & Edge-trans.\\
\midrule
SG300x22x1   & 2.69900  & 5  &  4.27235   &  300 & False\\ 
SG300x22x2   & 2.69900  & 5  &  4.64685   &  300 & False\\ 
\bottomrule
\end{tabular}
\end{center}

\begin{table}[ht]
\centering
\begin{minipage}[t]{0.42\textwidth}
\centering
{Distance Distribution}\\
\smallskip
\begin{tabular}{lc}
\toprule
Graph & Distance Distribution \\
\midrule
SG300x22x1   &  [1, 9, 72, 218] \\
SG300x22x2   &  [1, 9, 72, 218] \\
\bottomrule
\end{tabular}
\end{minipage}
\hspace{0.02\textwidth} % smaller gap
\begin{minipage}[t]{0.40\textwidth}
\centering
{Number of $k$-cycles ($k=3$ to $6$)}\\
\begin{tabular}{lrrrr}
\toprule
Graph & $C_3$ & $C_4$ & $C_5$ & $C_6$     \\
\midrule
SG300x22x1  &  0  &  0  &  3960  &  20500 \\
SG300x22x2  &  0  &  0  &  3360  &  19750 \\
\bottomrule
\end{tabular}
\end{minipage}
\end{table}

%%%%%%%%%%%%%%%%%%%%%%%%%%%%%%%%%%%%%%%%%%%%%%%%%%%%%%%%%%%%
\section{Degree 10}
%%%%%%%%%%%%%%%%%%%%%%%%%%%%%%%%%%%%%%%%%%%%%%%%%%%%%%%%%%%%

\subsection{\large  (10,2)=72}
Reference: \url{https://sagecell.sagemath.org/?q=hfvnnx}

\bigskip 

\noindent{\bf 8 non-isomorphic optimal Cayley graphs }\\
\texttt{ MC72 / SG(72,27)-> 48}. All are isomorphic to MC72\\  
\texttt{  SG(72,28)-> 192}. Two of these graphs are non-isomorphic to one another
and to MC96.\\
\texttt{  SG(72,30)-> 240}. Five of these graphs are non-isomorphic to one another and to the 
other cases.

\bigskip 
\begin{center}
\begin{tabular}{lcccccc}
\toprule
Graph & Avg.dist & Girth & Alg.conn. & Dom.num  & Aut.group & Edge-trans.\\
\midrule
MC72 / SG72x27 & 1.85915 & 4 & 6.00000 & 9 & 72 & False\\
SG72x28x1      & 1.85915 & 4 & 6.00000 & 10 & 72 & False \\
SG72x28x2      & 1.85915 & 4 & 6.00000 & 10 & 72 & False\\
SG72x30x1      & 1.85915 & 3 & 5.64110 & 10 & 72 & False\\
SG72x30x5      & 1.85915 & 3 & 5.64110 & 9 & 72 & False\\
SG72x30x9      & 1.85915 & 3 & 6.00000 & 9 & 72 & False\\
SG72x30x13     & 1.85915 & 3 & 5.64110 & 9 & 72 & False\\
SG72x30x17     & 1.85915 & 3 & 6.00000 & 9 & 72 & False\\
\bottomrule
\end{tabular}
\end{center}

\begin{center}
\begin{tabular}{lc}
\toprule
Graph & Distance Distribution \\
\midrule
MC72/ SG72x27/SG72x28/SG72x30  &  [1, 10, 61]\\
\bottomrule
\end{tabular}
\end{center}
\bigskip

\begin{center}
{Number of $k$-cycles ($k=3$ to $7$)}\\
\begin{tabular}{lrrrrr}
\toprule
Graph & $C_3$ & $C_4$ & $C_5$ & $C_6$ & $C_7$  \\
\midrule
MC72 / SG72x27 &   0 & 594 & 7776 & 47136 & 327528 \\
SG72x28x1      &   0 & 594 & 7776 & 47184 & 327240 \\
SG72x28x2      &   0 & 594 & 7776 & 47112 & 327456 \\
SG72x30x1      & 144 & 378 & 5904 & 47712 & 319680 \\
SG72x30x5      & 144 & 378 & 5904 & 47640 & 318672 \\
SG72x30x9      & 144 & 342 & 5904 & 47988 & 321696 \\
SG72x30x13     & 144 & 378 & 5904 & 47568 & 318672 \\
SG72x30x17     & 144 & 342 & 5904 & 47892 & 321696 \\
\bottomrule
\end{tabular}
\end{center}
\bigskip

%%%%%%%%%%%%%%%%%%%%%%%
\subsection{\large  (10,3)=420}
Reference: \url{https://sagecell.sagemath.org/?q=lhigxq}
\bigskip

Two non-isomorphic graphs not listed as optimal but that improve the current value in the table, obtained from SG(420,18) and semidirect product \texttt{$Z_{12}\rtimes_{2}Z_{35}$}.

\bigskip 
\begin{center}
\begin{tabular}{lccccc}
\toprule
Graph & Avg.dist & Girth & Alg.conn. &  Aut.group & Edge-trans.\\
\midrule
SG420x18                      & 2.73986  & 4  &  4.00000   &  840 & False \\ 
$Z_{12}\rtimes_{2}Z_{35}$     & 2.73747  & 5  &  4.00000   &  420 & False \\ 
\bottomrule
\end{tabular}
\end{center}

\begin{table}[ht]
\centering
\begin{minipage}[t]{0.42\textwidth}
\centering
{Distance Distribution}\\
\begin{tabular}{lc}
\toprule
Graph & Distance Distribution \\
\midrule
SG420x18                     & [1, 10, 89, 320] \\
$Z_{12}\rtimes_{2}Z_{35}$    & [1, 10, 90, 319] \\
\bottomrule
\end{tabular}
\end{minipage}
\hspace{0.02\textwidth} % smaller gap
\begin{minipage}[t]{0.40\textwidth}
\centering
{Number of $k$-cycles ($k=3$ to $5$)}\\
\begin{tabular}{lrrr}
\toprule
Graph    & $C_3$ & $C_4$ & $C_5$  \\
\midrule
SG420x18                    &   0  &  105  &  4200  \\
$Z_{12}\rtimes_{2}Z_{35}$   &  0  &    0  &  5460  \\
\bottomrule
\end{tabular}
\end{minipage}
\end{table}

%%%%%%%%%%%%%%%%%%%%%%%%%%%%%%%%%%%%%%%%%%%%%%%%%%%%%%%%%%%%
\section{Degree 11}
%%%%%%%%%%%%%%%%%%%%%%%%%%%%%%%%%%%%%%%%%%%%%%%%%%%%%%%%%%%%

\subsection{\large  (11,2)=84}
Reference: \url{https://sagecell.sagemath.org/?q=szbabj}
\bigskip 

\noindent{\bf 4 non-isomorphic optimal Cayley graphs }\\
\texttt{ SG(84,10)-> 252,}  Four are non-isomorphic to one another and one of them is isomorphic to MC84

\begin{center}
\begin{tabular}{lccccc}
\toprule
Graph & Avg.dist & Girth & Alg.conn. & Dom.num & Aut.group \\
\midrule
MC84 / SG84x10x1 & 1.86747 & 3 & 6.00000 & 11 & 336 \\
SG84x10x3  & 1.86747 & 4 & 6.00000 & 11 & 168 \\
SG84x10x5 & 1.86747 & 4 & 6.00000 & 11 & 168 \\
SG84x10x7 & 1.86747 & 4 & 6.00000 & 11 & 168 \\
\bottomrule
\end{tabular}
\end{center}

\begin{center}
\begin{tabular}{lc}
\toprule
Graph & Distance Distribution \\
\midrule
MC84/ SG84x10x3,5,7 & [1, 11, 72]\\
\bottomrule
\end{tabular}
\end{center}

\begin{center}
{Number of $k$-cycles ($k=3$ to $7$)}\\
\begin{tabular}{lrrrrr}
\toprule
Graph & $C_3$ & $C_4$ & $C_5$ & $C_6$ & $C_7$  \\
\midrule
MC84 / SG84x10x1  &  56 & 1008  & 11760  & 85736  & 670908 \\
SG84x10x3       &   0 & 1008  & 12768  & 86072  & 681744 \\
SG84x10x5       &   0 & 1008  & 12768  & 86072  & 681408 \\
SG84x10x7       &   0 & 1008  & 12768  & 86072  & 681576 \\
\bottomrule
\end{tabular}
\end{center}

%%%%%%%%%%%%%%%%%%%%%%%
\subsection{\large  (11,3)=504}
Reference: \url{https://sagecell.sagemath.org/?q=iuqfeb}
\bigskip

Three non-isomorphic graphs not listed as optimal but that improve the current value in the table, obtained from SG(504,52) and SG(504,196),

\bigskip 
\begin{center}
\begin{tabular}{lccccc}
\toprule
Graph & Avg.dist & Girth & Alg.conn. &  Aut.group & Edge-trans.\\
\midrule
SG504x52   & 2.73757  & 5  &  5.83916   &  504 & False \\ 
SG504x52b   & 2.73757  & 5  &  5.83916   &  504 & False \\ 
SG504x196  & 2.74950  & 3  &  5.57963   &  504 & False \\ 
\bottomrule
\end{tabular}
\end{center}

\begin{table}[ht]
\centering
\begin{minipage}[t]{0.42\textwidth}
\centering
{Distance Distribution}\\
\begin{tabular}{lc}
\toprule
Graph & Distance Distribution \\
\midrule
SG504x52    & [1, 11, 110, 382] \\
SG504x52b    & [1, 11, 110, 382] \\
SG504x196   & [1, 11, 104, 388] \\
\bottomrule
\end{tabular}
\end{minipage}
\hspace{0.02\textwidth} % smaller gap
\begin{minipage}[t]{0.40\textwidth}
\centering
{Number of $k$-cycles ($k=3$ to $5$)}\\
\begin{tabular}{lrrr}
\toprule
Graph    & $C_3$ & $C_4$ & $C_5$  \\
\midrule
SG504x52  &  0  &  0  &  10584  \\
SG504x52b  &   0  &  0  &  9072 \\
SG504x196  & 168  &  504  &  7056  \\
\bottomrule
\end{tabular}
\end{minipage}
\end{table}

%%%%%%%%%%%%%%%%%%%%%%%%%%%%%%%%%%%%%%%%%%%%%%%%%%%%%%%%%%%%
\section{Degree 12}
%%%%%%%%%%%%%%%%%%%%%%%%%%%%%%%%%%%%%%%%%%%%%%%%%%%%%%%%%%%%

\subsection{\large  (12,2)=96}
Reference: \url{https://sagecell.sagemath.org/?q=zrremt}

\noindent{\bf 64 non-isomorphic optimal Cayley graphs }\\
\texttt{SG(96,3)->576, SG(96,4)->96, SG(96,5)-> 288, SG(96,31)->3072, SG(96,34)->2304, SG(96,41)->384,\\
SG(96,43)->1536, SG(96,49)-128, SG(96,51)->1536, SG(96,62)->1792, SG(96,73)->384, SG(96,84)->768,\\ SG(96,86)->1536,
SG(96,171)->512, SG(96,173)->768, SG(96,182)->256, SG(96,185)->384,\\ 
SG(96,190)->192}, SG(96,191)->96. Of all these 15608 graphs, 64  are non-isomorphic to one another and one of them is MC96. 
\bigskip 
\bigskip\bigskip

All graphs have average distance 1.87368, distance distribution [1,12,83] and are not edge-transitive.

\begin{center}\footnotesize
\begin{tabular}{lcccc}
\toprule
Graph           & Girth & Alg.conn. & Dom.number & Aut.group \\
\midrule
MC96 / SG96x3x65 &   3    &   6.00000   &    11   &    192 \\ 
SG96x3x1         &   3    &   6.00000   &    11   &    192  \\  
SG96x3x385       &   3    &   6.00000   &    11   &    192  \\  
SG96x4x1         &   4    &   6.00000   &    12   &    96  \\  
SG96x5x1         &   4    &   6.70850   &    12   &    96  \\  
SG96x5x2         &   4    &   6.70850   &    12   &    96  \\  
SG96x5x3         &   4    &   6.70850   &    11   &    96  \\  
SG96x31x1        &   3    &   7.13630   &    12   &    96  \\  
SG96x31x3        &   3    &   7.13630   &    12   &    96  \\  
SG96x31x13       &   3    &   7.13630   &    12   &    96  \\  
SG96x31x15       &   3    &   7.13630   &    12   &    96  \\  
SG96x31x1537     &   3    &   7.00000   &    12   &    96  \\  
SG96x31x1539     &   3    &   7.00000   &    12   &    96  \\  
SG96x31x1549     &   3    &   7.00000   &    12   &    96  \\  
SG96x31x1551     &   3    &   7.00000   &    12   &    96  \\  
SG96x34x1        &   3    &   5.70850   &    12   &    96  \\  
SG96x34x3        &   3    &   5.70850   &    11   &    96  \\  
SG96x34x25       &   3    &   5.70850   &    12   &    96  \\  
SG96x34x27       &   3    &   5.70850   &    12   &    96  \\  
SG96x34x1537     &   4    &   6.00000   &    11   &    96  \\  
SG96x34x1549     &   4    &   6.00000   &    11   &    96  \\  
SG96x41x1        &   3    &   7.00000   &    12   &    384  \\  
SG96x41x193      &   3    &   7.00000   &    11   &    384  \\  
SG96x43x1        &   3    &   7.00000   &    12   &    96  \\  
SG96x43x3        &   3    &   7.00000   &    12   &    96  \\  
SG96x49x1        &   3    &   7.00000   &    12   &    384  \\  
SG96x49x65       &   3    &   7.00000   &    11   &    384  \\  
SG96x51x1        &   3    &   7.53590   &    12   &    192  \\  
SG96x51x129      &   3    &   7.26795   &    12   &    96  \\  
SG96x51x130      &   3    &   7.26795   &    12   &    96  \\  
SG96x51x133      &   3    &   7.00000   &    12   &    96  \\  
SG96x51x135      &   3    &   7.00000   &    12   &    96  \\  
SG96x51x1281     &   3    &   7.53590   &    11   &    192  \\  
SG96x62x1        &   3    &   6.00000   &    12   &    96  \\  
SG96x62x2        &   3    &   6.00000   &    12   &    96  \\  
SG96x62x65       &   3    &   6.00000   &    12   &    96  \\  
SG96x62x66       &   3    &   6.00000   &    12   &    96  \\  
SG96x62x129      &   3    &   6.00000   &    12   &    96  \\  
SG96x62x130      &   3    &   6.00000   &    12   &    96  \\  
SG96x62x145      &   3    &   6.00000   &    12   &    96  \\  
SG96x62x146      &   3    &   6.00000   &    11   &    96  \\  
SG96x62x193      &   3    &   6.00000   &    12   &    96  \\  
SG96x62x194      &   3    &   6.00000   &    11   &    96  \\  
SG96x62x321      &   3    &   6.00000   &    12   &    96  \\  
SG96x62x322      &   3    &   6.00000   &    12   &    96  \\  
SG96x62x329      &   3    &   6.00000   &    12   &    96  \\  
SG96x62x330      &   3    &   6.00000   &    12   &    96  \\  
SG96x73x1        &   3    &   6.00000   &    12   &    768  \\  
SG96x73x5        &   4    &   6.00000   &    11   &    768  \\  
SG96x73x9        &   4    &   6.00000   &    12   &    384  \\  
SG96x73x17       &   4    &   6.00000   &    12   &    192  \\  
SG96x73x49       &   4    &   6.00000   &    12   &    192  \\  
SG96x84x1        &   3    &   7.00000   &    11   &    384  \\  
SG96x86x1        &   3    &   6.00000   &    12   &    96  \\  
SG96x86x13       &   3    &   6.00000   &    12   &    96  \\  
SG96x86x769      &   3    &   6.00000   &    11   &    96  \\  
SG96x86x781      &   3    &   6.00000   &    11   &    96  \\  
SG96x171x1       &   3    &   6.00000   &    12   &    192  \\  
SG96x173x1       &   3    &   6.00000   &    12   &    192  \\  
SG96x173x385     &   3    &   6.00000   &    11   &    96  \\  
SG96x182x1       &   4    &   6.00000   &    11   &    288  \\  
SG96x182x2       &   4    &   6.00000   &    12   &    96  \\  
SG96x185x1       &   3    &   6.00000   &    12   &    768  \\  
SG96x185x17      &   4    &   6.00000   &    11   &    768  \\  
SG96x185x33      &   4    &   6.00000   &    12   &    384  \\  
SG96x190x1       &   3    &   6.00000   &    12   &    96  \\  
SG96x190x2       &   3    &   6.00000   &    12   &    96  \\  
SG96x191x1       &   3    &   6.00000   &    10   &    960  \\  
\bottomrule
\end{tabular}
\end{center}

\bigskip
\bigskip\bigskip\bigskip\bigskip

\begin{center}\footnotesize
\begin{tabular}{lrrrrrr}
{\small Number of $k$-cycles ($k=3$ to $8$)}\\
\toprule
Graph         & $C_3$ & $C_4$ & $C_5$ & $C_6$ & $C_7$  & $C_8$ \\
\midrule
MC96 / SG96x3x65      & 256 & 1152 & 16224 & 149280  & 1268832 & 11299932 \\
SG96x3x1              &  96 & 1680 & 17952 & 148336  & 1296384 & 11741916 \\
SG96x3x385            & 192 & 1344 & 17184 & 147616  & 1269696 & 11429052 \\
SG96x4x1              &   0 & 1656 & 20064 & 151168  & 1320096 & 12190020 \\
SG96x5x1              &   0 & 1464 & 20064 & 154048  & 1339296 & 12502308 \\
SG96x5x2              &   0 & 1656 & 19296 & 150400  & 1324608 & 12259476 \\
SG96x5x3              &   0 & 1464 & 20064 & 153952  & 1338048 & 12501252 \\
SG96x31x1             & 192 & 1128 & 16896 & 154864  & 1327680 & 11936064 \\
SG96x31x3             & 192 & 1128 & 16896 & 154864  & 1327392 & 11938656 \\
SG96x31x13            & 192 & 1128 & 16896 & 154864  & 1327392 & 11937792 \\
SG96x31x15            & 192 & 1128 & 16896 & 154864  & 1327584 & 11938272 \\
SG96x31x1537          & 288 & 1128 & 15360 & 147344  & 1263360 & 11193360 \\
SG96x31x1539          & 288 & 1128 & 15360 & 147344  & 1263456 & 11194224 \\
SG96x31x1549          & 288 & 1128 & 15360 & 147248  & 1263456 & 11197584 \\
SG96x31x1551          & 288 & 1128 & 15360 & 147248  & 1263360 & 11197200 \\
SG96x34x1             & 288 & 1008 & 16416 & 151936  & 1263552 & 11121468 \\
SG96x34x3             & 288 & 1104 & 16416 & 150016  & 1258848 & 11126676 \\
SG96x34x25            & 288 & 1104 & 16416 & 150016  & 1257120 & 11106420 \\
SG96x34x27            & 288 & 1104 & 16416 & 149728  & 1256064 & 11114532 \\
SG96x34x1537          &   0 & 1608 & 20256 & 152368  & 1323552 & 12243684 \\
SG96x34x1549          &   0 & 1608 & 20256 & 152272  & 1323840 & 12245172 \\
SG96x41x1             & 192 & 1272 & 17472 & 149296  & 1282368 & 11670456 \\
SG96x41x193           & 192 &  984 & 17472 & 155968  & 1319328 & 11969544 \\
SG96x43x1             &  96 & 1272 & 18720 & 154352  & 1331424 & 12232776 \\
SG96x43x3             &  96 & 1272 & 18720 & 154352  & 1331136 & 12230952 \\
SG96x51x1             & 192 & 1272 & 17088 & 149680  & 1291776 & 11708184 \\
SG96x51x129           &  96 & 1176 & 18816 & 156032  & 1338432 & 12353484 \\
SG96x51x130           &  96 & 1176 & 18816 & 156032  & 1339104 & 12355404 \\
SG96x51x133           &  96 & 1272 & 18720 & 154352  & 1331136 & 12233448 \\
SG96x51x135           &  96 & 1272 & 18720 & 154352  & 1331328 & 12231816 \\
SG96x51x1281          & 192 &  984 & 17280 & 155584  & 1324704 & 12004776 \\
SG96x62x1             & 288 & 1008 & 16032 & 151552  & 1279392 & 11319024 \\
SG96x62x2             & 288 & 1008 & 16032 & 151744  & 1280064 & 11315184 \\
SG96x62x65            & 288 & 1008 & 16032 & 151552  & 1279392 & 11319360 \\
SG96x62x66            & 288 & 1008 & 16032 & 151744  & 1280064 & 11315520 \\
SG96x62x129           & 288 & 1032 & 16128 & 151648  & 1273728 & 11251836 \\
SG96x62x130           & 288 & 1032 & 16224 & 150592  & 1271328 & 11272812 \\
SG96x62x145           & 288 & 1032 & 16128 & 151360  & 1274304 & 11270940 \\
SG96x62x146           & 288 & 1032 & 16224 & 150688  & 1271040 & 11264076 \\
SG96x62x193           & 288 &  984 & 16224 & 152176  & 1276992 & 11302380 \\
SG96x62x194           & 288 &  984 & 16128 & 152944  & 1278624 & 11290572 \\
SG96x62x321           & 288 & 1032 & 16032 & 151168  & 1272672 & 11252472 \\
SG96x62x322           & 288 &  936 & 16224 & 153280  & 1283328 & 11349336 \\
SG96x62x329           & 288 & 1032 & 16032 & 151456  & 1272768 & 11245800 \\
SG96x62x330           & 288 &  936 & 16224 & 153376  & 1283520 & 11347752 \\
SG96x73x1             &  64 & 1680 & 18624 & 149632  & 1301472 & 11791032 \\
SG96x73x5             &   0 & 1680 & 19968 & 150272  & 1318080 & 12151224 \\
SG96x73x9             &   0 & 1680 & 19968 & 150016  & 1316928 & 12155640 \\
SG96x73x17            &   0 & 1680 & 19968 & 150016  & 1317312 & 12157944 \\
SG96x73x49            &   0 & 1680 & 19968 & 150016  & 1317312 & 12158328 \\
SG96x86x1             & 288 &  960 & 16224 & 152464  & 1279872 & 11336688 \\
SG96x86x13            & 288 &  960 & 16032 & 152896  & 1285824 & 11365200 \\
SG96x86x769           & 288 &  960 & 16032 & 152896  & 1285824 & 11365344 \\
SG96x86x781           & 288 &  864 & 16224 & 154480  & 1282944 & 11327880 \\
SG96x171x1            & 288 &  960 & 16032 & 153376  & 1285440 & 11348592 \\
SG96x173x385          & 288 &  960 & 16032 & 153280  & 1286016 & 11353248 \\
SG96x182x1            &   0 & 1560 & 20448 & 153568  & 1326048 & 12304788 \\
SG96x182x2            &   0 & 1752 & 19680 & 148960  & 1315872 & 12084036 \\
SG96x185x1            &  64 & 1584 & 19008 & 151552  & 1308960 & 11905440 \\
SG96x185x17           &   0 & 1584 & 20352 & 152192  & 1324032 & 12263328 \\
SG96x185x33           &   0 & 1584 & 20352 & 151936  & 1323648 & 12273504 \\
SG96x190x1            & 128 & 1272 & 18336 & 152688  & 1310208 & 11983188 \\
SG96x190x2            & 128 & 1272 & 18336 & 152784  & 1309920 & 11978532 \\
SG96x191x1            & 288 & 1392 & 16416 & 138352  & 1174560 & 10521660 \\
\bottomrule
\end{tabular}
\end{center}
\bigskip

%%%%%%%%%%%%%%%%%%%%%%%
\subsection{\large  (12,3)=620}
Reference: \url{https://sagecell.sagemath.org/?q=tiicga}
\bigskip

Ten non-isomorphic graphs not listed as optimal but that improve the current value in the table.
Two are obtained from SG(620) and the other are 
semidirect product \texttt{$Z_{20}\rtimes_{(2,4,8,15,16,23,27,29)}Z_{31}$}.
They have these properties:

\bigskip 
\begin{center}
\begin{tabular}{lccccc}
\toprule
Graph & Avg.dist & Girth & Alg.conn. &  Aut.group & Edge-trans.\\
\midrule
SG620x1   & 2.74798  & 5  &  6.75519   &  620 & False\\
SG620x2   & 2.75444  & 4  &  6.35924   &  620 & False\\ 
$Z_{20}\rtimes_{2}Z_{31}$  & 2.75444  & 4  &  6.45525   &  620 & False \\ 
$Z_{20}\rtimes_{4}Z_{31}$   & 2.75444  & 4  &  6.08599   &  620 & False \\ 
$Z_{20}\rtimes_{8}Z_{31}$  & 2.74798  & 5  &  6.60233   &  620 & False \\ 
$Z_{20}\rtimes_{15}Z_{31}$   & 2.75444  & 4  &  6.52786   &  620 & False \\ 
$Z_{20}\rtimes_{16}Z_{31}$  & 2.75444  & 4  &  6.52909   &  620 & False \\ 
$Z_{20}\rtimes_{23}Z_{31}$  & 2.74798  & 5  &  6.50118   &  620 & False\\ 
$Z_{20}\rtimes_{27}Z_{31}$  & 2.74798  & 5  &  6.60943   &  620 & False\\ 
$Z_{20}\rtimes_{29}Z_{31}$   & 2.74798  & 5  &  6.50980   &  620 & False\\ 

\bottomrule
\end{tabular}
\end{center}

\begin{table}[ht]
\centering
\begin{minipage}[t]{0.42\textwidth}
\centering
{Distance Distribution}\\
\begin{tabular}{lc}
\toprule
Graph & Distance Distribution \\
\midrule
SG620x1   &   [1, 12, 132, 475]  \\
SG620x2   &   [1, 12, 128, 479] \\
$Z_{20}\rtimes_{2}Z_{31}$  & [1, 12, 128, 479] \\
$Z_{20}\rtimes_{4}Z_{31}$  & [1, 12, 128, 479] \\
$Z_{20}\rtimes_{8}Z_{31}$  & [1, 12, 132, 475] \\
$Z_{20}\rtimes_{15}Z_{31}$  & [1, 12, 128, 479]\\
$Z_{20}\rtimes_{16}Z_{31}$  & [1, 12, 128, 479] \\
$Z_{20}\rtimes_{23}Z_{31}$  & [1, 12, 132, 475]\\
$Z_{20}\rtimes_{27}Z_{31}$  & [1, 12, 132, 475]\\
$Z_{20}\rtimes_{29}Z_{31}$  & [1, 12, 132, 475]\\
\bottomrule
\end{tabular}
\end{minipage}
\hspace{0.02\textwidth} % smaller gap
\begin{minipage}[t]{0.40\textwidth}
\centering
{Number of $k$-cycles ($k=3$ to $6$)}\\
\begin{tabular}{lrrrr}
\toprule
Graph & $C_3$ & $C_4$ & $C_5$ & $C_6$    \\
\midrule
SG620x1   & 0  &  0  &  13888  &  146010  \\
SG620x2   &0  &  620  &  14508  &  138570 \\
$Z_{20}\rtimes_{2}Z_{31}$  & 0  &  620  &  12524  &  139810 \\
$Z_{20}\rtimes_{4}Z_{31}$  &  0  &  620  &  15748  &  146010 \\
$Z_{20}\rtimes_{8}Z_{31}$  & 0  &  0  &  16244  &  144150  \\
$Z_{20}\rtimes_{15}Z_{31}$  &  0  &  620  &  17484  &  130510 \\
$Z_{20}\rtimes_{16}Z_{31}$  &  0  &  620  &  16864  &  137950 \\
$Z_{20}\rtimes_{23}Z_{31}$  & 0  &  0  &  18848  &  144150  \\
$Z_{20}\rtimes_{27}Z_{31}$  &  0  &  0  &  17732  &  140430 \\
$Z_{20}\rtimes_{29}Z_{31}$  & 0  &  0  &  16864  &  145390  \\
\bottomrule
\end{tabular}
\end{minipage}
\end{table}

%%%%%%%%%%%%%%%%%%%%%%%%%%%%%%%%%%%%%%%%%%%%%%%%%%%%%%%%%%%%
\section{Degree 13}
%%%%%%%%%%%%%%%%%%%%%%%%%%%%%%%%%%%%%%%%%%%%%%%%%%%%%%%%%%%%

%%%%%%%%%%%%%%%%%%%%%%%
\subsection{\large  (13,2)=128}
Reference: \url{ https://sagecell.sagemath.org/?q=hmrdix }
\bigskip

Four non-isomorphic graphs to one another that improve the best current value of 112 in the table  (not listed as optimal).  Obtained from SG(128,17), SG(128,57), SG(128,553) and SG(128,1498) 

\begin{center}
\begin{tabular}{lccccc}
\toprule
Graph      & Avg.dist & Girth & Alg.conn. &  Aut.group & Edge-trans.\\
\midrule
SG128x17     & 1.89764  & 3  &  8.00000   &  512 & False \\ 
SG128x57     & 1.89764  & 3  &  8.00000   &  1024 & False \\ 
SG128x553x2    & 1.89764  & 3  &  8.00000   &  1024 & False \\ 
SG128x1498    & 1.89764  & 3  &  8.00000   &  128 & False \\ 
\bottomrule
\end{tabular}
\end{center}

\begin{table}[ht]
\centering
\begin{minipage}[t]{0.42\textwidth}
\centering
{Distance Distribution}\\
\smallskip
\begin{tabular}{lc}
\toprule
Graph & Dist. Distrib. \\
\midrule
SG128x17/57/553/1498   &   [1, 13, 114]  \\
\bottomrule
\end{tabular}
\end{minipage}
\hspace{0.02\textwidth} % smaller gap
\begin{minipage}[t]{0.40\textwidth}
\centering
{Number of $k$-cycles ($k=3$ to $6$)}\\
\begin{tabular}{lrrrr}
\toprule
Graph & $C_3$ & $C_4$ & $C_5$    & $C_6$  \\
\midrule
SG128x17x1    &   256  &  1472  &  26368  &  257024   \\
SG128x57x1    &   256  &  1472  &  26368  &  257024   \\
SG128x553x2   &   256  &  1472  &  26368  &  257024   \\
SG128x1498x1  &   256  &  1344  &  26112  &  262656   \\
\bottomrule
\end{tabular}
\end{minipage}
\end{table}

\subsection{\large  (13,3)=684}
Reference: \url{https://sagecell.sagemath.org/?q=jwbmpd}
\bigskip

Twenty-four non-isomorphic graphs to one another that improve the best current value of 680 in the table  (not listed as optimal).
Nine are obtained from SG(684,11) and all other are 
semidirect products
\texttt{$Z_{9}\rtimes_{(5,9,17,25,61,73)}Z_{76}$}.

\smallskip 

\begin{center}
\begin{tabular}{lcccc}
\toprule
Graph & Avg.dist & Girth & Alg.conn. & Aut.group  \\
\midrule
SG684x11xa                    & 2.73939  & 4  &  7.02964   &  684 \\ 
SG684x11xb                    & 2.74524  & 4  &  6.59290   &  684 \\ 
SG684x11xc                    & 2.73353  & 5  &  6.86303   &  684 \\ 
SG684x11xd                    & 2.73646  & 4  &  7.09437   &  684 \\ 
SG684x11xe                    & 2.73939  & 4  &  6.91724   &  684 \\ 
SG684x11xf                    & 2.73939  & 4  &  6.73054   &  684 \\ 
SG684x11xg                    & 2.73939  & 4  &  6.73227   &  684 \\ 
SG684x11xh                    & 2.73353  & 5  &  7.19245   &  684 \\ 
SG684x11xi                    & 2.75110  & 4  &  6.86389   &  684 \\ 
$Z_{9}\rtimes_{(5a)}Z_{76}$   & 2.75695  & 4  &  6.77332   &  684 \\ 
$Z_{9}\rtimes_{(5b)}Z_{76}$   & 2.75695  & 4  &  6.99450   &  684 \\ 
$Z_{9}\rtimes_{(5c)}Z_{76}$   & 2.75110  & 4  &  6.93582   &  684 \\ 
$Z_{9}\rtimes_{(9a)}Z_{76}$   & 2.75110  & 4  &  6.98545   &  684 \\ 
$Z_{9}\rtimes_{(9b)}Z_{76}$   & 2.75695  & 4  &  6.77332   &  684 \\ 
$Z_{9}\rtimes_{(17a)}Z_{76}$  & 2.75988  & 3  &  6.63081   &  684 \\ 
$Z_{9}\rtimes_{(17b)}Z_{76}$  & 2.75110  & 4  &  7.15898   &  684 \\ 
$Z_{9}\rtimes_{(25a)}Z_{76}$  & 2.75695  & 4  &  6.76203   &  684 \\ 
$Z_{9}\rtimes_{(25b)}Z_{76}$  & 2.75110  & 4  &  7.33029   &  684 \\ 
$Z_{9}\rtimes_{(61a)}Z_{76}$  & 2.75988  & 4  &  6.42275   &  684 \\ 
$Z_{9}\rtimes_{(61b)}Z_{76}$  & 2.75110  & 4  &  7.14228   &  684 \\ 
$Z_{9}\rtimes_{(61c)}Z_{76}$  & 2.75695  & 4  &  6.71388   &  684 \\ 
$Z_{9}\rtimes_{(61d)}Z_{76}$  & 2.75110  & 4  &  7.01993   &  684 \\ 
$Z_{9}\rtimes_{(73a)}Z_{76}$  & 2.75110  & 4  &  6.70914   &  684 \\ 
$Z_{9}\rtimes_{(73b)}Z_{76}$  & 2.75988  & 3  &  6.42275   &  684 \\ 
\bottomrule
\end{tabular}
\end{center}
\begin{table}[ht]
\centering
\begin{minipage}[t]{0.42\textwidth}
\centering
Distance Distribution
\medskip

\begin{tabular}{lc}
\toprule
Graph & Distance Distribution \\
\midrule
SG684x11xa                    &     [1, 13, 152, 518] \\ 
SG684x11xb                    &     [1, 13, 148, 522] \\ 
SG684x11xc                    &     [1, 13, 156, 514] \\ 
SG684x11xd                    &     [1, 13, 154, 516] \\ 
SG684x11xe                    &     [1, 13, 152, 518] \\ 
SG684x11xf                    &     [1, 13, 152, 518] \\ 
SG684x11xg                    &     [1, 13, 152, 518] \\ 
SG684x11xh                    &     [1, 13, 156, 514] \\ 
SG684x11xi                    &     [1, 13, 144, 526] \\ 
$Z_{9}\rtimes_{(5a)}Z_{76}$   &     [1, 13, 140, 530] \\ 
$Z_{9}\rtimes_{(5b)}Z_{76}$   &     [1, 13, 140, 530] \\ 
$Z_{9}\rtimes_{(5c)}Z_{76}$   &     [1, 13, 144, 526] \\ 
$Z_{9}\rtimes_{(9a)}Z_{76}$   &     [1, 13, 144, 526] \\ 
$Z_{9}\rtimes_{(9b)}Z_{76}$   &     [1, 13, 140, 530] \\ 
$Z_{9}\rtimes_{(17a)}Z_{76}$  &     [1, 13, 138, 532] \\ 
$Z_{9}\rtimes_{(17b)}Z_{76}$  &     [1, 13, 144, 526] \\ 
$Z_{9}\rtimes_{(25a)}Z_{76}$  &     [1, 13, 140, 530] \\ 
$Z_{9}\rtimes_{(25b)}Z_{76}$  &     [1, 13, 144, 526] \\ 
$Z_{9}\rtimes_{(61a)}Z_{76}$  &     [1, 13, 138, 532] \\ 
$Z_{9}\rtimes_{(61b)}Z_{76}$  &     [1, 13, 144, 526] \\ 
$Z_{9}\rtimes_{(61c)}Z_{76}$  &     [1, 13, 140, 530] \\ 
$Z_{9}\rtimes_{(61d)}Z_{76}$  &     [1, 13, 144, 526] \\ 
$Z_{9}\rtimes_{(73a)}Z_{76}$  &     [1, 13, 144, 526] \\ 
$Z_{9}\rtimes_{(73b)}Z_{76}$  &     [1, 13, 138, 532] \\ 
\bottomrule
\end{tabular}
\end{minipage}
\hspace{0.001\textwidth} % smaller gap
\begin{minipage}[t]{0.54\textwidth}
\centering
Number of $k$-cycles ($k=3$ to $6$)\\
\medskip

\begin{tabular}{lrrrr}
\toprule
Graph & $C_3$ & $C_4$ & $C_5$ & $C_6$ \\
\midrule
SG684x11xa                    &    0  &   684  &  23940  &  247266\\ 
SG684x11xb                    &    0  &   1368 &  29412  &  260262\\ 
SG684x11xc                    &    0  &     0  &  23256  &  259236\\ 
SG684x11xd                    &    0  &   342  &  24624  &  247266\\ 
SG684x11xe                    &    0  &   684  &  24624  &  257184\\ 
SG684x11xf                    &    0  &   684  &  22572  &  253422\\ 
SG684x11xg                    &    0  &   684  &  29412  &  237006\\ 
SG684x11xh                    &    0  &     0  &  23940  &  246924\\ 
SG684x11xi                    &    0  &  2052  &  25308  &  241794\\ 
$Z_{9}\rtimes_{(5a)}Z_{76}$   &    0  &  2736  &  23256  &  246012\\ 
$Z_{9}\rtimes_{(5b)}Z_{76}$   &    0  &  2736  &  25308  &  239856\\ 
$Z_{9}\rtimes_{(5c)}Z_{76}$   &    0  &  2052  &  21888  &  233700\\ 
$Z_{9}\rtimes_{(9a)}Z_{76}$   &    0  &  2052  &  26676  &  243504\\ 
$Z_{9}\rtimes_{(9b)}Z_{76}$   &    0  &  2736  &  23256  &  230736\\ 
$Z_{9}\rtimes_{(17a)}Z_{76}$  &  228  &  2736  &  21204  &  246696\\ 
$Z_{9}\rtimes_{(17b)}Z_{76}$  &    0  &  2052  &  17784  &  263568\\ 
$Z_{9}\rtimes_{(25a)}Z_{76}$  &    0  &  2736  &  22572  &  247608\\ 
$Z_{9}\rtimes_{(25b)}Z_{76}$  &    0  &  2052  &  20520  &  246696\\ 
$Z_{9}\rtimes_{(61a)}Z_{76}$  &    0  &  3420  &  22572  &  256500\\ 
$Z_{9}\rtimes_{(61b)}Z_{76}$  &    0  &  2052  &  21888  &  247152\\ 
$Z_{9}\rtimes_{(61c)}Z_{76}$  &    0  &  2736  &  28728  &  243276\\ 
$Z_{9}\rtimes_{(61d)}Z_{76}$  &    0  &  2052  &  21204  &  254676\\ 
$Z_{9}\rtimes_{(73a)}Z_{76}$  &    0  &  2052  &  19836  &  260832\\ 
$Z_{9}\rtimes_{(73b)}Z_{76}$  &  228  &  2736  &  25992  &  235980\\ 
\bottomrule
\end{tabular}
\end{minipage}
\end{table}

\clearpage

\newgeometry{
  left=1.8cm,
  right=1.8cm,
  top=1.8cm,
  bottom=1.8cm
}

\appendix

\section*{Appendix: Edge lists for non-isomorphic optimal Cayley graphs}
\addcontentsline{toc}{section}{Edge lists for non-isomorphic optimal Cayley graphs}

In this appendix, we present the edge lists of the non-isomorphic graphs found corresponding to the optimal degree--diameter Cayley graphs catalogued by Marston Conder on the Combinatorics Wiki: {\small \url{http://combinatoricswiki.org/wiki/Description_of_optimal_Cayley_graphs_found_by_Marston_Conder}}.

Exceptions are the case $(4,5)=216$, for which the graphs can be downloaded from \url{https://graphsym.net/}, and the case $(12,2)=96$, for which we include only six representative graphs out of the 42 distinct non-isomorphic examples.

Edge lists for all  cases considered in this paper are available in sparse6 format via the SageMath links provided for each degree--diameter pair.

\subsection*{(3,4)=24}

\noindent HG36335 / SG24x12x1

\begin{lstlisting}
{1,2},{1,3},{1,4},{2,7},{2,8},{3,6},{3,10},{4,5},{4,9},{5,12},{5,18},{6,11},{6,17},{7,14},{7,20},{8,13},
(8,19},{9,16},{9,21},{10,15},{10,21},{11,12},{11,13},{12,14},{13,16},{14,15},{15,23},{16,24},{17,20},
{17,24},{18,19},{18,23},{19,22},{20,22},{21,22},{23,24}
\end{lstlisting}

\noindent  HG36333 / SG24x12x2

\begin{lstlisting}
{1,2},{1,3},{1,4},{2,7},{2,8},{3,6},{3,9},{4,5},{4,9},{5,13},{5,14},{6,11},{6,15},{7,12},{7,13},{8,10},
{8,12},{9,16},{10,18},{10,21},{11,18},{11,20},{12,22},{13,19},{14,17},{14,19},{15,20},{15,24},{16,23},
{16,24},{17,18},{17,21},{19,20},{21,23},{22,23},{22,24}
\end{lstlisting}

\subsection*{(3,5)=60}

\noindent HG36684   / SG60x5x2

\begin{lstlisting}
{1,2},{1,3},{1,4},{2,6},{2,10},{3,5},{3,9},{4,7},{4,8},{5,13},{5,19},{6,12},{6,18},{7,15},{7,21},{8,14},
{8,20},{9,11},{9,17},{10,11},{10,16},{11,33},{12,24},{12,35},{13,25},{13,34},{14,25},{14,30},{15,24},
{15,29},{16,23},{16,37},{17,22},{17,36},{18,27},{18,32},{19,26},{19,31},{20,28},{20,39},{21,28},{21,38},
{22,42},{22,46},{23,43},{23,47},{24,41},{25,40},{26,56},{26,58},{27,55},{27,57},{28,48},{29,50},{29,56},
{30,49},{30,55},{31,52},{31,60},{32,51},{32,59},{33,53},{33,54},{34,35},{34,60},{35,59},{36,39},{36,58},
{37,38},{37,57},{38,44},{39,45},{40,41},{40,47},{41,46},{42,43},{42,45},{43,44},{44,60},{45,59},{46,50},
{47,49},{48,51},{48,52},{49,54},{50,53},{51,54},{52,53},{55,58},{56,57}
\end{lstlisting}

\noindent HG36698 / SG60x6x1  / $Z_{4}\rtimes_7  Z_{15}$ with generators $(1,1), (3,2), (2,0)$

\begin{lstlisting}
{1,2},{1,3},{1,4},{2,6},{2,9},{3,5},{3,7},{4,8},{4,10},{5,12},{5,21},{6,11},{6,20},{7,14},{7,18},{8,13},
{8,16},{9,15},{9,17},{10,19},{10,22},{11,23},{11,35},{12,24},{12,28},{13,27},{13,29},{14,25},{14,39},
{15,26},{15,37},{16,32},{16,40},{17,33},{17,45},{18,34},{18,44},{19,30},{19,31},{20,38},{20,42},{21,41},
{21,43},{22,36},{22,46},{23,28},{23,30},{24,26},{24,27},{25,29},{25,35},{26,36},{27,56},{28,57},{29,54},
{30,53},{31,33},{31,47},{32,34},{32,49},{33,49},{34,50},{35,51},{36,50},{37,48},{37,52},{38,47},{38,56},
{39,47},{39,48},{40,48},{40,57},{41,52},{41,55},{42,50},{42,55},{43,49},{43,51},{44,53},{44,59},{45,54},
{45,58},{46,51},{46,60},{52,53},{54,55},{56,59},{57,58},{58,60},{59,60}
\end{lstlisting}

\subsection*{(3,6)=72}

\noindent HG36950  / SG72x40 

\begin{lstlisting}
{1,2},{1,3},{1,4},{2,5},{2,6},{3,7},{3,8},{4,9},{4,10},{5,14},{5,22},{6,17},{6,19},{7,18},{7,20},{8,13},
{8,16},{9,11},{9,12},{10,15},{10,21},{11,29},{11,39},{12,24},{12,38},{13,25},{13,42},{14,26},{14,40},
{15,27},{15,28},{16,37},{16,43},{17,24},{17,31},{18,23},{18,31},{19,32},{19,41},{20,26},{20,34},{21,33},
{21,35},{22,30},{22,36},{23,30},{23,60},{24,25},{25,27},{26,55},{27,57},{28,45},{28,61},{29,44},{29,62},
{30,46},{31,33},{32,34},{32,58},{33,43},{34,59},{35,42},{35,47},{36,51},{36,59},{37,49},{37,56},{38,54},
{38,63},{39,45},{39,50},{40,48},{40,52},{41,53},{41,60},{42,65},{43,64},{44,55},{44,57},{45,46},{46,58},
{47,48},{47,51},{48,49},{49,50},{50,53},{51,72},{52,63},{52,71},{53,70},{54,61},{54,69},{55,69},{56,62},
{56,68},{57,67},{58,68},{59,67},{60,66},{61,66},{62,66},{63,67},{64,71},{64,72},{65,68},{65,70},{69,72},
{70,71}
\end{lstlisting}

\subsection*{(4,4)=84}

\noindent P84x182 / SG84x10

\begin{lstlisting}
{1,2},{1,3},{1,4},{1,5},{2,7},{2,15},{2,17},{3,6},{3,14},{3,16},{4,9},{4,10},{4,13},{5,8},{5,11},{5,12},
{6,19},{6,31},{6,53},{7,18},{7,30},{7,52},{8,35},{8,36},{8,39},{9,34},{9,37},{9,38},{10,21},{10,33},
{10,47},{11,20},{11,32},{11,46},{12,27},{12,41},{12,45},{13,26},{13,40},{13,44},{14,24},{14,28},{14,42},
{15,25},{15,29},{15,43},{16,23},{16,49},{16,51},{17,22},{17,48},{17,50},{18,32},{18,34},{18,56},{19,33},
{19,35},{19,55},{20,21},{20,22},{20,71},{21,23},{21,72},{22,28},{22,68},{23,29},{23,67},{24,25},{24,36},
{24,78},{25,37},{25,77},{26,31},{26,56},{26,68},{27,30},{27,55},{27,67},{28,62},{28,69},{29,61},{29,70},
{30,59},{30,64},{31,60},{31,63},{32,54},{32,80},{33,54},{33,79},{34,62},{34,73},{35,61},{35,74},{36,56},
{36,76},{37,55},{37,75},{38,60},{38,66},{38,70},{39,59},{39,65},{39,69},{40,57},{40,61},{40,64},{41,58},
{41,62},{41,63},{42,54},{42,64},{42,66},{43,54},{43,63},{43,65},{44,67},{44,69},{44,80},{45,68},{45,70},
{45,79},{46,57},{46,60},{46,78},{47,58},{47,59},{47,77},{48,55},{48,57},{48,76},{49,56},{49,58},{49,75},
{50,58},{50,66},{50,84},{51,57},{51,65},{51,83},{52,72},{52,74},{52,82},{53,71},{53,73},{53,81},{59,60},
{61,62},{63,72},{64,71},{65,73},{66,74},{67,84},{68,83},{69,82},{70,81},{71,75},{72,76},{73,84},{74,83},
{75,82},{76,81},{77,80},{77,83},{78,79},{78,84},{79,82},{80,81}
\end{lstlisting}

\subsection*{(5,4)=210}
\def\-{\raisebox{.75pt}{-}}

\noindent SG(210,2) /  $Z_{6}\rtimes_4  Z_{35}$ with generators $[1,0]\-[5,0];[2,1]\-[4,24];[3,21]$

\begin{lstlisting}
{1,8},{1,22},{1,27},{1,103},{1,144},{2,14},{2,41},{2,46},{2,127},{2,161},{3,20},{3,57},{3,62},{3,72},
{3,186},{4,7},{4,24},{4,67},{4,166},{4,194},{5,6},{5,28},{5,49},{5,77},{5,88},{6,23},{6,31},{6,59},
{6,186},{7,26},{7,34},{7,94},{7,142},{8,45},{8,75},{8,81},{8,83},{9,39},{9,75},{9,82},{9,93},{9,142},
{10,13},{10,43},{10,99},{10,185},{10,201},{11,12},{11,47},{11,66},{11,98},{11,107},{12,24},{12,42},
{12,49},{12,144},{13,45},{13,51},{13,141},{13,156},{14,62},{14,96},{14,105},{14,110},{15,55},{15,83},
{15,102},{15,125},{15,170},{16,19},{16,59},{16,122},{16,137},{16,209},{17,18},{17,31},{17,63},{17,129},
{17,134},{18,58},{18,64},{18,119},{18,170},{19,43},{19,61},{19,66},{19,161},{20,26},{20,127},{20,131},
{20,132},{21,23},{21,69},{21,76},{21,172},{21,198},{22,68},{22,72},{22,124},{22,209},{23,75},{23,160},
{23,192},{24,31},{24,33},{24,140},{25,26},{25,68},{25,78},{25,110},{25,154},{26,81},{26,134},{27,39},
{27,76},{27,84},{27,86},{28,34},{28,106},{28,147},{28,149},{29,85},{29,105},{29,116},{29,126},{29,156},
{30,75},{30,86},{30,90},{30,94},{30,104},{31,84},{31,132},{32,36},{32,87},{32,96},{32,169},{32,178},
{33,53},{33,88},{33,103},{33,109},{34,37},{34,72},{34,80},{35,71},{35,91},{35,105},{35,111},{35,141},
{36,94},{36,140},{36,159},{36,192},{37,38},{37,86},{37,97},{37,148},{38,43},{38,73},{38,98},{38,127},
{39,80},{39,143},{39,164},{40,42},{40,101},{40,106},{40,189},{40,202},{41,69},{41,100},{41,103},
{41,194},{42,105},{42,191},{42,199},{43,49},{43,168},{44,45},{44,100},{44,108},{44,133},{44,164},
{45,88},{45,110},{46,85},{46,91},{46,106},{46,112},{47,51},{47,129},{47,162},{47,172},{48,54},{48,74},
{48,113},{48,143},{48,178},{49,83},{49,112},{50,80},{50,89},{50,114},{50,119},{50,125},{51,79},{51,90},
{51,103},{52,70},{52,117},{52,132},{52,143},{52,169},{53,119},{53,165},{53,183},{53,204},{54,64},
{54,114},{54,121},{54,149},{55,61},{55,111},{55,144},{55,175},{56,58},{56,124},{56,128},{56,147},
{56,210},{57,123},{57,125},{57,181},{57,204},{58,101},{58,127},{58,201},{59,64},{59,65},{59,185},
{60,61},{60,81},{60,123},{60,130},{60,190},{61,114},{61,131},{62,98},{62,128},{62,133},{63,66},{63,76},
{63,187},{63,188},{64,111},{64,133},{65,72},{65,79},{65,115},{65,134},{66,96},{66,125},{67,77},
{67,130},{67,137},{67,210},{68,140},{68,163},{68,196},{69,72},{69,74},{69,150},{70,141},{70,168},
{70,184},{70,199},{71,140},{71,142},{71,160},{71,167},{73,92},{73,143},{73,203},{73,207},{74,118},
{74,144},{74,166},{75,137},{76,147},{76,190},{77,97},{77,150},{77,152},{78,84},{78,99},{78,137},
{78,139},{79,91},{79,151},{79,179},{80,152},{80,157},{81,150},{81,188},{82,113},{82,117},{82,153},
{82,162},{83,154},{83,171},{84,95},{84,146},{85,90},{85,178},{85,190},{86,156},{86,176},{87,97},
{87,109},{87,116},{87,161},{88,121},{88,169},{89,136},{89,156},{89,168},{89,191},{90,158},{90,173},
{91,152},{91,170},{92,150},{92,160},{92,162},{92,196},{93,101},{93,138},{93,161},{93,185},{94,98},
{94,193},{95,96},{95,138},{95,152},{95,163},{96,164},{97,146},{97,202},{98,153},{99,107},{99,166},
{99,198},{100,168},{100,195},{100,200},{101,103},{101,174},{102,139},{102,169},{102,193},{102,207},
{104,155},{104,170},{104,181},{104,183},{105,166},{106,154},{106,172},{107,151},{107,158},{107,174},
{108,112},{108,123},{108,159},{108,166},{109,117},{109,120},{109,175},{110,149},{110,174},{111,146},
{111,176},{112,145},{112,157},{113,116},{113,119},{113,154},{114,142},{114,158},{115,135},{115,178},
{115,193},{115,203},{116,180},{116,188},{117,151},{117,186},{118,171},{118,181},{118,187},{118,206},
{119,194},{120,131},{120,139},{120,176},{120,182},{121,128},{121,173},{121,205},{122,129},{122,182},
{122,183},{122,206},{123,185},{123,202},{124,125},{124,126},{124,189},{126,137},{126,177},{126,186},
{127,183},{128,176},{128,187},{129,164},{129,189},{130,133},{130,183},{130,184},{131,173},{131,189},
{132,145},{132,190},{133,162},{134,141},{134,180},{135,174},{135,188},{135,191},{135,200},{136,150},
{136,163},{136,173},{136,192},{138,193},{138,205},{138,208},{139,172},{139,194},{140,147},{141,207},
{142,199},{143,201},{144,203},{145,158},{145,177},{145,195},{146,155},{146,196},{147,184},{148,159},
{148,175},{148,180},{148,197},{149,165},{149,198},{151,157},{151,210},{152,206},{153,163},{153,171},
{153,179},{154,182},{155,174},{155,195},{155,199},{156,204},{157,167},{157,200},{158,196},{159,189},
{159,201},{160,161},{160,197},{162,202},{163,201},{164,197},{165,197},{165,203},{165,208},{167,182},
{167,187},{167,204},{168,172},{169,181},{170,192},{171,180},{171,205},{173,206},{175,179},{175,198},
{176,200},{177,197},{177,205},{177,207},{178,209},{179,184},{179,208},{180,195},{181,198},{182,203},
{184,209},{185,187},{186,191},{188,210},{190,208},{191,208},{192,200},{193,202},{194,205},{195,209},
{196,204},{199,206},{207,210}
\end{lstlisting}

\subsection*{(7,3)= 168}

\noindent SG168x43x1 

\begin{lstlisting}
{1,2},{1,5},{1,6},{1,41},{1,55},{1,120},{1,159},{2,11},{2,12},{2,17},{2,70},{2,88},{2,144},{3,9},{3,12},
{3,16},{3,17},{3,20},{3,147},{3,163},{4,10},{4,16},{4,20},{4,21},{4,61},{4,127},{4,157},{5,15},{5,16},
{5,25},{5,54},{5,63},{5,116},{6,9},{6,28},{6,29},{6,74},{6,100},{6,155},{7,32},{7,34},{7,35},{7,66},
{7,71},{7,93},{7,131},{8,23},{8,38},{8,39},{8,45},{8,69},{8,111},{8,139},{9,24},{9,41},{9,43},{9,134},
{9,167},{10,24},{10,28},{10,43},{10,56},{10,99},{10,141},{11,24},{11,39},{11,40},{11,44},{11,84},
{11,110},{12,21},{12,47},{12,77},{12,91},{12,168},{13,22},{13,50},{13,51},{13,76},{13,102},{13,113},
{13,138},{14,26},{14,37},{14,53},{14,54},{14,68},{14,112},{14,165},{15,28},{15,40},{15,55},{15,56},
{15,154},{15,162},{16,30},{16,58},{16,112},{16,113},{17,61},{17,62},{17,123},{17,128},{17,164},{18,34},
{18,65},{18,66},{18,67},{18,78},{18,84},{18,98},{19,36},{19,39},{19,53},{19,68},{19,69},{19,82},
{19,149},{20,55},{20,71},{20,81},{20,108},{20,125},{21,29},{21,40},{21,50},{21,74},{21,148},{22,25},
{22,49},{22,67},{22,77},{22,78},{22,104},{23,79},{23,80},{23,85},{23,106},{23,114},{23,129},{24,48},
{24,69},{24,71},{24,81},{25,64},{25,70},{25,71},{25,85},{25,165},{26,34},{26,45},{26,88},{26,128},
{26,141},{26,145},{27,31},{27,76},{27,83},{27,90},{27,133},{27,166},{27,168},{28,34},{28,47},{28,91},
{28,135},{29,94},{29,95},{29,105},{29,136},{29,138},{30,35},{30,51},{30,86},{30,96},{30,130},{30,146},
{31,32},{31,58},{31,64},{31,73},{31,99},{31,127},{32,44},{32,100},{32,101},{32,103},{32,128},{33,50},
{33,57},{33,101},{33,102},{33,103},{33,116},{33,126},{34,102},{34,119},{34,148},{35,48},{35,57},{35,78},
{35,79},{35,94},{36,80},{36,105},{36,106},{36,107},{36,119},{36,132},{37,72},{37,105},{37,107},{37,108},
{37,117},{37,130},{38,42},{38,78},{38,99},{38,105},{38,123},{38,163},{39,54},{39,90},{39,96},{39,164},
{40,61},{40,70},{40,137},{40,153},{41,56},{41,94},{41,104},{41,109},{41,142},{42,52},{42,54},{42,59},
{42,111},{42,112},{42,151},{43,58},{43,70},{43,80},{43,96},{43,113},{44,55},{44,97},{44,113},{44,117},
{44,139},{45,50},{45,104},{45,120},{45,132},{45,157},{46,49},{46,60},{46,93},{46,122},{46,136},{46,146},
{46,155},{47,123},{47,124},{47,129},{47,131},{47,166},{48,51},{48,106},{48,118},{48,125},{48,133},
{49,81},{49,97},{49,99},{49,115},{49,127},{50,66},{50,87},{50,135},{51,65},{51,101},{51,107},{51,123},
{52,79},{52,87},{52,108},{52,129},{52,130},{52,145},{53,94},{53,101},{53,127},{53,129},{53,167},
{54,122},{54,125},{54,142},{55,75},{55,84},{55,130},{56,62},{56,102},{56,122},{56,123},{57,58},{57,77},
{57,97},{57,104},{57,131},{58,84},{58,132},{58,151},{59,66},{59,82},{59,100},{59,117},{59,134},{59,153},
{60,64},{60,86},{60,124},{60,126},{60,135},{60,142},{61,66},{61,90},{61,94},{61,109},{62,79},{62,91},
{62,95},{62,101},{62,136},{63,67},{63,92},{63,103},{63,108},{63,118},{63,133},{64,65},{64,121},{64,137},
{64,162},{65,81},{65,100},{65,128},{65,138},{66,114},{66,122},{67,91},{67,105},{67,110},{67,131},
{68,74},{68,79},{68,131},{68,137},{68,159},{69,92},{69,112},{69,135},{69,155},{70,92},{70,106},{70,116},
{71,117},{71,126},{71,139},{72,77},{72,88},{72,102},{72,114},{72,141},{72,151},{73,76},{73,89},{73,135},
{73,142},{73,150},{73,160},{74,103},{74,109},{74,143},{74,161},{75,78},{75,101},{75,110},{75,112},
{75,140},{75,152},{76,84},{76,126},{76,144},{76,159},{77,82},{77,128},{77,155},{78,92},{78,109},
{79,144},{79,162},{80,108},{80,110},{80,122},{80,142},{81,116},{81,145},{81,149},{82,85},{82,102},
{82,147},{82,154},{83,95},{83,97},{83,98},{83,118},{83,148},{83,164},{84,145},{84,149},{85,95},{85,113},
{85,124},{85,145},{86,88},{86,110},{86,133},{86,157},{86,161},{87,104},{87,119},{87,134},{87,153},
{87,165},{88,153},{88,158},{88,167},{89,97},{89,98},{89,122},{89,143},{89,155},{89,158},{90,102},
{90,108},{90,121},{90,155},{91,111},{91,156},{91,160},{92,101},{92,150},{92,158},{93,98},{93,115},
{93,116},{93,144},{93,159},{94,156},{94,160},{95,117},{95,150},{95,161},{96,131},{96,138},{96,140},
{96,152},{97,137},{97,162},{98,113},{98,163},{98,167},{99,129},{99,153},{99,159},{100,104},{100,114},
{100,168},{103,110},{103,129},{103,138},{104,164},{105,127},{105,167},{106,130},{106,148},{106,168},
{107,109},{107,111},{107,137},{107,159},{108,164},{109,124},{109,166},{110,146},{111,138},{111,144},
{111,162},{112,148},{112,168},{113,165},{114,120},{114,149},{114,157},{115,121},{115,148},{115,152},
{115,161},{115,164},{116,132},{116,151},{117,124},{117,132},{118,120},{118,141},{118,146},{118,160},
{119,128},{119,139},{119,147},{119,154},{120,140},{120,154},{120,163},{121,126},{121,140},{121,156},
{121,168},{122,168},{123,143},{123,161},{124,152},{124,160},{125,131},{125,138},{125,150},{125,158},
{126,163},{126,167},{127,144},{127,154},{128,142},{129,163},{130,135},{130,155},{132,136},{132,166},
{133,134},{133,143},{133,154},{134,141},{134,144},{134,152},{135,164},{136,140},{136,143},{136,145},
{137,141},{137,163},{139,143},{139,151},{139,156},{140,141},{140,150},{142,148},{143,165},{145,166},
{146,147},{146,153},{146,156},{147,150},{147,157},{147,159},{149,160},{149,161},{149,165},{150,153},
{151,160},{151,161},{152,154},{152,158},{156,165},{156,166},{157,158},{157,162},{158,166},{162,167}
\end{lstlisting}

\noindent SG168x43x2 

\begin{lstlisting}
{1,2},{1,5},{1,6},{1,43},{1,61},{1,68},{1,75},{2,11},{2,12},{2,20},{2,56},{2,92},{2,111},{3,9},{3,16},
{3,17},{3,28},{3,38},{3,70},{3,125},{4,10},{4,11},{4,20},{4,21},{4,30},{4,41},{4,105},{5,10},{5,16},
{5,25},{5,33},{5,44},{5,142},{6,28},{6,29},{6,40},{6,46},{6,123},{6,134},{7,24},{7,32},{7,34},{7,35},
{7,77},{7,122},{7,138},{8,23},{8,38},{8,39},{8,53},{8,91},{8,104},{8,130},{9,21},{9,24},{9,41},{9,53},
{9,55},{9,96},{10,17},{10,28},{10,43},{10,48},{10,129},{11,18},{11,24},{11,25},{11,44},{11,164},{12,15},
{12,21},{12,27},{12,47},{12,94},{12,147},{13,16},{13,22},{13,50},{13,51},{13,90},{13,100},{13,131},
{14,37},{14,38},{14,53},{14,54},{14,74},{14,106},{14,128},{15,24},{15,40},{15,55},{15,56},{15,63},
{15,79},{16,40},{16,58},{16,81},{16,168},{17,61},{17,62},{17,74},{17,83},{17,88},{18,65},{18,66},
{18,67},{18,101},{18,104},{18,148},{19,36},{19,68},{19,69},{19,77},{19,94},{19,108},{19,111},{20,55},
{20,65},{20,71},{20,90},{20,113},{21,62},{21,74},{21,115},{21,153},{22,34},{22,70},{22,77},{22,78},
{22,103},{22,164},{23,47},{23,66},{23,79},{23,80},{23,107},{23,112},{24,58},{24,81},{24,155},{25,71},
{25,85},{25,87},{25,132},{25,162},{26,34},{26,45},{26,84},{26,88},{26,100},{26,137},{26,147},{27,31},
{27,83},{27,90},{27,93},{27,135},{27,154},{28,73},{28,91},{28,109},{28,154},{29,37},{29,50},{29,56},
{29,94},{29,95},{29,166},{30,35},{30,67},{30,86},{30,96},{30,107},{30,150},{31,64},{31,99},{31,103},
{31,115},{31,160},{31,167},{32,50},{32,55},{32,67},{32,100},{32,101},{32,142},{33,57},{33,78},{33,102},
{33,103},{33,128},{33,135},{34,36},{34,47},{34,102},{34,114},{35,63},{35,65},{35,78},{35,97},{35,156},
{36,54},{36,105},{36,106},{36,109},{36,129},{37,69},{37,79},{37,102},{37,107},{37,108},{38,92},{38,99},
{38,105},{38,144},{39,52},{39,54},{39,90},{39,140},{39,148},{39,149},{40,61},{40,70},{40,107},{40,110},
{41,56},{41,60},{41,91},{41,109},{41,120},{42,52},{42,68},{42,80},{42,100},{42,111},{42,112},{42,123},
{43,57},{43,70},{43,71},{43,113},{43,122},{44,113},{44,117},{44,119},{44,137},{44,145},{45,50},{45,77},
{45,116},{45,120},{45,134},{45,162},{46,49},{46,60},{46,76},{46,122},{46,148},{46,153},{47,61},{47,123},
{47,124},{47,136},{48,51},{48,79},{48,103},{48,118},{48,125},{48,152},{49,67},{49,73},{49,97},{49,127},
{49,161},{49,163},{50,52},{50,66},{50,72},{51,57},{51,64},{51,101},{51,110},{51,143},{52,62},{52,105},
{52,129},{52,130},{53,75},{53,127},{53,129},{53,159},{54,122},{54,135},{54,151},{54,158},{55,84},
{55,116},{55,135},{56,89},{56,123},{56,141},{57,66},{57,104},{57,131},{57,155},{58,72},{58,84},{58,85},
{58,127},{58,132},{59,66},{59,71},{59,82},{59,99},{59,120},{59,128},{59,134},{60,64},{60,90},{60,98},
{60,135},{60,157},{61,94},{61,121},{61,157},{62,91},{62,102},{62,124},{62,136},{63,67},{63,92},{63,129},
{63,133},{63,140},{64,89},{64,137},{64,144},{64,156},{65,102},{65,128},{65,138},{65,168},{66,109},
{66,119},{67,131},{67,160},{68,79},{68,96},{68,137},{68,167},{69,112},{69,122},{69,135},{69,139},
{69,150},{70,84},{70,116},{70,148},{71,139},{71,149},{71,163},{72,77},{72,114},{72,141},{72,154},
{72,167},{73,76},{73,89},{73,142},{73,147},{73,168},{74,100},{74,109},{74,143},{74,160},{75,78},
{75,110},{75,131},{75,133},{75,140},{76,126},{76,137},{76,138},{76,144},{76,166},{77,91},{77,128},
{78,98},{78,124},{78,125},{79,99},{79,144},{80,86},{80,108},{80,132},{80,142},{80,155},{81,99},{81,114},
{81,116},{81,117},{81,145},{82,88},{82,102},{82,104},{82,113},{82,127},{82,147},{83,97},{83,122},
{83,126},{83,141},{83,148},{84,139},{84,149},{84,159},{85,95},{85,106},{85,143},{85,145},{85,152},
{86,88},{86,120},{86,125},{86,133},{86,151},{87,102},{87,104},{87,119},{87,144},{87,153},{87,157},
{88,89},{88,118},{88,153},{89,98},{89,155},{89,164},{90,112},{90,155},{91,156},{91,161},{92,101},
{92,138},{92,146},{92,158},{93,98},{93,115},{93,131},{93,136},{93,159},{93,162},{94,104},{94,143},
{94,160},{95,98},{95,101},{95,109},{95,139},{95,161},{96,101},{96,118},{96,138},{96,152},{97,121},
{97,143},{97,159},{97,162},{98,127},{98,163},{99,126},{99,159},{100,104},{100,114},{101,126},{103,110},
{103,138},{103,161},{105,110},{105,137},{105,167},{106,130},{106,133},{106,164},{106,168},{107,111},
{107,127},{107,159},{108,118},{108,145},{108,164},{108,168},{109,166},{110,146},{110,158},{111,125},
{111,162},{111,163},{112,148},{112,152},{112,165},{113,151},{113,165},{113,167},{114,153},{114,157},
{114,163},{115,121},{115,134},{115,155},{115,164},{116,144},{116,151},{116,165},{117,124},{117,130},
{117,132},{117,150},{117,156},{118,120},{118,146},{118,149},{119,128},{119,141},{119,154},{119,159},
{120,121},{120,154},{121,126},{121,142},{121,168},{122,168},{123,128},{123,156},{123,161},{124,126},
{124,160},{124,165},{125,131},{125,150},{126,167},{127,144},{129,162},{129,163},{130,142},{130,146},
{130,155},{131,166},{132,136},{132,158},{132,161},{133,134},{133,147},{133,165},{134,141},{134,146},
{135,164},{136,138},{136,143},{136,149},{137,163},{139,143},{139,146},{139,151},{140,141},{140,145},
{140,150},{140,157},{141,158},{142,148},{145,160},{145,166},{146,147},{147,157},{149,161},{149,165},
{150,153},{150,154},{151,160},{151,166},{152,153},{152,154},{152,158},{156,165},{156,166},{157,158},
{162,167}
\end{lstlisting}

\noindent SG168x43x4 

\begin{lstlisting}
{1,2},{1,5},{1,6},{1,67},{1,68},{1,75},{1,89},{2,11},{2,12},{2,92},{2,103},{2,111},{2,121},{3,9},
{3,16},{3,17},{3,35},{3,38},{3,115},{3,125},{4,10},{4,20},{4,21},{4,30},{4,83},{4,105},{4,131},{5,16},
{5,25},{5,33},{5,77},{5,107},{5,142},{6,28},{6,29},{6,36},{6,46},{6,118},{6,134},{7,24},{7,32},{7,34},
{7,35},{7,112},{7,122},{7,160},{8,23},{8,38},{8,39},{8,56},{8,91},{8,104},{8,152},{9,24},{9,41},
{9,51},{9,53},{9,73},{9,96},{10,28},{10,43},{10,48},{10,60},{10,129},{10,138},{11,18},{11,24},{11,44},
{11,79},{11,100},{11,164},{12,21},{12,27},{12,47},{12,52},{12,86},{12,147},{13,16},{13,22},{13,50},
{13,51},{13,69},{13,90},{13,161},{14,37},{14,53},{14,54},{14,61},{14,74},{14,128},{14,150},{15,40},
{15,46},{15,55},{15,56},{15,63},{15,78},{15,79},{16,58},{16,104},{16,129},{16,168},{17,23},{17,61},
{17,62},{17,83},{17,88},{17,146},{18,54},{18,65},{18,66},{18,67},{18,148},{18,156},{19,21},{19,36},
{19,68},{19,69},{19,77},{19,94},{19,158},{20,34},{20,53},{20,55},{20,65},{20,71},{20,90},{21,74},
{21,115},{21,153},{21,158},{22,70},{22,77},{22,78},{22,130},{22,164},{22,166},{23,47},{23,66},{23,79},
{23,80},{23,146},{24,81},{24,105},{24,128},{24,155},{25,71},{25,85},{25,87},{25,122},{25,147},{25,162},
{26,34},{26,45},{26,84},{26,88},{26,91},{26,126},{26,137},{27,31},{27,40},{27,83},{27,90},{27,101},
{27,154},{28,42},{28,73},{28,91},{28,140},{28,154},{29,37},{29,50},{29,94},{29,95},{29,97},{29,114},
{30,35},{30,86},{30,96},{30,107},{30,144},{30,149},{31,64},{31,87},{31,99},{31,103},{31,109},{31,160},
{32,55},{32,100},{32,101},{32,106},{32,136},{32,142},{33,39},{33,57},{33,102},{33,103},{33,135},
{33,143},{34,36},{34,47},{34,53},{34,102},{35,78},{35,97},{35,115},{35,156},{36,105},{36,106},
{36,109},{36,118},{37,41},{37,102},{37,107},{37,108},{37,133},{38,43},{38,50},{38,92},{38,99},{38,105},
{39,54},{39,90},{39,140},{39,143},{39,149},{40,61},{40,70},{40,101},{40,107},{40,110},{41,56},{41,60},
{41,109},{41,120},{41,133},{42,52},{42,100},{42,111},{42,112},{42,123},{42,140},{43,50},{43,57},
{43,70},{43,113},{43,122},{44,90},{44,113},{44,117},{44,119},{44,134},{44,137},{45,50},{45,74},{45,98},
{45,116},{45,120},{45,162},{46,49},{46,60},{46,78},{46,122},{46,153},{47,64},{47,72},{47,123},{47,124},
{48,51},{48,79},{48,118},{48,125},{48,151},{48,159},{49,62},{49,67},{49,97},{49,119},{49,127},{49,161},
{50,52},{50,66},{51,64},{51,73},{51,101},{51,143},{52,62},{52,86},{52,129},{52,130},{53,75},{53,127},
{53,129},{54,122},{54,151},{54,156},{54,158},{55,66},{55,84},{55,111},{55,135},{56,89},{56,123},
{56,141},{56,152},{57,104},{57,108},{57,124},{57,131},{57,155},{58,72},{58,84},{58,120},{58,127},
{58,132},{58,148},{59,66},{59,71},{59,76},{59,82},{59,94},{59,99},{59,134},{60,64},{60,135},{60,138},
{60,157},{61,94},{61,121},{61,150},{61,157},{62,91},{62,102},{62,119},{62,136},{63,67},{63,92},
{63,129},{63,133},{63,139},{63,167},{64,72},{64,137},{64,156},{65,80},{65,95},{65,128},{65,138},
{65,168},{66,109},{66,111},{67,89},{67,131},{67,160},{68,70},{68,79},{68,96},{68,102},{68,137},
{69,112},{69,135},{69,139},{69,150},{69,161},{70,102},{70,116},{70,148},{71,139},{71,154},{71,163},
{71,164},{72,77},{72,114},{72,141},{72,167},{73,76},{73,89},{73,142},{73,147},{74,98},{74,100},
{74,109},{74,143},{75,78},{75,99},{75,110},{75,132},{75,140},{76,94},{76,126},{76,138},{76,144},
{76,166},{77,91},{77,107},{77,128},{78,98},{78,124},{79,100},{79,144},{80,86},{80,95},{80,108},
{80,132},{80,142},{81,88},{81,99},{81,114},{81,116},{81,135},{81,145},{82,93},{82,102},{82,113},
{82,123},{82,127},{82,147},{83,97},{83,131},{83,141},{83,148},{84,149},{84,155},{84,157},{84,159},
{85,95},{85,106},{85,125},{85,137},{85,145},{85,152},{86,88},{86,133},{86,151},{87,104},{87,109},
{87,119},{87,144},{87,153},{88,89},{88,135},{88,153},{89,98},{89,155},{90,134},{90,155},{91,126},
{91,156},{92,101},{92,127},{92,145},{92,158},{93,98},{93,115},{93,123},{93,131},{93,136},{93,159},
{94,104},{94,160},{95,98},{95,101},{95,161},{96,117},{96,138},{96,152},{96,162},{97,114},{97,143},
{97,162},{98,163},{99,132},{99,159},{100,104},{100,114},{101,126},{103,110},{103,121},{103,138},
{103,161},{104,129},{105,110},{105,128},{105,167},{106,130},{106,133},{106,136},{106,168},{107,111},
{107,159},{108,118},{108,124},{108,145},{108,164},{109,166},{110,146},{110,163},{110,165},{111,125},
{111,162},{112,148},{112,152},{112,160},{112,165},{113,142},{113,153},{113,165},{113,167},{114,157},
{114,163},{115,121},{115,134},{115,164},{116,141},{116,144},{116,151},{116,168},{117,124},{117,130},
{117,132},{117,150},{117,162},{118,120},{118,146},{118,149},{119,128},{119,154},{119,159},{120,121},
{120,148},{120,154},{121,126},{121,168},{122,147},{122,168},{123,128},{123,161},{124,126},{124,160},
{125,131},{125,137},{125,150},{126,167},{127,144},{127,145},{129,163},{130,146},{130,155},{130,166},
{131,166},{132,136},{132,158},{133,134},{133,165},{134,141},{135,164},{136,138},{136,143},{137,163},
{139,143},{139,146},{139,151},{139,167},{140,141},{140,145},{140,150},{141,168},{142,148},{142,153},
{144,149},{145,166},{146,147},{147,157},{149,161},{149,165},{150,153},{151,159},{151,160},{152,154},
{152,158},{154,164},{155,157},{156,165},{156,166},{157,158},{162,167},{163,165}
\end{lstlisting}

\noindent SG168x43x5 

\begin{lstlisting}
{1,2},{1,5},{1,6},{1,109},{1,116},{1,122},{1,133},{2,11},{2,12},{2,62},{2,84},{2,90},{2,146},{3,9},
{3,16},{3,17},{3,47},{3,86},{3,113},{3,148},{4,10},{4,20},{4,21},{4,81},{4,94},{4,150},{4,164},{5,16},
{5,25},{5,56},{5,78},{5,137},{5,165},{6,28},{6,29},{6,45},{6,55},{6,107},{6,143},{7,32},{7,34},{7,35},
{7,83},{7,91},{7,119},{7,125},{8,23},{8,38},{8,39},{8,75},{8,135},{8,161},{8,163},{9,24},{9,29},{9,41},
{9,71},{9,118},{9,135},{10,28},{10,43},{10,58},{10,123},{10,142},{10,152},{11,24},{11,44},{11,61},
{11,101},{11,139},{11,162},{12,21},{12,26},{12,47},{12,70},{12,79},{12,156},{13,22},{13,50},{13,51},
{13,60},{13,74},{13,87},{13,96},{14,37},{14,53},{14,54},{14,92},{14,148},{14,160},{14,167},{15,40},
{15,44},{15,55},{15,56},{15,91},{15,140},{15,155},{16,21},{16,58},{16,99},{16,131},{16,151},{17,20},
{17,61},{17,62},{17,82},{17,129},{17,161},{18,46},{18,65},{18,66},{18,67},{18,75},{18,94},{18,114},
{19,36},{19,68},{19,69},{19,90},{19,125},{19,143},{19,159},{20,35},{20,55},{20,71},{20,117},{20,167},
{21,53},{21,74},{21,95},{21,114},{22,47},{22,63},{22,77},{22,78},{22,82},{22,89},{23,30},{23,79},
{23,80},{23,136},{23,162},{23,168},{24,28},{24,81},{24,127},{24,138},{24,149},{25,40},{25,71},{25,80},
{25,85},{25,102},{25,124},{26,34},{26,45},{26,65},{26,88},{26,105},{26,158},{27,31},{27,33},{27,83},
{27,90},{27,94},{27,108},{27,144},{28,38},{28,72},{28,91},{28,124},{29,32},{29,89},{29,94},{29,95},
{29,132},{30,35},{30,57},{30,86},{30,96},{30,147},{30,160},{31,53},{31,64},{31,99},{31,134},{31,145},
{31,168},{32,59},{32,100},{32,101},{32,110},{32,121},{33,57},{33,72},{33,92},{33,102},{33,103},
{33,123},{34,54},{34,67},{34,81},{34,102},{34,141},{35,66},{35,78},{35,118},{35,159},{36,63},{36,95},
{36,105},{36,106},{36,127},{36,142},{37,48},{37,107},{37,108},{37,137},{37,155},{37,166},{38,49},
{38,87},{38,99},{38,105},{38,112},{39,50},{39,54},{39,71},{39,90},{39,111},{39,115},{40,61},{40,70},
{40,74},{40,158},{40,168},{41,43},{41,56},{41,59},{41,105},{41,109},{41,160},{42,52},{42,96},{42,111},
{42,112},{42,122},{42,144},{42,156},{43,51},{43,70},{43,85},{43,113},{43,163},{44,66},{44,95},{44,108},
{44,113},{44,117},{45,50},{45,57},{45,120},{45,129},{45,140},{46,49},{46,60},{46,80},{46,122},{46,123},
{46,159},{47,121},{47,123},{47,124},{47,145},{48,51},{48,65},{48,118},{48,125},{48,134},{48,161},
{49,97},{49,127},{49,132},{49,147},{49,155},{50,58},{50,66},{50,103},{50,157},{51,86},{51,101},
{51,102},{51,144},{52,99},{52,110},{52,124},{52,129},{52,130},{52,164},{53,69},{53,119},{53,127},
{53,129},{54,68},{54,73},{54,113},{54,122},{55,67},{55,84},{55,144},{55,145},{56,111},{56,119},
{56,123},{56,136},{57,104},{57,109},{57,115},{57,131},{58,84},{58,106},{58,132},{58,166},{59,66},
{59,79},{59,82},{59,134},{59,152},{60,64},{60,91},{60,130},{60,135},{60,167},{61,68},{61,87},{61,94},
{61,166},{62,65},{62,85},{62,91},{62,115},{62,136},{63,67},{63,92},{63,120},{63,133},{63,156},{64,88},
{64,111},{64,117},{64,137},{64,142},{65,73},{65,128},{65,138},{66,112},{66,153},{67,131},{67,146},
{67,163},{68,72},{68,79},{68,97},{68,137},{69,84},{69,89},{69,102},{69,112},{69,135},{70,103},{70,116},
{70,132},{70,159},{71,128},{71,139},{71,156},{72,77},{72,114},{72,118},{72,141},{73,76},{73,89},
{73,99},{73,109},{73,142},{74,83},{74,109},{74,143},{74,149},{75,78},{75,110},{75,140},{75,154},
{75,166},{76,107},{76,126},{76,135},{76,144},{76,151},{76,153},{77,88},{77,108},{77,116},{77,128},
{77,131},{78,104},{78,127},{78,158},{79,93},{79,130},{79,144},{80,100},{80,108},{80,129},{80,142},
{81,116},{81,130},{81,136},{81,145},{82,102},{82,107},{82,147},{82,150},{83,97},{83,106},{83,148},
{83,163},{84,100},{84,149},{84,160},{85,95},{85,97},{85,145},{85,153},{86,88},{86,114},{86,133},
{86,155},{87,104},{87,119},{87,146},{87,153},{88,110},{88,149},{88,153},{89,98},{89,155},{89,162},
{90,126},{90,152},{90,155},{91,151},{91,156},{92,101},{92,136},{92,153},{92,158},{93,98},{93,115},
{93,148},{93,149},{93,154},{93,159},{94,160},{94,165},{95,96},{95,161},{96,138},{96,141},{96,152},
{97,120},{97,162},{97,164},{98,105},{98,122},{98,139},{98,157},{98,163},{99,101},{99,159},{100,104},
{100,114},{100,120},{100,138},{101,128},{101,140},{102,154},{103,110},{103,133},{103,138},{103,167},
{104,113},{104,130},{104,134},{105,108},{105,167},{106,107},{106,128},{106,130},{106,168},{107,111},
{107,159},{108,164},{109,117},{109,166},{110,143},{110,146},{111,114},{111,162},{112,116},{112,121},
{112,148},{113,143},{113,165},{114,157},{115,121},{115,127},{115,164},{116,151},{116,161},{117,124},
{117,132},{117,154},{118,120},{118,146},{118,168},{119,128},{119,133},{119,154},{120,151},{120,154},
{121,126},{121,137},{121,168},{122,150},{122,168},{123,139},{123,161},{124,125},{124,160},{125,131},
{125,150},{125,157},{126,129},{126,141},{126,165},{126,167},{127,144},{128,147},{129,163},{130,155},
{131,152},{131,162},{132,136},{132,141},{133,134},{133,164},{134,139},{134,141},{135,158},{135,164},
{136,143},{137,138},{137,163},{138,150},{139,143},{139,151},{140,141},{140,148},{140,150},{142,146},
{142,148},{145,157},{145,166},{146,147},{147,157},{147,165},{149,161},{149,165},{150,153},{151,160},
{152,154},{152,158},{156,165},{156,166},{157,158},{162,167}
\end{lstlisting}

\subsection*{(8,2)= 48}

\noindent SG48x5 

\begin{lstlisting}
{1,2},{1,5},{1,7},{1,15},{1,18},{1,23},{1,38},{1,42},{2,3},{2,9},{2,11},{2,16},{2,22},{2,31},{2,46},
{3,12},{3,19},{3,20},{3,21},{3,26},{3,35},{3,47},{4,8},{4,14},{4,16},{4,17},{4,18},{4,19},{4,36},{4,46},
{5,9},{5,10},{5,19},{5,24},{5,29},{5,30},{5,32},{6,9},{6,14},{6,16},{6,23},{6,24},{6,30},{6,35},{6,44},
{7,12},{7,13},{7,14},{7,19},{7,28},{7,33},{7,48},{8,11},{8,12},{8,21},{8,23},{8,24},{8,29},{8,42},
{9,12},{9,17},{9,25},{9,36},{9,37},{10,16},{10,17},{10,21},{10,23},{10,28},{10,37},{10,40},{11,25},
{11,27},{11,28},{11,32},{11,33},{11,44},{12,34},{12,39},{12,40},{12,45},{13,20},{13,24},{13,25},{13,27},
{13,40},{13,45},{13,46},{14,21},{14,22},{14,31},{14,32},{14,37},{15,21},{15,29},{15,31},{15,36},{15,37},
{15,44},{15,45},{16,38},{16,42},{16,45},{16,48},{17,20},{17,29},{17,31},{17,33},{17,38},{18,25},{18,26},
{18,32},{18,34},{18,35},{18,40},{19,27},{19,41},{19,43},{19,44},{20,32},{20,34},{20,41},{20,42},{20,44},
{21,25},{21,30},{21,38},{22,29},{22,30},{22,36},{22,38},{22,40},{22,41},{23,31},{23,41},{23,46},{23,47},
{24,26},{24,31},{24,36},{24,38},{25,41},{25,43},{25,48},{26,33},{26,37},{26,39},{26,41},{26,48},{27,34},
{27,35},{27,37},{27,38},{27,47},{28,34},{28,35},{28,36},{28,39},{28,41},{29,35},{29,46},{29,48},{30,33},
{30,34},{30,42},{30,46},{31,34},{31,43},{32,39},{32,45},{32,47},{33,43},{33,45},{33,47},{34,48},{35,43},
{35,45},{36,42},{36,47},{37,42},{37,46},{38,39},{39,43},{39,44},{39,46},{40,43},{40,44},{40,47},{41,45},
{42,43},{44,48},{47,48}
\end{lstlisting}

\subsection*{(9,2)= 60}

\noindent SG60x9 

\begin{lstlisting}
{1,2},{1,3},{1,7},{1,8},{1,19},{1,28},{1,47},{1,48},{1,57},{2,20},{2,22},{2,25},{2,27},{2,30},{2,31},
{2,37},{2,53},{3,4},{3,10},{3,15},{3,35},{3,40},{3,42},{3,50},{3,60},{4,5},{4,9},{4,17},{4,22},{4,33},
{4,46},{4,51},{4,52},{5,8},{5,21},{5,26},{5,36},{5,40},{5,44},{5,53},{5,54},{6,7},{6,11},{6,19},{6,20},
{6,31},{6,44},{6,50},{6,52},{6,60},{7,12},{7,22},{7,41},{7,42},{7,54},{7,58},{7,59},{8,9},{8,11},{8,16},
{8,25},{8,32},{8,34},{8,38},{9,12},{9,13},{9,19},{9,24},{9,37},{9,39},{9,43},{10,11},{10,20},{10,24},
{10,29},{10,45},{10,54},{10,55},{10,57},{11,12},{11,17},{11,23},{11,27},{11,36},{11,51},{12,15},{12,26},
{12,33},{12,49},{12,53},{12,55},{13,14},{13,15},{13,17},{13,28},{13,31},{13,36},{13,38},{13,54},{14,20},
{14,22},{14,34},{14,35},{14,40},{14,47},{14,49},{14,51},{15,16},{15,23},{15,25},{15,44},{15,45},{15,46},
{16,17},{16,20},{16,26},{16,39},{16,52},{16,56},{16,58},{17,18},{17,30},{17,48},{17,53},{17,59},{18,19},
{18,21},{18,25},{18,32},{18,33},{18,35},{18,47},{18,54},{19,29},{19,30},{19,40},{19,45},{19,56},{20,21},
{20,33},{20,41},{20,43},{21,23},{21,28},{21,37},{21,42},{21,46},{21,49},{22,23},{22,29},{22,32},{22,44},
{22,56},{23,24},{23,38},{23,39},{23,40},{23,48},{24,25},{24,26},{24,30},{24,47},{24,58},{24,60},{25,50},
{25,51},{25,55},{25,59},{26,29},{26,31},{26,35},{26,48},{26,51},{27,35},{27,43},{27,46},{27,49},{27,54},
{27,56},{27,58},{28,29},{28,44},{28,51},{28,55},{28,58},{28,60},{29,38},{29,46},{29,50},{29,59},{30,34},
{30,41},{30,42},{30,44},{30,49},{31,32},{31,40},{31,42},{31,46},{31,57},{32,41},{32,43},{32,45},{32,49},
{32,60},{33,34},{33,36},{33,38},{33,57},{33,58},{34,37},{34,46},{34,48},{34,54},{34,60},{35,37},{35,38},
{35,41},{35,44},{36,37},{36,41},{36,47},{36,50},{36,56},{37,45},{37,52},{37,59},{38,42},{38,52},{38,53},
{39,41},{39,50},{39,51},{39,53},{39,54},{39,57},{40,55},{40,58},{40,59},{41,46},{41,55},{42,43},{42,51},
{42,56},{43,44},{43,47},{43,48},{43,59},{44,57},{45,48},{45,51},{45,53},{45,58},{46,47},{47,52},{47,53},
{48,50},{48,55},{49,50},{49,52},{49,57},{50,58},{52,54},{52,55},{53,60},{55,56},{56,57},{56,60},{57,59},
{59,60}
\end{lstlisting}

\subsection*{(10,2)= 72}

\noindent SG72x28x1 

\begin{lstlisting}
{1,2},{1,4},{1,9},{1,22},{1,30},{1,35},{1,37},{1,39},{1,50},{1,65},{2,5},{2,6},{2,16}{2,26},{2,32},
{2,38},{2,46},{2,47},{2,56},{3,4},{3,16},{3,25},{3,32},{3,35},{3,49},{3,53},{3,59},{3,61},{3,71},{4,5},
{4,20},{4,24},{4,27},{4,40},{4,55},{4,57},{4,63},{5,15},{5,17},{5,23},{5,25},{5,33},{5,43},{5,44},
{5,67},{6,7},{6,19},{6,31},{6,36},{6,39},{6,40},{6,45},{6,53},{6,68},{7,9},{7,11},{7,18},{7,33},{7,42},
{7,49},{7,55},{7,57},{7,67},{8,9},{8,17},{8,23},{8,33},{8,36},{8,45},{8,46},{8,52},{8,53},{8,63},{9,16},
{9,20},{9,25},{9,44},{9,62},{9,69},{9,70},{10,12},{10,13},{10,25},{10,39},{10,47},{10,52},{10,53},
{10,57},{10,58},{10,60},{11,12},{11,16},{11,17},{11,21},{11,24},{11,39},{11,43},{11,58},{11,63},{12,35},
{12,36},{12,38},{12,40},{12,42},{12,44},{12,46},{12,65},{13,14},{13,15},{13,19},{13,30},{13,32},{13,41},
{13,49},{13,63},{13,70},{14,20},{14,21},{14,37},{14,38},{14,51},{14,53},{14,61},{14,66},{14,67},{15,16},
{15,18},{15,20},{15,26},{15,34},{15,36},{15,50},{15,71},{16,27},{16,52},{16,64},{16,66},{16,72},{17,19},
{17,30},{17,55},{17,60},{17,61},{17,69},{17,70},{18,28},{18,37},{18,40},{18,46},{18,48},{18,53},{18,58},
{18,61},{19,22},{19,27},{19,28},{19,35},{19,51},{19,62},{19,65},{20,29},{20,35},{20,47},{20,58},{20,60},
{20,68},{21,22},{21,25},{21,34},{21,36},{21,40},{21,47},{21,48},{21,49},{22,31},{22,42},{22,43},{22,46},
{22,52},{22,58},{22,71},{23,39},{23,40},{23,42},{23,49},{23,58},{23,62},{23,66},{23,72},{24,26},{24,33},
{24,37},{24,50},{24,52},{24,62},{24,68},{24,70},{25,26},{25,28},{25,37},{25,41},{25,68},{26,31},{26,35},
{26,42},{26,45},{26,51},{26,69},{27,37},{27,42},{27,44},{27,45},{27,47},{27,54},{27,71},{28,29},{28,38},
{28,50},{28,59},{28,60},{28,63},{28,66},{29,32},{29,36},{29,39},{29,42},{29,43},{29,49},{29,52},{29,61},
{30,31},{30,38},{30,42},{30,48},{30,59},{30,64},{30,68},{31,44},{31,54},{31,60},{31,61},{31,63},{31,66},
{32,33},{32,34},{32,37},{32,44},{32,58},{32,69},{33,35},{33,47},{33,64},{33,66},{33,71},{34,53},{34,54},
{34,55},{34,59},{34,60},{34,62},{34,65},{35,48},{35,54},{35,72},{36,37},{36,41},{36,57},{36,59},{37,60},
{37,72},{38,45},{38,54},{38,55},{38,62},{38,71},{39,48},{39,51},{39,54},{39,71},{40,52},{40,64},{40,69},
{40,70},{41,50},{41,54},{41,55},{41,56},{41,58},{41,61},{41,64},{42,53},{42,56},{43,45},{43,53},{43,57},
{43,64},{43,70},{43,72},{44,48},{44,50},{44,51},{44,68},{45,48},{45,49},{45,58},{45,65},{46,49},{46,51},
{46,54},{46,57},{46,68},{47,50},{47,59},{47,61},{47,72},{48,52},{48,56},{48,57},{49,50},{49,60},{50,53},
{50,69},{51,55},{51,56},{51,59},{51,64},{52,55},{52,67},{54,67},{54,70},{55,66},{55,72},{56,60},{56,63},
{56,65},{56,70},{56,71},{57,62},{57,66},{57,69},{58,59},{59,67},{59,70},{60,64},{61,62},{61,65},{62,63},
{62,64},{63,67},{63,72},{64,67},{65,66},{65,67},{65,68},{66,70},{67,69},{68,71},{68,72},{69,71},{69,72}
\end{lstlisting}

\noindent SG72x28x3 

\begin{lstlisting}
{1,2},{1,3},{1,4},{1,10},{1,31},{1,32},{1,34},{1,43},{1,49},{1,50},{2,5},{2,9},{2,22},{2,37},{2,38},
{2,39},{2,51},{2,63},{2,69},{3,22},{3,23},{3,27},{3,42},{3,54},{3,55},{3,57},{3,58},{3,70},{4,5},{4,6},
{4,16},{4,18},{4,23},{4,29},{4,40},{4,48},{4,67},{5,25},{5,26},{5,27},{5,42},{5,47},{5,52},{5,64},
{5,68},{6,7},{6,9},{6,22},{6,34},{6,42},{6,44},{6,56},{6,62},{6,71},{7,24},{7,31},{7,33},{7,39},{7,43},
{7,57},{7,60},{7,68},{7,69},{8,9},{8,10},{8,12},{8,16},{8,25},{8,27},{8,30},{8,45},{8,46},{8,68},{9,26},
{9,32},{9,35},{9,40},{9,61},{9,65},{9,70},{10,11},{10,24},{10,28},{10,39},{10,41},{10,47},{10,53},
{10,56},{11,26},{11,38},{11,40},{11,42},{11,43},{11,63},{11,67},{11,70},{11,71},{12,13},{12,15},{12,22},
{12,41},{12,43},{12,48},{12,54},{12,59},{12,64},{13,21},{13,26},{13,28},{13,29},{13,34},{13,39},{13,45},
{13,58},{13,66},{14,15},{14,17},{14,23},{14,34},{14,38},{14,53},{14,54},{14,57},{14,65},{14,68},{15,18},
{15,21},{15,26},{15,31},{15,33},{15,37},{15,55},{15,56},{16,17},{16,33},{16,38},{16,52},{16,57},{16,60},
{16,66},{16,72},{17,19},{17,35},{17,39},{17,42},{17,43},{17,49},{17,50},{17,56},{18,24},{18,35},{18,46},
{18,58},{18,61},{18,63},{18,68},{18,72},{19,20},{19,22},{19,28},{19,31},{19,33},{19,40},{19,58},{19,59},
{19,68},{20,26},{20,30},{20,34},{20,41},{20,46},{20,48},{20,51},{20,55},{20,57},{21,22},{21,24},{21,30},
{21,36},{21,40},{21,49},{21,52},{21,57},{22,26},{22,53},{22,67},{22,72},{23,26},{23,28},{23,30},{23,32},
{23,59},{23,69},{23,71},{24,38},{24,42},{24,50},{24,51},{24,59},{24,65},{25,28},{25,29},{25,36},{25,43},
{25,55},{25,62},{25,65},{25,72},{26,50},{26,60},{27,33},{27,34},{27,36},{27,50},{27,61},{27,67},{27,71},
{28,37},{28,44},{28,54},{28,57},{28,61},{29,33},{29,35},{29,41},{29,46},{29,51},{29,53},{29,70},{30,37},
{30,43},{30,47},{30,53},{30,58},{30,62},{31,35},{31,38},{31,45},{31,52},{31,62},{31,71},{32,33},{32,36},
{32,42},{32,46},{32,64},{32,66},{32,72},{33,47},{33,63},{33,65},{34,37},{34,40},{34,59},{34,72},{35,36},
{35,47},{35,57},{35,59},{35,67},{36,38},{36,48},{36,56},{36,60},{36,68},{37,42},{37,60},{37,64},{37,67},
{37,70},{38,41},{38,44},{38,58},{39,46},{39,48},{39,55},{39,62},{39,71},{40,54},{40,55},{40,64},{40,69},
{41,42},{41,49},{41,61},{41,62},{41,69},{42,45},{43,44},{43,51},{43,61},{44,45},{44,46},{44,47},{44,49},
{44,50},{44,55},{45,48},{45,51},{45,53},{45,63},{45,69},{46,52},{46,54},{46,67},{47,48},{47,54},{47,66},
{47,72},{48,61},{48,65},{48,70},{49,60},{49,63},{49,65},{49,68},{49,71},{50,62},{50,64},{50,69},{50,70},
{51,54},{51,56},{51,66},{51,71},{52,56},{52,59},{52,61},{52,65},{52,70},{53,55},{53,60},{53,61},{53,64},
{54,60},{54,62},{55,59},{55,66},{56,58},{56,63},{56,69},{57,63},{57,64},{58,60},{58,64},{58,65},{59,60},
{59,63},{61,66},{62,63},{62,66},{64,71},{65,67},{66,67},{66,68},{67,69},{68,70},{69,72},{70,72},{71,72}
\end{lstlisting}

\noindent SG72x30x1 

\begin{lstlisting}
{1,2},{1,9},{1,10},{1,13},{1,29},{1,45},{1,49},{1,54},{1,64},{1,70},{2,3},{2,6},{2,11},{2,14},{2,17},
{2,34},{2,48},{2,60},{2,64},{3,4},{3,17},{3,35},{3,41},{3,44},{3,46},{3,51},{3,65},{3,72},{4,5},{4,20},
{4,24},{4,27},{4,29},{4,37},{4,39},{4,41},{4,43},{5,6},{5,8},{5,18},{5,24},{5,45},{5,47},{5,62},{5,66},
{5,70},{6,11},{6,22},{6,32},{6,33},{6,36},{6,46},{6,62},{6,63},{7,8},{7,13},{7,20},{7,21},{7,30},{7,31},
{7,36},{7,46},{7,50},{7,60},{8,10},{8,17},{8,21},{8,34},{8,42},{8,56},{8,70},{8,71},{9,10},{9,12},
{9,15},{9,31},{9,33},{9,41},{9,43},{9,55},{9,62},{10,20},{10,45},{10,51},{10,56},{10,61},{10,63},
{10,67},{11,13},{11,33},{11,37},{11,47},{11,51},{11,56},{11,59},{11,69},{12,13},{12,14},{12,15},{12,22},
{12,24},{12,42},{12,57},{12,58},{12,65},{13,18},{13,36},{13,37},{13,49},{13,54},{13,65},{14,17},{14,19},
{14,20},{14,31},{14,40},{14,58},{14,66},{14,69},{15,16},{15,17},{15,24},{15,36},{15,47},{15,53},{15,60},
{15,68},{16,23},{16,28},{16,29},{16,34},{16,36},{16,56},{16,65},{16,66},{16,72},{17,25},{17,26},{17,53},
{17,54},{17,71},{18,19},{18,26},{18,43},{18,45},{18,48},{18,53},{18,55},{18,65},{19,20},{19,33},{19,40},
{19,55},{19,64},{19,68},{19,71},{19,72},{20,29},{20,36},{20,38},{20,59},{20,67},{21,22},{21,26},{21,34},
{21,38},{21,40},{21,41},{21,47},{21,49},{22,26},{22,29},{22,46},{22,58},{22,61},{22,68},{22,71},{23,24},
{23,29},{23,31},{23,38},{23,40},{23,46},{23,48},{23,54},{23,56},{24,27},{24,38},{24,50},{24,51},{24,64},
{25,26},{25,30},{25,32},{25,33},{25,37},{25,38},{25,45},{25,57},{25,66},{26,29},{26,32},{26,51},{26,55},
{26,69},{27,28},{27,49},{27,55},{27,57},{27,60},{27,63},{27,69},{27,71},{28,36},{28,41},{28,45},{28,48},
{28,50},{28,58},{28,69},{28,71},{29,30},{29,59},{29,70},{30,33},{30,40},{30,50},{30,60},{30,63},{30,65},
{30,70},{31,32},{31,43},{31,44},{31,46},{31,66},{31,69},{32,36},{32,39},{32,42},{32,44},{32,63},{32,64},
{33,38},{33,43},{33,71},{33,72},{34,43},{34,52},{34,57},{34,64},{34,65},{34,67},{35,36},{35,38},{35,44},
{35,49},{35,52},{35,55},{35,56},{35,58},{35,70},{36,52},{37,41},{37,56},{37,61},{37,64},{37,66},{37,68},
{38,48},{38,49},{38,67},{39,42},{39,43},{39,47},{39,54},{39,58},{39,60},{39,67},{39,72},{40,41},{40,47},
{40,52},{40,54},{40,63},{41,42},{41,48},{41,62},{42,48},{42,57},{42,59},{42,70},{42,72},{43,52},{43,53},
{43,58},{44,47},{44,59},{44,61},{44,64},{44,65},{44,71},{45,46},{45,47},{45,57},{45,58},{46,57},{46,68},
{46,72},{47,59},{47,60},{48,53},{48,60},{48,61},{49,53},{49,63},{49,66},{49,72},{50,51},{50,53},{50,58},
{50,59},{50,62},{50,64},{51,52},{51,61},{51,69},{51,72},{52,54},{52,57},{52,61},{52,70},{53,59},{53,63},
{53,68},{54,62},{54,67},{54,71},{55,56},{55,57},{55,60},{55,62},{56,58},{56,63},{57,59},{59,62},{60,61},
{61,66},{61,71},{62,66},{62,67},{63,65},{64,68},{65,67},{66,72},{67,68},{67,69},{68,69},{68,70},{69,70}
\end{lstlisting}

\noindent SG72x30x2 

\begin{lstlisting}
{1,4},{1,6},{1,8},{1,13},{1,32},{1,35},{1,39},{1,46},{1,57},{1,71},{2,3},{2,7},{2,13},{2,29},{2,36},
{2,38},{2,48},{2,55},{2,64},{2,68},{3,11},{3,19},{3,31},{3,43},{3,44},{3,45},{3,56},{3,68},{3,71},{4,5},
{4,12},{4,37},{4,38},{4,40},{4,44},{4,57},{4,62},{4,70},{5,14},{5,25},{5,36},{5,58},{5,59},{5,61},
{5,68},{5,70},{5,72},{6,7},{6,14},{6,16},{6,34},{6,41},{6,49},{6,53},{6,65},{6,71},{7,10},{7,12},{7,27},
{7,28},{7,41},{7,47},{7,59},{7,64},{8,9},{8,17},{8,18},{8,28},{8,29},{8,32},{8,52},{8,68},{8,72},{9,22},
{9,26},{9,34},{9,44},{9,52},{9,59},{9,60},{9,64},{9,66},{10,11},{10,12},{10,17},{10,20},{10,37},{10,39},
{10,43},{10,59},{10,69},{11,26},{11,31},{11,49},{11,52},{11,54},{11,57},{11,59},{11,63},{12,23},{12,28},
{12,30},{12,44},{12,54},{12,55},{12,67},{13,35},{13,38},{13,39},{13,49},{13,51},{13,61},{13,66},{13,67},
{14,23},{14,29},{14,31},{14,33},{14,34},{14,43},{14,58},{14,66},{15,16},{15,19},{15,21},{15,25},{15,29},
{15,35},{15,37},{15,44},{15,47},{15,49},{16,21},{16,22},{16,23},{16,46},{16,49},{16,62},{16,68},{16,69},
{17,18},{17,24},{17,25},{17,34},{17,38},{17,45},{17,69},{17,72},{18,27},{18,30},{18,33},{18,35},{18,57},
{18,61},{18,65},{18,68},{19,27},{19,32},{19,34},{19,36},{19,37},{19,54},{19,56},{19,67},{20,21},{20,26},
{20,30},{20,32},{20,41},{20,43},{20,61},{20,64},{20,70},{21,27},{21,32},{21,38},{21,44},{21,53},{21,58},
{21,63},{22,24},{22,27},{22,43},{22,48},{22,59},{22,62},{22,67},{22,71},{23,24},{23,26},{23,32},{23,33},
{23,50},{23,67},{23,68},{24,32},{24,33},{24,45},{24,47},{24,48},{24,49},{24,70},{25,26},{25,47},{25,51},
{25,55},{25,61},{25,71},{25,72},{26,34},{26,50},{26,51},{26,57},{26,64},{27,47},{27,56},{27,57},{27,58},
{27,67},{28,29},{28,51},{28,54},{28,56},{28,62},{28,63},{28,70},{29,35},{29,42},{29,43},{29,48},{29,57},
{30,31},{30,34},{30,35},{30,55},{30,60},{30,62},{30,70},{31,32},{31,47},{31,53},{31,62},{31,66},{31,72},
{32,36},{32,47},{33,37},{33,63},{33,64},{33,65},{33,66},{33,71},{34,38},{34,54},{34,60},{35,45},{35,46},
{35,50},{35,59},{36,42},{36,51},{36,55},{36,59},{36,65},{36,69},{37,39},{37,41},{37,52},{37,62},{37,64},
{38,40},{38,50},{38,62},{38,63},{39,48},{39,51},{39,53},{39,58},{39,60},{39,68},{40,43},{40,46},{40,47},
{40,50},{40,52},{40,54},{40,60},{40,65},{41,42},{41,45},{41,50},{41,52},{41,61},{41,62},{42,45},{42,57},
{42,60},{42,63},{42,67},{42,69},{42,72},{43,65},{43,67},{43,71},{44,45},{44,53},{44,65},{44,66},{45,46},
{45,51},{45,58},{46,54},{46,58},{46,64},{46,69},{46,72},{47,60},{47,66},{48,53},{48,54},{48,55},{48,57},
{48,61},{49,52},{49,65},{49,67},{49,70},{50,52},{50,53},{50,56},{50,59},{51,58},{51,62},{51,65},{52,55},
{52,58},{53,56},{53,65},{53,72},{54,61},{54,63},{55,58},{55,69},{55,71},{56,66},{56,69},{56,70},{56,71},
{57,69},{59,63},{60,63},{60,68},{60,70},{61,66},{61,68},{63,71},{64,70},{64,72},{66,69},{67,72}
\end{lstlisting}

\subsection*{(11,2)= 84}

\noindent SG84x10x3

\begin{lstlisting}
{1,2},{1,13},{1,18},{1,19},{1,26},{1,45},{1,47},{1,66},{1,72},{1,74},{1,75},{2,3},{2,12},{2,23},{2,30},
{2,34},{2,49},{2,54},{2,57},{2,80},{2,81},{3,4},{3,32},{3,38},{3,48},{3,50},{3,52},{3,59},{3,62},{3,74},
{3,83},{4,9},{4,10},{4,16},{4,18},{4,22},{4,23},{4,41},{4,46},{4,56},{4,66},{5,15},{5,21},{5,24},{5,25},
{5,43},{5,48},{5,66},{5,70},{5,75},{5,77},{5,81},{6,7},{6,16},{6,26},{6,27},{6,35},{6,37},{6,50},{6,52},
{6,73},{6,81},{6,84},{7,9},{7,12},{7,20},{7,32},{7,38},{7,44},{7,62},(7,75},{7,77},{7,83},{8,9},{8,14},
{8,45},{8,48},{8,52},{8,53},{8,54},{8,62},{8,73},{8,78},{8,80},{9,24},{9,26},{9,30},{9,35},{9,49},
{9,61},{9,67},{9,74},{10,11},{10,20},{10,27},{10,29},{10,40},{10,64},{10,70},{10,75},{10,76},{10,80},
{11,14},{11,17},{11,19},{11,26},{11,31},{11,38},{11,41},{11,46},{11,69},{11,81},{12,13},{12,19},{12,43},
{12,53},{12,56},{12,70},{12,72},{12,76},{12,82},{13,16},{13,21},{13,28},{13,29},{13,40},{13,54},{13,59},
{13,67},{13,81},{14,18},{14,23},{14,32},{14,37},{14,43},{14,49},{14,58},{14,59},{14,63},{15,16},{15,19},
{15,27},{15,40},{15,47},{15,49},{15,57},{15,62},{15,65},{15,72},{16,17},{16,36},{16,63},{16,64},{16,68},
{16,69},{16,80},{17,18},{17,20},{17,28},{17,43},{17,51},{17,52},{17,57},{17,58},{17,61},{18,30},{18,42},
{18,44},{18,50},{18,65},{18,70},{18,71},{19,22},{19,24},{19,25},{19,34},{19,36},{19,52},{19,84},{20,21},
{20,22},{20,23},{20,31},{20,46},{20,48},{20,63},{20,72},{21,34},{21,41},{21,49},{21,53},{21,60},{21,65},
{21,83},{21,84},{22,33},{22,37},{22,45},{22,51},{22,58},{22,68},{22,71},{22,81},{23,24},{23,40},{23,64},
{23,68},{23,69},{23,73},{23,79},{24,28},{24,47},{24,50},{24,55},{24,59},{24,65},{24,80},{25,28},{25,32},
{25,42},{25,49},{25,54},{25,56},{25,64},{25,72},{25,73},{26,33},{26,48},{26,56},{26,62},{26,71},{26,79},
{26,83},{27,28},{27,30},{27,43},{27,45},{27,67},{27,79},{27,82},{27,83},{28,34},{28,37},{28,39},{28,62},
{28,66},{28,76},{29,33},{29,35},{29,36},{29,43},{29,50},{29,57},{29,62},{29,73},{29,77},{30,31},{30,33},
{30,36},{30,59},{30,60},{30,73},{30,77},{31,40},{31,50},{31,56},{31,58},{31,66},{31,68},{31,74},{31,80},
{32,33},{32,47},{32,60},{32,61},{32,67},{32,79},{32,80},{33,52},{33,65},{33,66},{33,69},{33,78},{33,82},
{34,35},{34,44},{34,47},{34,56},{34,63},{34,75},{34,78},{35,38},{35,42},{35,48},{35,58},{35,60},{35,72},
{35,79},{36,37},{36,38},{36,39},{36,48},{36,61},{36,79},{36,83},{37,40},{37,41},{37,46},{37,57},{37,70},
{37,74},{38,45},{38,65},{38,66},{38,67},{38,71},{38,73},{39,47},{39,49},{39,52},{39,56},{39,69},{39,71},
{39,72},{39,75},{39,81},{40,42},{40,44},{40,52},{40,61},{40,71},{41,44},{41,58},{41,68},{41,72},{41,77},
{41,80},{41,82},{42,46},{42,53},{42,63},{42,68},{42,69},{42,80},{42,83},{43,44},{43,46},{43,68},{43,71},
{43,74},{44,48},{44,51},{44,55},{44,64},{44,69},{45,50},{45,56},{45,60},{45,61},{45,69},{45,77},{46,47},
{46,51},{46,54},{46,60},{46,78},{47,53},{47,58},{47,73},{47,76},{48,67},{48,76},{48,82},{49,50},{49,51},
{49,66},{49,76},{50,78},{50,79},{50,82},{51,59},{51,70},{51,73},{51,74},{51,79},{51,80},{52,60},{52,67},
{52,77},{53,57},{53,59},{53,64},{53,66},{53,79},{53,81},{54,55},{54,58},{54,65},{54,75},{54,79},{54,84},
{55,57},{55,60},{55,63},{55,66},{55,72},{55,76},{55,81},{55,83},{56,57},{56,59},{56,65},{57,67},{57,75},
{57,84},{58,64},{58,70},{58,83},{59,72},{59,75},{59,84},{60,62},{60,64},{60,82},{61,62},{61,78},{61,81},
{61,82},{61,84},{62,68},{62,70},{63,70},{63,71},{63,74},{63,77},{63,82},{64,71},{64,74},{64,78},{65,74},
{65,76},{65,81},{66,84},{67,68},{67,70},{67,78},{68,75},{68,76},{69,70},{69,74},{69,84},{71,80},{71,84},
{72,78},{73,82},{73,83},{75,82},{76,77},{76,84},{77,78},{77,79},{78,83}
\end{lstlisting}

\noindent SG84x10x5

\begin{lstlisting}
{1,2},{1,11},{1,20},{1,35},{1,36},{1,38},{1,47},{1,48},{1,54},{1,59},{1,69},{2,12},{2,17},{2,39},{2,42},
{2,51},{2,53},{2,57},{2,61},{2,67},{2,82},{3,4},{3,9},{3,12},{3,15},{3,30},{3,34},{3,36},{3,54},{3,61},
{3,70},{3,84},{4,17},{4,21},{4,40},{4,44},{4,47},{4,52},{4,57},{4,69},{4,72},{4,73},{5,6},{5,8},{5,27},
{5,47},{5,48},{5,52},{5,59},{5,66},{5,68},{5,70},{5,82},{6,17},{6,32},{6,35},{6,45},{6,50},{6,51},
{6,53},{6,54},{6,62},{6,78},{7,10},{7,16},{7,18},{7,33},{7,42},{7,45},{7,54},{7,60},{7,72},{7,79},
{7,82},{8,21},{8,24},{8,31},{8,35},{8,36},{8,39},{8,55},{8,61},{8,72},{8,77},{9,10},{9,19},{9,39},
{9,48},{9,57},{9,59},{9,66},{9,68},{9,77},{9,78},{10,13},{10,17},{10,22},{10,23},{10,26},{10,32},
{10,47},{10,61},{10,83},{11,22},{11,30},{11,52},{11,62},{11,68},{11,72},{11,73},{11,77},{11,79},{11,80},
{12,13},{12,24},{12,28},{12,46},{12,47},{12,50},{12,60},{12,64},{12,77},{13,14},{13,35},{13,43},{13,44},
{13,49},{13,58},{13,66},{13,73},{13,84},{14,15},{14,16},{14,28},{14,36},{14,40},{14,46},{14,53},{14,68},
{14,71},{14,79},{15,22},{15,27},{15,35},{15,37},{15,42},{15,58},{15,64},{15,80},{15,81},{16,17},{16,24},
{16,25},{16,30},{16,51},{16,56},{16,59},{16,62},{16,67},{17,37},{17,63},{17,74},{17,75},{17,77},{17,84},
{18,27},{18,28},{18,32},{18,44},{18,48},{18,51},{18,61},{18,74},{18,75},{18,80},{19,24},{19,25},{19,37},
{19,44},{19,45},{19,47},{19,51},{19,53},{19,62},{19,84},{20,21},{20,25},{20,27},{20,32},{20,40},{20,70},
{20,72},{20,73},{20,77},{20,78},{21,26},{21,33},{21,51},{21,64},{21,66},{21,68},{21,76},{21,80},{22,24},
{22,40},{22,41},{22,51},{22,67},{22,74},{22,76},{22,82},{23,24},{23,29},{23,34},{23,38},{23,42},{23,51},
{23,53},{23,73},{23,75},{23,82},{24,32},{24,54},{24,57},{24,63},{24,69},{25,35},{25,50},{25,52},{25,61},
{25,63},{25,74},{25,75},{25,82},{26,27},{26,28},{26,41},{26,45},{26,54},{26,62},{26,63},{26,81},{26,82},
{27,30},{27,34},{27,43},{27,50},{27,57},{27,84},{28,29},{28,35},{28,37},{28,52},{28,70},{28,76},{28,78},
{29,50},{29,58},{29,59},{29,60},{29,64},{29,67},{29,71},{29,72},{29,84},{30,31},{30,53},{30,58},{30,60},
{30,65},{30,76},{30,83},{31,40},{31,42},{31,47},{31,59},{31,66},{31,74},{31,78},{31,79},{31,81},{32,36},
{32,41},{32,42},{32,52},{32,65},{32,67},{33,34},{33,40},{33,52},{33,53},{33,55},{33,58},{33,59},{33,69},
{33,77},{34,35},{34,37},{34,46},{34,49},{34,56},{34,67},{34,79},{35,57},{35,60},{35,65},{36,37},{36,43},
{36,49},{36,50},{36,58},{36,82},{37,60},{37,66},{37,71},{37,73},{37,83},{38,40},{38,45},{38,52},{38,56},
{38,66},{38,77},{38,80},{38,83},{38,84},{39,40},{39,50},{39,52},{39,54},{39,56},{39,65},{39,73},{39,80},
{40,48},{40,60},{40,62},{41,44},{41,53},{41,56},{41,59},{41,60},{41,61},{41,75},{41,79},{42,44},{42,62},
{42,63},{42,68},{42,76},{43,46},{43,51},{43,52},{43,60},{43,63},{43,64},{43,71},{43,78},{44,54},{44,67},
{44,77},{44,82},{44,83},{45,46},{45,58},{45,61},{45,67},{45,73},{45,74},{46,59},{46,65},{46,72},{46,75},
{46,76},{46,83},{47,55},{47,56},{47,71},{47,80},{48,56},{48,63},{48,64},{48,73},{48,76},{48,81},{49,57},
{49,60},{49,62},{49,70},{49,74},{49,76},{49,80},{49,81},{50,69},{50,76},{50,79},{50,83},{51,70},{51,81},
{52,81},{53,63},{53,72},{53,74},{54,71},{54,74},{54,75},{55,67},{55,68},{55,73},{55,75},{55,78},{55,81},
{55,83},{55,84},{56,58},{56,68},{56,72},{56,78},{57,58},{57,71},{57,79},{57,83},{58,70},{58,75},{59,73},
{59,80},{60,68},{61,62},{61,69},{61,71},{62,64},{62,75},{63,67},{63,70},{63,80},{64,65},{64,74},{64,79},
{64,83},{65,69},{65,70},{65,71},{65,82},{65,84},{66,67},{66,69},{66,72},{66,75},{68,69},{68,74},{69,78},
{69,81},{70,79},{70,83},{71,76},{71,77},{72,81},{76,84},{77,81},{78,80},{78,82},{79,84}
\end{lstlisting}

\noindent SG84x10x7

\begin{lstlisting}
{1,2},{1,18},{1,25},{1,34},{1,63},{1,66},{1,67},{1,73},{1,79},{1,80},{1,84},{2,6},{2,8},{2,12},{2,13},
{2,30},{2,32},{2,39},{2,58},{2,65},{2,83},{3,4},{3,13},{3,22},{3,25},{3,35},{3,37},{3,38},{3,43},{3,62},
{3,71},{3,81},{4,8},{4,11},{4,19},{4,23},{4,24},{4,45},{4,48},{4,64},{4,69},{4,79},{5,6},{5,10},{5,15},
{5,17},{5,22},{5,26},{5,35},{5,36},{5,44},{5,63},{5,64},{6,7},{6,24},{6,38},{6,40},{6,43},{6,52},{6,71},
{6,76},{6,79},{7,10},{7,26},{7,34},{7,37},{7,41},{7,45},{7,50},{7,57},{7,64},{7,81},{8,9},{8,15},{8,29},
{8,41},{8,53},{8,56},{8,67},{8,76},{8,78},{9,12},{9,14},{9,32},{9,34},{9,38},{9,42},{9,44},{9,65},
{9,72},{9,75},{10,11},{10,33},{10,38},{10,43},{10,48},{10,51},{10,56},{10,58},{10,67},{11,12},{11,15},
{11,16},{11,31},{11,52},{11,70},{11,77},{11,80},{11,82},{12,22},{12,28},{12,41},{12,51},{12,55},{12,63},
{12,69},{12,78},{13,14},{13,15},{13,17},{13,28},{13,34},{13,51},{13,57},{13,60},{13,78},{14,18},{14,31},
{14,40},{14,45},{14,46},{14,63},{14,67},{14,68},{14,73},{15,18},{15,27},{15,42},{15,50},{15,55},{15,73},
{15,83},{16,17},{16,21},{16,24},{16,27},{16,37},{16,39},{16,53},{16,56},{16,63},{16,75},{17,20},{17,41},
{17,59},{17,66},{17,69},{17,70},{17,72},{17,81},{18,19},{18,21},{18,41},{18,43},{18,47},{18,61},{18,65},
{18,78},{19,20},{19,32},{19,34},{19,40},{19,51},{19,52},{19,56},{19,63},{19,84},{20,37},{20,39},{20,55},
{20,58},{20,67},{20,75},{20,76},{20,77},{20,79},{21,23},{21,28},{21,30},{21,33},{21,36},{21,72},{21,76},
{21,77},{21,81},{22,23},{22,27},{22,45},{22,56},{22,60},{22,61},{22,70},{22,79},{23,26},{23,38},{23,43},
{23,49},{23,57},{23,58},{23,59},{23,73},{24,26},{24,35},{24,44},{24,50},{24,51},{24,55},{24,61},{24,67},
{25,26},{25,42},{25,45},{25,52},{25,56},{25,58},{25,64},{25,69},{25,76},{26,32},{26,46},{26,47},{26,54},
{26,77},{26,78},{27,32},{27,34},{27,40},{27,48},{27,65},{27,67},{27,69},{27,84},{28,29},{28,40},{28,47},
{28,54},{28,64},{28,67},{28,73},{28,84},{29,32},{29,34},{29,42},{29,43},{29,59},{29,61},{29,63},{29,77},
{29,83},{30,31},{30,37},{30,42},{30,44},{30,46},{30,48},{30,52},{30,66},{30,70},{31,36},{31,54},{31,59},
{31,61},{31,75},{31,76},{31,81},{31,84},{32,33},{32,62},{32,73},{32,81},{32,82},{33,35},{33,45},{33,50},
{33,59},{33,64},{33,71},{33,78},{33,79},{34,35},{34,36},{34,49},{34,70},{35,41},{35,46},{35,54},{35,58},
{35,74},{35,82},{36,37},{36,39},{36,53},{36,55},{36,60},{36,62},{36,69},{37,59},{37,72},{37,73},{37,78},
{37,82},{38,39},{38,46},{38,50},{38,61},{38,66},{38,84},{39,42},{39,45},{39,47},{39,68},{39,74},{39,80},
{40,42},{40,58},{40,62},{40,66},{40,78},{40,80},{41,42},{41,60},{41,62},{41,73},{41,84},{42,49},{42,51},
{42,79},{43,44},{43,54},{43,68},{43,69},{43,75},{44,45},{44,57},{44,58},{44,71},{44,80},{44,84},{45,54},
{45,66},{45,83},{46,56},{46,64},{46,69},{46,79},{46,80},{46,83},{47,48},{47,53},{47,59},{47,70},{47,71},
{47,75},{47,82},{48,55},{48,60},{48,63},{48,72},{48,74},{48,76},{49,63},{49,65},{49,67},{49,71},{49,78},
{49,81},{49,82},{49,83},{50,54},{50,58},{50,60},{50,63},{50,69},{50,70},{51,53},{51,65},{51,73},{51,76},
{51,83},{52,53},{52,59},{52,60},{52,67},{52,72},{52,74},{53,58},{53,66},{53,68},{53,79},{53,81},{54,56},
{54,65},{54,72},{54,79},{55,56},{55,59},{55,66},{55,68},{55,81},{56,57},{56,71},{57,61},{57,62},{57,69},
{57,74},{57,79},{57,82},{58,61},{59,65},{59,80},{60,75},{60,77},{60,80},{60,82},{61,64},{61,71},{61,72},
{62,63},{62,65},{62,67},{62,70},{62,83},{64,65},{64,68},{64,75},{65,77},{66,75},{66,77},{66,82},{68,70},
{68,72},{68,77},{68,82},{68,84},{69,71},{70,74},{70,76},{71,73},{71,77},{72,80},{72,83},{73,74},{74,75},
{74,77},{74,78},{74,81},{75,83},{76,80},{76,82},{78,84},{80,81},{83,84}
\end{lstlisting}

\subsection*{(12,2)= 96}

\noindent SG96x4x1  

\begin{lstlisting}
{1,4},{1,8},{1,12},{1,49},{1,50},{1,54},{1,62},{1,71},{1,74},{1,80},{1,89},{1,95},{2,5},{2,9},{2,10},
{2,50},{2,51},{2,52},{2,63},{2,72},{2,75},{2,81},{2,90},{2,96},{3,6},{3,7},{3,11},{3,49},{3,51},{3,53},
{3,61},{3,70},{3,73},{3,79},{3,88},{3,94},{4,9},{4,11},{4,51},{4,52},{4,53},{4,65},{4,68},{4,77},{4,83},
{4,86},{4,92},{5,7},{5,12},{5,49},{5,53},{5,54},{5,66},{5,69},{5,78},{5,84},{5,87},{5,93},{6,8},{6,10},
{6,50},{6,52},{6,54},{6,64},{6,67},{6,76},{6,82},{6,85},{6,91},{7,10},{7,55},{7,56},{7,60},{7,62},
{7,68},{7,77},{7,80},{7,86},{7,95},{8,11},{8,56},{8,57},{8,58},{8,63},{8,69},{8,78},{8,81},{8,87},
{8,96},{9,12},{9,55},{9,57},{9,59},{9,61},{9,67},{9,76},{9,79},{9,85},{9,94},{10,57},{10,58},{10,59},
{10,65},{10,71},{10,74},{10,83},{10,89},{10,92},{11,55},{11,59},{11,60},{11,66},{11,72},{11,75},{11,84},
{11,90},{11,93},{12,56},{12,58},{12,60},{12,64},{12,70},{12,73},{12,82},{12,88},{12,91},{13,16},{13,20},
{13,24},{13,50},{13,56},{13,61},{13,62},{13,66},{13,74},{13,83},{13,86},{13,92},{14,17},{14,21},{14,22},
{14,51},{14,57},{14,62},{14,63},{14,64},{14,75},{14,84},{14,87},{14,93},{15,18},{15,19},{15,23},{15,49},
{15,55},{15,61},{15,63},{15,65},{15,73},{15,82},{15,85},{15,91},{16,21},{16,23},{16,53},{16,59},{16,63},
{16,64},{16,65},{16,77},{16,80},{16,89},{16,95},{17,19},{17,24},{17,54},{17,60},{17,61},{17,65},{17,66},
{17,78},{17,81},{17,90},{17,96},{18,20},{18,22},{18,52},{18,58},{18,62},{18,64},{18,66},{18,76},{18,79},
{18,88},{18,94},{19,22},{19,53},{19,56},{19,67},{19,68},{19,72},{19,74},{19,80},{19,89},{19,92},{20,23},
{20,54},{20,57},{20,68},{20,69},{20,70},{20,75},{20,81},{20,90},{20,93},{21,24},{21,52},{21,55},{21,67},
{21,69},{21,71},{21,73},{21,79},{21,88},{21,91},{22,50},{22,59},{22,69},{22,70},{22,71},{22,77},{22,83},
{22,86},{22,95},{23,51},{23,60},{23,67},{23,71},{23,72},{23,78},{23,84},{23,87},{23,96},{24,49},{24,58},
{24,68},{24,70},{24,72},{24,76},{24,82},{24,85},{24,94},{25,28},{25,32},{25,36},{25,50},{25,59},{25,62},
{25,68},{25,73},{25,74},{25,78},{25,86},{25,95},{26,29},{26,33},{26,34},{26,51},{26,60},{26,63},{26,69},
{26,74},{26,75},{26,76},{26,87},{26,96},{27,30},{27,31},{27,35},{27,49},{27,58},{27,61},{27,67},{27,73},
{27,75},{27,77},{27,85},{27,94},{28,33},{28,35},{28,53},{28,56},{28,65},{28,71},{28,75},{28,76},{28,77},
{28,89},{28,92},{29,31},{29,36},{29,54},{29,57},{29,66},{29,72},{29,73},{29,77},{29,78},{29,90},{29,93},
{30,32},{30,34},{30,52},{30,55},{30,64},{30,70},{30,74},{30,76},{30,78},{30,88},{30,91},{31,34},{31,50},
{31,56},{31,65},{31,68},{31,79},{31,80},{31,84},{31,86},{31,92},{32,35},{32,51},{32,57},{32,66},{32,69},
{32,80},{32,81},{32,82},{32,87},{32,93},{33,36},{33,49},{33,55},{33,64},{33,67},{33,79},{33,81},{33,83},
{33,85},{33,91},{34,53},{34,59},{34,62},{34,71},{34,81},{34,82},{34,83},{34,89},{34,95},{35,54},{35,60},
{35,63},{35,72},{35,79},{35,83},{35,84},{35,90},{35,96},{36,52},{36,58},{36,61},{36,70},{36,80},{36,82},
{36,84},{36,88},{36,94},{37,40},{37,44},{37,48},{37,53},{37,56},{37,62},{37,71},{37,74},{37,80},{37,85},
{37,86},{37,90},{38,41},{38,45},{38,46},{38,54},{38,57},{38,63},{38,72},{38,75},{38,81},{38,86},{38,87},
{38,88},{39,42},{39,43},{39,47},{39,52},{39,55},{39,61},{39,70},{39,73},{39,79},{39,85},{39,87},{39,89},
{40,45},{40,47},{40,50},{40,59},{40,65},{40,68},{40,77},{40,83},{40,87},{40,88},{40,89},{41,43},{41,48},
{41,51},{41,60},{41,66},{41,69},{41,78},{41,84},{41,85},{41,89},{41,90},{42,44},{42,46},{42,49},{42,58},
{42,64},{42,67},{42,76},{42,82},{42,86},{42,88},{42,90},{43,46},{43,53},{43,59},{43,62},{43,68},{43,77},
{43,80},{43,91},{43,92},{43,96},{44,47},{44,54},{44,60},{44,63},{44,69},{44,78},{44,81},{44,92},{44,93},
{44,94},{45,48},{45,52},{45,58},{45,61},{45,67},{45,76},{45,79},{45,91},{45,93},{45,95},{46,50},{46,56},
{46,65},{46,71},{46,74},{46,83},{46,93},{46,94},{46,95},{47,51},{47,57},{47,66},{47,72},{47,75},{47,84},
{47,91},{47,95},{47,96},{48,49},{48,55},{48,64},{48,70},{48,73},{48,82},{48,92},{48,94},{48,96},{49,52},
{49,57},{49,59},{50,53},{50,55},{50,60},{51,54},{51,56},{51,58},{52,56},{52,60},{53,57},{53,58},{54,55},
{54,59},{55,58},{56,59},{57,60},{61,64},{61,69},{61,71},{62,65},{62,67},{62,72},{63,66},{63,68},{63,70},
{64,68},{64,72},{65,69},{65,70},{66,67},{66,71},{67,70},{68,71},{69,72},{73,76},{73,81},{73,83},{74,77},
{74,79},{74,84},{75,78},{75,80},{75,82},{76,80},{76,84},{77,81},{77,82},{78,79},{78,83},{79,82},{80,83},
{81,84},{85,88},{85,93},{85,95},{86,89},{86,91},{86,96},{87,90},{87,92},{87,94},{88,92},{88,96},{89,93},
{89,94},{90,91},{90,95},{91,94},{92,95},{93,96}
\end{lstlisting}

\noindent SG96x41x1 

\begin{lstlisting}
{1,4},{1,7},{1,11},{1,12},{1,51},{1,58},{1,63},{1,67},{1,76},{1,84},{1,87},{1,94},{2,5},{2,8},{2,10},
{2,12},{2,49},{2,59},{2,61},{2,68},{2,77},{2,82},{2,85},{2,95},{3,6},{3,9},{3,10},{3,11},{3,50},{3,60},
{3,62},{3,69},{3,78},{3,83},{3,86},{3,96},{4,8},{4,9},{4,10},{4,54},{4,55},{4,66},{4,70},{4,73},{4,81},
{4,90},{4,91},{5,7},{5,9},{5,11},{5,52},{5,56},{5,64},{5,71},{5,74},{5,79},{5,88},{5,92},{6,7},{6,8},
{6,12},{6,53},{6,57},{6,65},{6,72},{6,75},{6,80},{6,89},{6,93},{7,10},{7,52},{7,57},{7,61},{7,69},
{7,78},{7,82},{7,88},{7,93},{8,11},{8,53},{8,55},{8,62},{8,67},{8,76},{8,83},{8,89},{8,91},{9,12},
{9,54},{9,56},{9,63},{9,68},{9,77},{9,84},{9,90},{9,92},{10,49},{10,60},{10,64},{10,72},{10,75},{10,79},
{10,85},{10,96},{11,50},{11,58},{11,65},{11,70},{11,73},{11,80},{11,86},{11,94},{12,51},{12,59},{12,66},
{12,71},{12,74},{12,81},{12,87},{12,95},{13,16},{13,19},{13,23},{13,24},{13,51},{13,55},{13,63},{13,70},
{13,75},{13,82},{13,88},{13,96},{14,17},{14,20},{14,22},{14,24},{14,49},{14,56},{14,61},{14,71},{14,73},
{14,83},{14,89},{14,94},{15,18},{15,21},{15,22},{15,23},{15,50},{15,57},{15,62},{15,72},{15,74},{15,84},
{15,90},{15,95},{16,20},{16,21},{16,22},{16,54},{16,58},{16,66},{16,67},{16,78},{16,79},{16,85},{16,93},
{17,19},{17,21},{17,23},{17,52},{17,59},{17,64},{17,68},{17,76},{17,80},{17,86},{17,91},{18,19},{18,20},
{18,24},{18,53},{18,60},{18,65},{18,69},{18,77},{18,81},{18,87},{18,92},{19,22},{19,49},{19,57},{19,64},
{19,69},{19,76},{19,81},{19,90},{19,94},{20,23},{20,50},{20,55},{20,65},{20,67},{20,77},{20,79},{20,88},
{20,95},{21,24},{21,51},{21,56},{21,66},{21,68},{21,78},{21,80},{21,89},{21,96},{22,52},{22,60},{22,61},
{22,72},{22,73},{22,84},{22,87},{22,91},{23,53},{23,58},{23,62},{23,70},{23,74},{23,82},{23,85},{23,92},
{24,54},{24,59},{24,63},{24,71},{24,75},{24,83},{24,86},{24,93},{25,28},{25,31},{25,35},{25,36},{25,52},
{25,60},{25,66},{25,67},{25,75},{25,82},{25,90},{25,94},{26,29},{26,32},{26,34},{26,36},{26,53},{26,58},
{26,64},{26,68},{26,73},{26,83},{26,88},{26,95},{27,30},{27,33},{27,34},{27,35},{27,54},{27,59},{27,65},
{27,69},{27,74},{27,84},{27,89},{27,96},{28,32},{28,33},{28,34},{28,49},{28,57},{28,63},{28,70},{28,78},
{28,79},{28,87},{28,91},{29,31},{29,33},{29,35},{29,50},{29,55},{29,61},{29,71},{29,76},{29,80},{29,85},
{29,92},{30,31},{30,32},{30,36},{30,51},{30,56},{30,62},{30,72},{30,77},{30,81},{30,86},{30,93},{31,34},
{31,54},{31,58},{31,61},{31,72},{31,76},{31,81},{31,88},{31,96},{32,35},{32,52},{32,59},{32,62},{32,70},
{32,77},{32,79},{32,89},{32,94},{33,36},{33,53},{33,60},{33,63},{33,71},{33,78},{33,80},{33,90},{33,95},
{34,51},{34,55},{34,64},{34,69},{34,73},{34,84},{34,85},{34,93},{35,49},{35,56},{35,65},{35,67},{35,74},
{35,82},{35,86},{35,91},{36,50},{36,57},{36,66},{36,68},{36,75},{36,83},{36,87},{36,92},{37,40},{37,43},
{37,47},{37,48},{37,54},{37,55},{37,64},{37,72},{37,78},{37,82},{37,87},{37,94},{38,41},{38,44},{38,46},
{38,48},{38,52},{38,56},{38,65},{38,70},{38,76},{38,83},{38,85},{38,95},{39,42},{39,45},{39,46},{39,47},
{39,53},{39,57},{39,66},{39,71},{39,77},{39,84},{39,86},{39,96},{40,44},{40,45},{40,46},{40,51},{40,58},
{40,61},{40,69},{40,75},{40,79},{40,90},{40,91},{41,43},{41,45},{41,47},{41,49},{41,59},{41,62},{41,67},
{41,73},{41,80},{41,88},{41,92},{42,43},{42,44},{42,48},{42,50},{42,60},{42,63},{42,68},{42,74},{42,81},
{42,89},{42,93},{43,46},{43,49},{43,60},{43,66},{43,70},{43,76},{43,84},{43,88},{43,93},{44,47},{44,50},
{44,58},{44,64},{44,71},{44,77},{44,82},{44,89},{44,91},{45,48},{45,51},{45,59},{45,65},{45,72},{45,78},
{45,83},{45,90},{45,92},{46,52},{46,57},{46,63},{46,67},{46,73},{46,81},{46,85},{46,96},{47,53},{47,55},
{47,61},{47,68},{47,74},{47,79},{47,86},{47,94},{48,54},{48,56},{48,62},{48,69},{48,75},{48,80},{48,87},
{48,95},{49,52},{49,56},{49,57},{49,58},{50,53},{50,55},{50,57},{50,59},{51,54},{51,55},{51,56},{51,60},
{52,55},{52,59},{52,60},{53,56},{53,58},{53,60},{54,57},{54,58},{54,59},{55,58},{56,59},{57,60},{61,64},
{61,68},{61,69},{61,70},{62,65},{62,67},{62,69},{62,71},{63,66},{63,67},{63,68},{63,72},{64,67},{64,71},
{64,72},{65,68},{65,70},{65,72},{66,69},{66,70},{66,71},{67,70},{68,71},{69,72},{73,76},{73,80},{73,81},
{73,82},{74,77},{74,79},{74,81},{74,83},{75,78},{75,79},{75,80},{75,84},{76,79},{76,83},{76,84},{77,80},
{77,82},{77,84},{78,81},{78,82},{78,83},{79,82},{80,83},{81,84},{85,88},{85,92},{85,93},{85,94},{86,89},
{86,91},{86,93},{86,95},{87,90},{87,91},{87,92},{87,96},{88,91},{88,95},{88,96},{89,92},{89,94},{89,96},
{90,93},{90,94},{90,95},{91,94},{92,95},{93,96}
\end{lstlisting}

\noindent SG96x51x3 

\begin{lstlisting}
{1,5},{1,6},{1,14},{1,22},{1,53},{1,55},{1,58},{1,67},{1,75},{1,78},{1,81},{1,95},{2,5},{2,6},{2,13},
{2,21},{2,54},{2,56},{2,57},{2,68},{2,76},{2,77},{2,82},{2,96},{3,7},{3,8},{3,16},{3,24},{3,54},{3,55},
{3,60},{3,66},{3,74},{3,80},{3,83},{3,94},{4,7},{4,8},{4,15},{4,23},{4,53},{4,56},{4,59},{4,65},{4,73},
{4,79},{4,84},{4,93},{5,10},{5,18},{5,49},{5,51},{5,62},{5,71},{5,74},{5,79},{5,85},{5,91},{6,9},
{6,17},{6,50},{6,52},{6,61},{6,72},{6,73},{6,80},{6,86},{6,92},{7,12},{7,20},{7,50},{7,51},{7,64},
{7,70},{7,76},{7,78},{7,87},{7,90},{8,11},{8,19},{8,49},{8,52},{8,63},{8,69},{8,75},{8,77},{8,88},
{8,89},{9,13},{9,14},{9,22},{9,51},{9,61},{9,63},{9,66},{9,79},{9,83},{9,86},{9,89},{10,13},{10,14},
{10,21},{10,52},{10,62},{10,64},{10,65},{10,80},{10,84},{10,85},{10,90},{11,15},{11,16},{11,24},{11,50},
{11,62},{11,63},{11,68},{11,78},{11,82},{11,88},{11,91},{12,15},{12,16},{12,23},{12,49},{12,61},{12,64},
{12,67},{12,77},{12,81},{12,87},{12,92},{13,18},{13,55},{13,57},{13,59},{13,70},{13,75},{13,82},{13,87},
{13,93},{14,17},{14,56},{14,58},{14,60},{14,69},{14,76},{14,81},{14,88},{14,94},{15,20},{15,54},{15,58},
{15,59},{15,72},{15,74},{15,84},{15,86},{15,95},{16,19},{16,53},{16,57},{16,60},{16,71},{16,73},{16,83},
{16,85},{16,96},{17,21},{17,22},{17,50},{17,59},{17,69},{17,71},{17,73},{17,87},{17,91},{17,94},{18,21},
{18,22},{18,49},{18,60},{18,70},{18,72},{18,74},{18,88},{18,92},{18,93},{19,23},{19,24},{19,52},{19,58},
{19,70},{19,71},{19,75},{19,86},{19,90},{19,96},{20,23},{20,24},{20,51},{20,57},{20,69},{20,72},{20,76},
{20,85},{20,89},{20,95},{21,54},{21,63},{21,65},{21,67},{21,77},{21,83},{21,90},{21,95},{22,53},{22,64},
{22,66},{22,68},{22,78},{22,84},{22,89},{22,96},{23,56},{23,62},{23,66},{23,67},{23,79},{23,82},{23,92},
{23,94},{24,55},{24,61},{24,65},{24,68},{24,80},{24,81},{24,91},{24,93},{25,29},{25,30},{25,38},{25,46},
{25,51},{25,53},{25,58},{25,71},{25,77},{25,80},{25,82},{25,92},{26,29},{26,30},{26,37},{26,45},{26,52},
{26,54},{26,57},{26,72},{26,78},{26,79},{26,81},{26,91},{27,31},{27,32},{27,40},{27,48},{27,50},{27,55},
{27,60},{27,70},{27,77},{27,79},{27,84},{27,89},{28,31},{28,32},{28,39},{28,47},{28,49},{28,56},{28,59},
{28,69},{28,78},{28,80},{28,83},{28,90},{29,34},{29,42},{29,49},{29,55},{29,62},{29,67},{29,73},{29,76},
{29,86},{29,96},{30,33},{30,41},{30,50},{30,56},{30,61},{30,68},{30,74},{30,75},{30,85},{30,95},{31,36},
{31,44},{31,51},{31,54},{31,64},{31,66},{31,73},{31,75},{31,88},{31,93},{32,35},{32,43},{32,52},{32,53},
{32,63},{32,65},{32,74},{32,76},{32,87},{32,94},{33,37},{33,38},{33,46},{33,55},{33,59},{33,61},{33,66},
{33,76},{33,85},{33,88},{33,90},{34,37},{34,38},{34,45},{34,56},{34,60},{34,62},{34,65},{34,75},{34,86},
{34,87},{34,89},{35,39},{35,40},{35,48},{35,54},{35,58},{35,63},{35,68},{35,73},{35,85},{35,87},{35,92},
{36,39},{36,40},{36,47},{36,53},{36,57},{36,64},{36,67},{36,74},{36,86},{36,88},{36,91},{37,42},{37,51},
{37,57},{37,63},{37,70},{37,80},{37,81},{37,84},{37,94},{38,41},{38,52},{38,58},{38,64},{38,69},{38,79},
{38,82},{38,83},{38,93},{39,44},{39,50},{39,59},{39,62},{39,72},{39,77},{39,81},{39,83},{39,96},{40,43},
{40,49},{40,60},{40,61},{40,71},{40,78},{40,82},{40,84},{40,95},{41,45},{41,46},{41,50},{41,63},{41,67},
{41,69},{41,74},{41,84},{41,93},{41,96},{42,45},{42,46},{42,49},{42,64},{42,68},{42,70},{42,73},{42,83},
{42,94},{42,95},{43,47},{43,48},{43,52},{43,62},{43,66},{43,71},{43,76},{43,81},{43,93},{43,95},{44,47},
{44,48},{44,51},{44,61},{44,65},{44,72},{44,75},{44,82},{44,94},{44,96},{45,54},{45,59},{45,65},{45,71},
{45,78},{45,88},{45,89},{45,92},{46,53},{46,60},{46,66},{46,72},{46,77},{46,87},{46,90},{46,91},{47,56},
{47,58},{47,67},{47,70},{47,80},{47,85},{47,89},{47,91},{48,55},{48,57},{48,68},{48,69},{48,79},{48,86},
{48,90},{48,92},{49,53},{49,54},{49,61},{49,69},{50,53},{50,54},{50,62},{50,70},{51,55},{51,56},{51,63},
{51,71},{52,55},{52,56},{52,64},{52,72},{53,57},{53,65},{54,58},{54,66},{55,59},{55,67},{56,60},{56,68},
{57,61},{57,62},{57,69},{58,61},{58,62},{58,70},{59,63},{59,64},{59,71},{60,63},{60,64},{60,72},{61,65},
{62,66},{63,67},{64,68},{65,69},{65,70},{66,69},{66,70},{67,71},{67,72},{68,71},{68,72},{73,77},{73,78},
{73,85},{73,93},{74,77},{74,78},{74,86},{74,94},{75,79},{75,80},{75,87},{75,95},{76,79},{76,80},{76,88},
{76,96},{77,81},{77,89},{78,82},{78,90},{79,83},{79,91},{80,84},{80,92},{81,85},{81,86},{81,93},{82,85},
{82,86},{82,94},{83,87},{83,88},{83,95},{84,87},{84,88},{84,96},{85,89},{86,90},{87,91},{88,92},{89,93},
{89,94},{90,93},{90,94},{91,95},{91,96},{92,95},{92,96}
\end{lstlisting}

\noindent SG96x182x1 

\begin{lstlisting}
{1,2},{1,15},{1,25},{1,28},{1,32},{1,36},{1,44},{1,47},{1,50},{1,51},{1,86},{1,87},{2,16},{2,26},{2,27},
{2,31},{2,35},{2,43},{2,48},{2,49},{2,52},{2,85},{2,88},{3,4},{3,14},{3,25},{3,27},{3,29},{3,33},{3,41},
{3,46},{3,50},{3,52},{3,86},{3,88},{4,13},{4,26},{4,28},{4,30},{4,34},{4,42},{4,45},{4,49},{4,51},
{4,85},{4,87},{5,6},{5,19},{5,28},{5,29},{5,32},{5,36},{5,39},{5,48},{5,54},{5,55},{5,90},{5,91},{6,20},
{6,27},{6,30},{6,31},{6,35},{6,40},{6,47},{6,53},{6,56},{6,89},{6,92},{7,8},{7,18},{7,25},{7,29},{7,31},
{7,33},{7,38},{7,45},{7,54},{7,56},{7,90},{7,92},{8,17},{8,26},{8,30},{8,32},{8,34},{8,37},{8,46},
{8,53},{8,55},{8,89},{8,91},{9,10},{9,23},{9,28},{9,32},{9,33},{9,36},{9,40},{9,43},{9,58},{9,59},
{9,94},{9,95},{10,24},{10,27},{10,31},{10,34},{10,35},{10,39},{10,44},{10,57},{10,60},{10,93},{10,96},
{11,12},{11,22},{11,25},{11,29},{11,33},{11,35},{11,37},{11,42},{11,58},{11,60},{11,94},{11,96},{12,21},
{12,26},{12,30},{12,34},{12,36},{12,38},{12,41},{12,57},{12,59},{12,93},{12,95},{13,14},{13,31},{13,36},
{13,37},{13,40},{13,44},{13,48},{13,62},{13,63},{13,73},{13,76},{14,32},{14,35},{14,38},{14,39},{14,43},
{14,47},{14,61},{14,64},{14,74},{14,75},{15,16},{15,30},{15,33},{15,37},{15,39},{15,41},{15,45},{15,62},
{15,64},{15,73},{15,75},{16,29},{16,34},{16,38},{16,40},{16,42},{16,46},{16,61},{16,63},{16,74},{16,76},
{17,18},{17,28},{17,35},{17,40},{17,41},{17,44},{17,48},{17,66},{17,67},{17,77},{17,80},{18,27},{18,36},
{18,39},{18,42},{18,43},{18,47},{18,65},{18,68},{18,78},{18,79},{19,20},{19,25},{19,34},{19,37},{19,41},
{19,43},{19,45},{19,66},{19,68},{19,77},{19,79},{20,26},{20,33},{20,38},{20,42},{20,44},{20,46},{20,65},
{20,67},{20,78},{20,80},{21,22},{21,27},{21,32},{21,40},{21,44},{21,45},{21,48},{21,70},{21,71},{21,81},
{21,84},{22,28},{22,31},{22,39},{22,43},{22,46},{22,47},{22,69},{22,72},{22,82},{22,83},{23,24},{23,26},
{23,29},{23,37},{23,41},{23,45},{23,47},{23,70},{23,72},{23,81},{23,83},{24,25},{24,30},{24,38},{24,42},
{24,46},{24,48},{24,69},{24,71},{24,82},{24,84},{25,26},{25,40},{25,62},{25,64},{25,74},{25,76},{26,39},
{26,61},{26,63},{26,73},{26,75},{27,28},{27,37},{27,61},{27,64},{27,73},{27,76},{28,38},{28,62},{28,63},
{28,74},{28,75},{29,30},{29,44},{29,66},{29,68},{29,78},{29,80},{30,43},{30,65},{30,67},{30,77},{30,79},
{31,32},{31,41},{31,65},{31,68},{31,77},{31,80},{32,42},{32,66},{32,67},{32,78},{32,79},{33,34},{33,48},
{33,70},{33,72},{33,82},{33,84},{34,47},{34,69},{34,71},{34,81},{34,83},{35,36},{35,45},{35,69},{35,72},
{35,81},{35,84},{36,46},{36,70},{36,71},{36,82},{36,83},{37,38},{37,49},{37,51},{37,86},{37,88},{38,50},
{38,52},{38,85},{38,87},{39,40},{39,50},{39,51},{39,85},{39,88},{40,49},{40,52},{40,86},{40,87},{41,42},
{41,53},{41,55},{41,90},{41,92},{42,54},{42,56},{42,89},{42,91},{43,44},{43,54},{43,55},{43,89},{43,92},
{44,53},{44,56},{44,90},{44,91},{45,46},{45,57},{45,59},{45,94},{45,96},{46,58},{46,60},{46,93},{46,95},
{47,48},{47,58},{47,59},{47,93},{47,96},{48,57},{48,60},{48,94},{48,95},{49,50},{49,64},{49,74},{49,75},
{49,78},{49,82},{49,90},{49,93},{50,63},{50,73},{50,76},{50,77},{50,81},{50,89},{50,94},{51,52},{51,61},
{51,74},{51,76},{51,80},{51,84},{51,92},{51,95},{52,62},{52,73},{52,75},{52,79},{52,83},{52,91},{52,96},
{53,54},{53,68},{53,74},{53,78},{53,79},{53,82},{53,85},{53,94},{54,67},{54,73},{54,77},{54,80},{54,81},
{54,86},{54,93},{55,56},{55,65},{55,76},{55,78},{55,80},{55,84},{55,87},{55,96},{56,66},{56,75},{56,77},
{56,79},{56,83},{56,88},{56,95},{57,58},{57,72},{57,74},{57,78},{57,82},{57,83},{57,86},{57,89},{58,71},
{58,73},{58,77},{58,81},{58,84},{58,85},{58,90},{59,60},{59,69},{59,76},{59,80},{59,82},{59,84},{59,88},
{59,91},{60,70},{60,75},{60,79},{60,81},{60,83},{60,87},{60,92},{61,62},{61,77},{61,82},{61,86},{61,87},
{61,90},{61,94},{62,78},{62,81},{62,85},{62,88},{62,89},{62,93},{63,64},{63,79},{63,84},{63,86},{63,88},
{63,92},{63,96},{64,80},{64,83},{64,85},{64,87},{64,91},{64,95},{65,66},{65,74},{65,81},{65,86},{65,90},
{65,91},{65,94},{66,73},{66,82},{66,85},{66,89},{66,92},{66,93},{67,68},{67,76},{67,83},{67,88},{67,90},
{67,92},{67,96},{68,75},{68,84},{68,87},{68,89},{68,91},{68,95},{69,70},{69,73},{69,78},{69,86},{69,90},
{69,94},{69,95},{70,74},{70,77},{70,85},{70,89},{70,93},{70,96},{71,72},{71,75},{71,80},{71,88},{71,92},
{71,94},{71,96},{72,76},{72,79},{72,87},{72,91},{72,93},{72,95},{73,74},{73,87},{74,88},{75,76},{75,86},
{76,85},{77,78},{77,91},{78,92},{79,80},{79,90},{80,89},{81,82},{81,95},{82,96},{83,84},{83,94},{84,93},
{85,86},{87,88},{89,90},{91,92},{93,94},{95,96}
\end{lstlisting}

\noindent SG96x185x1  

\begin{lstlisting}
{1,4},{1,5},{1,6},{1,9},{1,11},{1,14},{1,26},{1,38},{1,51},{1,62},{1,74},{1,87},{2,3},{2,5},{2,6},{2,10}
{2,12},{2,13},{2,25},{2,37},{2,52},{2,61},{2,73},{2,88},{3,7},{3,8},{3,9},{3,11},{3,16},{3,28},{3,40},
{3,49},{3,64},{3,76},{3,85},{4,7},{4,8},{4,10},{4,12},{4,15},{4,27},{4,39},{4,50},{4,63},{4,75},{4,86},
{5,7},{5,9},{5,12},{5,20},{5,32},{5,44},{5,54},{5,68},{5,80},{5,90},{6,8},{6,10},{6,11},{6,19},{6,31},
{6,43},{6,53},{6,67},{6,79},{6,89},{7,10},{7,11},{7,18},{7,30},{7,42},{7,56},{7,66},{7,78},{7,92},{8,9},
{8,12},{8,17},{8,29},{8,41},{8,55},{8,65},{8,77},{8,91},{9,10},{9,23},{9,35},{9,47},{9,60},{9,71},
{9,83},{9,96},{10,24},{10,36},{10,48},{10,59},{10,72},{10,84},{10,95},{11,12},{11,21},{11,33},{11,45},
{11,58},{11,69},{11,81},{11,94},{12,22},{12,34},{12,46},{12,57},{12,70},{12,82},{12,93},{13,16},{13,17},
{13,18},{13,21},{13,23},{13,26},{13,38},{13,50},{13,63},{13,75},{13,86},{14,15},{14,17},{14,18},{14,22},
{14,24},{14,25},{14,37},{14,49},{14,64},{14,76},{14,85},{15,19},{15,20},{15,21},{15,23},{15,28},{15,40},
{15,52},{15,61},{15,73},{15,88},{16,19},{16,20},{16,22},{16,24},{16,27},{16,39},{16,51},{16,62},{16,74},
{16,87},{17,19},{17,21},{17,24},{17,32},{17,44},{17,56},{17,66},{17,78},{17,92},{18,20},{18,22},{18,23},
{18,31},{18,43},{18,55},{18,65},{18,77},{18,91},{19,22},{19,23},{19,30},{19,42},{19,54},{19,68},{19,80},
{19,90},{20,21},{20,24},{20,29},{20,41},{20,53},{20,67},{20,79},{20,89},{21,22},{21,35},{21,47},{21,59},
{21,72},{21,84},{21,95},{22,36},{22,48},{22,60},{22,71},{22,83},{22,96},{23,24},{23,33},{23,45},{23,57},
{23,70},{23,82},{23,93},{24,34},{24,46},{24,58},{24,69},{24,81},{24,94},{25,28},{25,29},{25,30},{25,33},
{25,35},{25,38},{25,51},{25,62},{25,75},{25,86},{26,27},{26,29},{26,30},{26,34},{26,36},{26,37},{26,52},
{26,61},{26,76},{26,85},{27,31},{27,32},{27,33},{27,35},{27,40},{27,49},{27,64},{27,73},{27,88},{28,31},
{28,32},{28,34},{28,36},{28,39},{28,50},{28,63},{28,74},{28,87},{29,31},{29,33},{29,36},{29,44},{29,54},
{29,68},{29,78},{29,92},{30,32},{30,34},{30,35},{30,43},{30,53},{30,67},{30,77},{30,91},{31,34},{31,35},
{31,42},{31,56},{31,66},{31,80},{31,90},{32,33},{32,36},{32,41},{32,55},{32,65},{32,79},{32,89},{33,34},
{33,47},{33,60},{33,71},{33,84},{33,95},{34,48},{34,59},{34,72},{34,83},{34,96},{35,36},{35,45},{35,58},
{35,69},{35,82},{35,93},{36,46},{36,57},{36,70},{36,81},{36,94},{37,40},{37,41},{37,42},{37,45},{37,47},
{37,50},{37,63},{37,74},{37,87},{38,39},{38,41},{38,42},{38,46},{38,48},{38,49},{38,64},{38,73},{38,88},
{39,43},{39,44},{39,45},{39,47},{39,52},{39,61},{39,76},{39,85},{40,43},{40,44},{40,46},{40,48},{40,51},
{40,62},{40,75},{40,86},{41,43},{41,45},{41,48},{41,56},{41,66},{41,80},{41,90},{42,44},{42,46},{42,47},
{42,55},{42,65},{42,79},{42,89},{43,46},{43,47},{43,54},{43,68},{43,78},{43,92},{44,45},{44,48},{44,53},
{44,67},{44,77},{44,91},{45,46},{45,59},{45,72},{45,83},{45,96},{46,60},{46,71},{46,84},{46,95},{47,48},
{47,57},{47,70},{47,81},{47,94},{48,58},{48,69},{48,82},{48,93},{49,52},{49,53},{49,54},{49,57},{49,59},
{49,63},{49,75},{49,87},{50,51},{50,53},{50,54},{50,58},{50,60},{50,64},{50,76},{50,88},{51,55},{51,56},
{51,57},{51,59},{51,61},{51,73},{51,85},{52,55},{52,56},{52,58},{52,60},{52,62},{52,74},{52,86},{53,55},
{53,57},{53,60},{53,66},{53,78},{53,90},{54,56},{54,58},{54,59},{54,65},{54,77},{54,89},{55,58},{55,59},
{55,68},{55,80},{55,92},{56,57},{56,60},{56,67},{56,79},{56,91},{57,58},{57,72},{57,84},{57,96},{58,71},
{58,83},{58,95},{59,60},{59,70},{59,82},{59,94},{60,69},{60,81},{60,93},{61,64},{61,65},{61,66},{61,69},
{61,71},{61,75},{61,87},{62,63},{62,65},{62,66},{62,70},{62,72},{62,76},{62,88},{63,67},{63,68},{63,69},
{63,71},{63,73},{63,85},{64,67},{64,68},{64,70},{64,72},{64,74},{64,86},{65,67},{65,69},{65,72},{65,78},
{65,90},{66,68},{66,70},{66,71},{66,77},{66,89},{67,70},{67,71},{67,80},{67,92},{68,69},{68,72},{68,79},
{68,91},{69,70},{69,84},{69,96},{70,83},{70,95},{71,72},{71,82},{71,94},{72,81},{72,93},{73,76},{73,77},
{73,78},{73,81},{73,83},{73,87},{74,75},{74,77},{74,78},{74,82},{74,84},{74,88},{75,79},{75,80},{75,81},
{75,83},{75,85},{76,79},{76,80},{76,82},{76,84},{76,86},{77,79},{77,81},{77,84},{77,90},{78,80},{78,82},
{78,83},{78,89},{79,82},{79,83},{79,92},{80,81},{80,84},{80,91},{81,82},{81,96},{82,95},{83,84},{83,94},
{84,93},{85,88},{85,89},{85,90},{85,93},{85,95},{86,87},{86,89},{86,90},{86,94},{86,96},{87,91},{87,92},
{87,93},{87,95},{88,91},{88,92},{88,94},{88,96},{89,91},{89,93},{89,96},{90,92},{90,94},{90,95},{91,94},
{91,95},{92,93},{92,96},{93,94},{95,96}
\end{lstlisting}

\noindent SG96x191x1 

\begin{lstlisting}
{1,2},{1,5},{1,6},{1,32},{1,34},{1,40},{1,41},{1,44},{1,68},{1,71},{1,91},{1,96},{2,5},{2,6},{2,31},
{2,33},{2,39},{2,42},{2,43},{2,67},{2,72},{2,92},{2,95},{3,4},{3,7},{3,8},{3,29},{3,36},{3,37},{3,41},
{3,43},{3,65},{3,70},{3,90},{3,93},{4,7},{4,8},{4,30},{4,35},{4,38},{4,42},{4,44},{4,66},{4,69},{4,89},
{4,94},{5,6},{5,28},{5,36},{5,38},{5,45},{5,47},{5,67},{5,71},{5,92},{5,96},{6,27},{6,35},{6,37},{6,46},
{6,48},{6,68},{6,72},{6,91},{6,95},{7,8},{7,25},{7,33},{7,40},{7,46},{7,47},{7,66},{7,70},{7,89},{7,93},
{8,26},{8,34},{8,39},{8,45},{8,48},{8,65},{8,69},{8,90},{8,94},{9,10},{9,11},{9,12},{9,25},{9,32},
{9,37},{9,42},{9,45},{9,54},{9,56},{9,77},{9,79},{10,11},{10,12},{10,26},{10,31},{10,38},{10,41},
{10,46},{10,53},{10,55},{10,78},{10,80},{11,12},{11,27},{11,29},{11,39},{11,44},{11,47},{11,53},{11,56},
{11,78},{11,79},{12,28},{12,30},{12,40},{12,43},{12,48},{12,54},{12,55},{12,77},{12,80},{13,14},{13,15},
{13,16},{13,28},{13,29},{13,34},{13,42},{13,46},{13,49},{13,52},{13,74},{13,75},{14,15},{14,16},{14,27},
{14,30},{14,33},{14,41},{14,45},{14,50},{14,51},{14,73},{14,76},{15,16},{15,25},{15,31},{15,36},{15,44},
{15,48},{15,49},{15,51},{15,74},{15,76},{16,26},{16,32},{16,35},{16,43},{16,47},{16,50},{16,52},{16,73},
{16,75},{17,18},{17,23},{17,24},{17,26},{17,28},{17,33},{17,37},{17,44},{17,59},{17,61},{17,84},{17,86},
{18,23},{18,24},{18,25},{18,27},{18,34},{18,38},{18,43},{18,60},{18,62},{18,83},{18,85},{19,20},{19,21},
{19,22},{19,25},{19,28},{19,35},{19,39},{19,41},{19,58},{19,63},{19,81},{19,88},{20,21},{20,22},{20,26},
{20,27},{20,36},{20,40},{20,42},{20,57},{20,64},{20,82},{20,87},{21,22},{21,30},{21,31},{21,34},{21,37},
{21,47},{21,58},{21,64},{21,81},{21,87},{22,29},{22,32},{22,33},{22,38},{22,48},{22,57},{22,63},{22,82},
{22,88},{23,24},{23,30},{23,32},{23,36},{23,39},{23,46},{23,60},{23,61},{23,83},{23,86},{24,29},{24,31},
{24,35},{24,40},{24,45},{24,59},{24,62},{24,84},{24,85},{25,26},{25,29},{25,30},{25,68},{25,71},{25,92},
{25,95},{26,29},{26,30},{26,67},{26,72},{26,91},{26,96},{27,28},{27,31},{27,32},{27,65},{27,70},{27,89},
{27,94},{28,31},{28,32},{28,66},{28,69},{28,90},{28,93},{29,30},{29,67},{29,71},{29,91},{29,95},{30,68},
{30,72},{30,92},{30,96},{31,32},{31,66},{31,70},{31,90},{31,94},{32,65},{32,69},{32,89},{32,93},{33,34},
{33,35},{33,36},{33,54},{33,56},{33,78},{33,80},{34,35},{34,36},{34,53},{34,55},{34,77},{34,79},{35,36},
{35,53},{35,56},{35,77},{35,80},{36,54},{36,55},{36,78},{36,79},{37,38},{37,39},{37,40},{37,49},{37,52},
{37,73},{37,76},{38,39},{38,40},{38,50},{38,51},{38,74},{38,75},{39,40},{39,49},{39,51},{39,73},{39,75},
{40,50},{40,52},{40,74},{40,76},{41,42},{41,47},{41,48},{41,59},{41,61},{41,83},{41,85},{42,47},{42,48},
{42,60},{42,62},{42,84},{42,86},{43,44},{43,45},{43,46},{43,58},{43,63},{43,82},{43,87},{44,45},{44,46},
{44,57},{44,64},{44,81},{44,88},{45,46},{45,58},{45,64},{45,82},{45,88},{46,57},{46,63},{46,81},{46,87},
{47,48},{47,60},{47,61},{47,84},{47,85},{48,59},{48,62},{48,83},{48,86},{49,50},{49,51},{49,52},{49,80},
{49,82},{49,85},{49,89},{49,96},{50,51},{50,52},{50,79},{50,81},{50,86},{50,90},{50,95},{51,52},{51,77},
{51,84},{51,87},{51,91},{51,93},{52,78},{52,83},{52,88},{52,92},{52,94},{53,54},{53,55},{53,56},{53,76},
{53,82},{53,86},{53,92},{53,93},{54,55},{54,56},{54,75},{54,81},{54,85},{54,91},{54,94},{55,56},{55,73},
{55,84},{55,88},{55,89},{55,95},{56,74},{56,83},{56,87},{56,90},{56,96},{57,58},{57,63},{57,64},{57,73},
{57,77},{57,85},{57,90},{57,92},{58,63},{58,64},{58,74},{58,78},{58,86},{58,89},{58,91},{59,60},{59,61},
{59,62},{59,75},{59,79},{59,87},{59,89},{59,92},{60,61},{60,62},{60,76},{60,80},{60,88},{60,90},{60,91},
{61,62},{61,74},{61,77},{61,82},{61,94},{61,95},{62,73},{62,78},{62,81},{62,93},{62,96},{63,64},{63,76},
{63,79},{63,84},{63,94},{63,96},{64,75},{64,80},{64,83},{64,93},{64,95},{65,66},{65,69},{65,70},{65,74},
{65,80},{65,81},{65,84},{65,92},{66,69},{66,70},{66,73},{66,79},{66,82},{66,83},{66,91},{67,68},{67,71},
{67,72},{67,76},{67,77},{67,81},{67,83},{67,89},{68,71},{68,72},{68,75},{68,78},{68,82},{68,84},{68,90},
{69,70},{69,76},{69,78},{69,85},{69,87},{69,95},{70,75},{70,77},{70,86},{70,88},{70,96},{71,72},{71,73},
{71,80},{71,86},{71,87},{71,94},{72,74},{72,79},{72,85},{72,88},{72,93},{73,74},{73,75},{73,76},{74,75},
{74,76},{75,76},{77,78},{77,79},{77,80},{78,79},{78,80},{79,80},{81,82},{81,87},{81,88},{82,87},{82,88},
{83,84},{83,85},{83,86},{84,85},{84,86},{85,86},{87,88},{89,90},{89,93},{89,94},{90,93},{90,94},{91,92},
{91,95},{91,96},{92,95},{92,96},{93,94},{95,96}
\end{lstlisting}

\subsection*{(12,3)= 620}

\noindent $Z_{20}\rtimes_{2}Z_{31}$    generators $[18,17]\-[2,25 ];[14,13]\-[6,5 ];[16,16]\-[4,23 ];[19,15]\-[1,1 ];[3,17]\-[17,25 ];[5,15]\-[15,16 ]$\\
\noindent $Z_{20}\rtimes_{4}Z_{31}$    generators $[8,16]\-[12,23 ];[9,20]\-[11,13 ];[12,25]\-[8,12 ];[11,24]\-[9,25 ];[2,0]\-[18,0 ];[14,20]\-[6,13 ]$\\
\noindent $Z_{20}\rtimes_{8}Z_{31}$    generators $[18,3]\-[2,25 ];[16,23]\-[4,1 ];[5,1]\-[15,30 ];[19,30]\-[1,8 ];[17,27]\-[3,2 ];[7,16]\-[13,23 ]$\\
\noindent $Z_{20}\rtimes_{15}Z_{31}$  generators $[3,7]\-[17,25 ];[7,8]\-[13,1 ];[3,24]\-[17,6 ];[19,28]\-[1,14 ];[16,10]\-[4,11 ];[6,5]\-[14,21 ]$\\
\noindent $Z_{20}\rtimes_{16}Z_{31}$  generators $[19,20]\-[1,21 ];[9,24]\-[11,19 ];[8,10]\-[12,13 ];[18,3]\-[2,7 ];[14,21]\-[6,5 ];[17,22]\-[3,5 ]$\\
\noindent $Z_{20}\rtimes_{23}Z_{31}$  generators $[12,26]\-[8,18 ];[12,19]\-[8,6 ];[2,23]\-[18,4 ];[11,9]\-[9,5 ];[3,10]\-[17,20 ];[1,18]\-[19,10 ]$\\
\noindent $Z_{20}\rtimes_{27}Z_{31}$  generators $[12,25]\-[8,12 ];[12,21]\-[8,20 ];[1,18]\-[19,20 ];[18,1]\-[2,15 ];[13,19]\-[7,25 ];[16,15]\-[4,4 ]$\\
\noindent $Z_{20}\rtimes_{29}Z_{31}$  generators $[2,22]\-[18,10 ];[11,20]\-[9,10 ];[3,9]\-[17,5 ];[1,16]\-[19,8 ];[16,10]\-[4,26 ];[3,11]\-[17,13 ]$\\

\subsection*{(13,2)=112}

\noindent $Z_{2}\rtimes_{43}Z_{56}$ generators $[1,4]\-[1,52 ];[1,54]\-[1,30 ];[1,21]\-[1,49 ];[0,18]\-[0,38 ];[0,35]\-[0,21 ];[1,8]\-[1,48 ];[1,42]$\\
\noindent $Z_{2}\rtimes_{43}Z_{56}$ generators $[0,54]\-[0,2 ];[1,34]\-[1,50 ];[0,35]\-[0,21 ];[1,32]\-[1,24 ];[1,49]\-[1,21 ];[1,44]\-[1,12 ];[1,14]$\\
\noindent $Z_{2}\rtimes_{43}Z_{56}$ generators $[1,38]\-[1,46 ];[0,49]\-[0,7 ];[1,49]\-[1,21 ];[1,2]\-[1,26 ];[0,54]\-[0,2 ];[1,48]\-[1,8 ];[1,0]$\\
\noindent SG(120,2) 
\begin{lstlisting}
{1,24},{1,37},{1,41},{1,53},{1,57},{1,58},{1,59},{1,61},{1,64},{1,66},{1,82},{1,99},{1,108},{2,21},
{2,25},{2,26},{2,27},{2,29},{2,32},{2,40},{2,52},{2,68},{2,72},{2,81},{2,90},{2,108},{3,14},{3,25},
{3,35},{3,56},{3,67},{3,75},{3,83},{3,87},{3,88},{3,89},{3,91},{3,94},{3,99},{4,7},{4,44},{4,58},
{4,59},{4,60},{4,62},{4,63},{4,67},{4,71},{4,81},{4,86},{4,96},{4,105},{5,22},{5,23},{5,28},{5,40},
{5,49},{5,56},{5,62},{5,64},{5,65},{5,68},{5,72},{5,83},{5,87},{6,14},{6,21},{6,22},{6,23},{6,25},
{6,33},{6,61},{6,63},{6,75},{6,80},{6,90},{6,95},{6,104},{7,10},{7,13},{7,17},{7,20},{7,22},{7,59},
{7,63},{7,64},{7,66},{7,81},{7,99},{7,105},{8,12},{8,13},{8,14},{8,32},{8,41},{8,46},{8,62},{8,78},
{8,93},{8,97},{8,100},{8,107},{8,109},{9,12},{9,13},{9,31},{9,40},{9,44},{9,45},{9,61},{9,79},
{9,94},{9,98},{9,101},{9,106},{9,110},{10,20},{10,41},{10,42},{10,43},{10,45},{10,48},{10,51},
{10,60},{10,71},{10,83},{10,90},{10,96},{11,14},{11,26},{11,27},{11,28},{11,30},{11,31},{11,51},
{11,56},{11,75},{11,82},{11,87},{11,95},{11,105},{12,17},{12,24},{12,30},{12,32},{12,44},{12,53},
{12,57},{12,71},{12,80},{12,83},{12,96},{13,29},{13,31},{13,37},{13,41},{13,49},{13,50},{13,60},
{13,86},{13,99},{13,104},{14,27},{14,31},{14,32},{14,33},{14,36},{14,51},{14,81},{14,87},{14,108},
{15,30},{15,35},{15,36},{15,37},{15,48},{15,63},{15,72},{15,77},{15,89},{15,91},{15,102},{15,106},
{15,109},{16,29},{16,35},{16,36},{16,47},{16,64},{16,71},{16,75},{16,76},{16,88},{16,92},{16,103},
{16,107},{16,110},{17,36},{17,44},{17,57},{17,68},{17,75},{17,86},{17,105},{17,106},{17,107},
{17,109},{17,112},{18,21},{18,28},{18,36},{18,40},{18,66},{18,72},{18,88},{18,89},{18,90},{18,92},
{18,93},{18,96},{18,104},{19,24},{19,33},{19,34},{19,37},{19,41},{19,54},{19,55},{19,60},{19,68},
{19,86},{19,92},{19,94},{19,105},{20,44},{20,53},{20,54},{20,55},{20,57},{20,65},{20,71},{20,91},
{20,93},{20,95},{20,108},{21,35},{21,49},{21,52},{21,54},{21,66},{21,75},{21,89},{21,93},{21,94},
{21,96},{22,25},{22,30},{22,47},{22,64},{22,73},{22,76},{22,92},{22,98},{22,100},{22,111},{23,24},
{23,28},{23,29},{23,48},{23,63},{23,74},{23,77},{23,91},{23,97},{23,101},{23,112},{24,34},{24,50},
{24,58},{24,63},{24,64},{24,67},{24,82},{24,96},{24,108},{25,35},{25,51},{25,58},{25,61},{25,62},
{25,65},{25,68},{25,80},{25,104},{26,30},{26,40},{26,41},{26,47},{26,63},{26,70},{26,76},{26,84},
{26,91},{26,103},{26,112},{27,29},{27,44},{27,48},{27,64},{27,69},{27,77},{27,85},{27,92},
{27,102},{27,111},{28,33},{28,36},{28,49},{28,52},{28,59},{28,61},{28,62},{28,83},{28,95},{29,41},
{29,44},{29,69},{29,73},{29,84},{29,88},{29,102},{29,112},{30,41},{30,44},{30,70},{30,74},{30,85},
{30,89},{30,103},{30,111},{31,38},{31,40},{31,42},{31,54},{31,58},{31,98},{31,101},{31,107},
{31,109},{32,39},{32,40},{32,43},{32,55},{32,59},{32,97},{32,100},{32,106},{32,110},{33,53},{33,60},
{33,72},{33,73},{33,74},{33,76},{33,79},{33,95},{33,105},{34,37},{34,42},{34,43},{34,44},{34,46},
{34,47},{34,52},{34,57},{34,87},{34,99},{34,104},{35,40},{35,46},{35,48},{35,50},{35,53},{35,75},
{35,95},{35,105},{36,45},{36,47},{36,56},{36,68},{36,72},{36,86},{36,108},{37,43},{37,47},{37,48},
{37,49},{37,51},{37,57},{37,67},{37,90},{38,46},{38,59},{38,61},{38,66},{38,67},{38,68},{38,79},
{38,88},{38,91},{38,96},{38,101},{38,111},{39,45},{39,58},{39,62},{39,66},{39,67},{39,78},{39,89},
{39,92},{39,95},{39,99},{39,100},{39,112},{40,52},{40,65},{40,82},{40,90},{40,105},{41,50},{41,66},
{41,80},{41,83},{41,86},{42,46},{42,54},{42,61},{42,71},{42,72},{42,78},{42,85},{42,94},{42,100},
{42,112},{43,45},{43,55},{43,62},{43,75},{43,79},{43,84},{43,93},{43,101},{43,111},{44,49},{44,53},
{44,67},{44,104},{45,54},{45,58},{45,72},{45,75},{45,84},{45,94},{45,97},{45,111},{46,55},{46,59},
{46,72},{46,75},{46,85},{46,93},{46,98},{46,112},{47,69},{47,71},{47,73},{47,89},{47,91},{47,107},
{47,110},{48,70},{48,71},{48,74},{48,88},{48,92},{48,106},{48,109},{49,87},{49,96},{49,97},{49,98},
{49,100},{49,103},{49,104},{50,53},{50,60},{50,72},{50,95},{50,106},{50,107},{50,108},{50,110},
{50,111},{51,56},{51,84},{51,85},{51,90},{51,96},{51,104},{51,110},{51,112},{52,71},{52,83},{52,84},
{52,85},{52,87},{52,99},{52,109},{52,111},{53,72},{53,82},{53,84},{53,107},{53,111},{53,112},
{54,57},{54,62},{54,74},{54,76},{54,94},{54,102},{54,110},{55,56},{55,60},{55,61},{55,73},{55,77},
{55,93},{55,103},{55,109},{56,67},{56,72},{56,80},{56,88},{56,93},{56,94},{56,99},{57,65},{57,81},
{57,88},{57,91},{57,92},{57,95},{58,62},{58,69},{58,79},{58,93},{58,103},{58,109},{59,61},{59,70},
{59,78},{59,94},{59,102},{59,110},{60,68},{60,71},{60,82},{60,89},{60,91},{60,92},{61,69},{61,78},
{61,103},{61,106},{62,70},{62,79},{62,102},{62,107},{63,74},{63,77},{63,84},{63,88},{63,98},
{63,100},{64,73},{64,76},{64,85},{64,89},{64,97},{64,101},{65,68},{65,73},{65,74},{65,75},{65,77},
{65,78},{65,86},{65,108},{66,71},{66,77},{66,79},{66,87},{66,99},{66,104},{67,76},{67,78},{67,83},
{67,90},{67,96},{68,74},{68,78},{68,79},{68,82},{69,77},{69,81},{69,82},{69,83},{69,93},{69,103},
{69,105},{69,110},{70,76},{70,81},{70,82},{70,94},{70,102},{70,104},{70,108},{70,109},{71,80},
{71,87},{72,81},{73,77},{73,94},{73,95},{73,96},{73,102},{73,109},{74,76},{74,93},{74,99},{74,103},
{74,110},{75,82},{75,86},{76,93},{76,96},{76,99},{76,106},{77,94},{77,96},{77,99},{77,107},{78,84},
{78,89},{78,92},{78,95},{78,97},{79,85},{79,88},{79,91},{79,95},{79,98},{80,83},{80,90},{80,97},
{80,98},{80,99},{80,101},{80,102},{81,86},{81,95},{81,101},{81,103},{81,108},{82,100},{82,102},
{82,105},{83,98},{83,102},{83,103},{84,87},{84,92},{84,101},{84,112},{85,86},{85,90},{85,91},
{85,100},{85,111},{86,106},{86,111},{86,112},{87,106},{87,109},{87,110},{88,92},{88,100},{88,111},
{89,91},{89,101},{89,112},{90,107},{90,109},{90,110},{91,97},{92,98},{93,106},{94,107},{97,101},
{97,104},{97,105},{97,111},{98,100},{98,108},{98,112},{100,105},{100,108},{101,105},{101,108},
{102,104},{102,106},{103,104},{103,107},{106,110},{107,109}
\end{lstlisting}

\subsection*{(13,3)= 684}

\noindent $Z_{9}\rtimes_{5}Z_{76}$ generators $[6,5]\-[3,59];[5,64]\-[4,52];[0,5]\-[0,71];[0,30]\-[0,46];[1,72]\-[8,16];[8,38]\-[1,38];[0,38]$\\
\noindent $Z_{9}\rtimes_{9}Z_{76}$ generators $[1,5]\-[8,67];[2,11]\-[7,13];[4,5]\-[5,15];[1,74]\-[8,34];[2,33]\-[7,39];[0,68]\-[0,8];[0,38]$\\
\noindent $Z_{9}\rtimes_{17}Z_{76}$ generators $[6,60]\-[3,24];[4,74]\-[5,50];[4,32]\-[5,36];[8,18]\-[1,74];[0,15]\-[0,61];[6,19]\-[3,57];[0,38]$\\
\noindent $Z_{9}\rtimes_{25}Z_{76}$ generators $[2,48]\-[7,24];[1,22]\-[8,66];[8,0]\-[1,0];[5,40]\-[4,68];[3,5]\-[6,59];[4,63]\-[5,65];[0,38]$\\
\noindent $Z_{9}\rtimes_{61}Z_{76}$ generators $[5,58]\-[4,10];[8,71]\-[1,1];[5,43]\-[4,69];[6,42]\-[3,10];[0,72]\-[0,4];[6,65]\-[3,39];[0,38]$\\
\noindent $Z_{9}\rtimes_{73}Z_{76}$ generators $[5,62]\-[4,70];[2,32]\-[7,64];[8,26]\-[1,2];[7,12]\-[2,44];[1,31]\-[8,61];[8,43]\-[1,53];[0,38]$\\


\smallskip


\clearpage
\begin{thebibliography}{99}

\bibitem{BeEiOB05}
\newblock H. U. Besche, B. Eick, E. A. O'Brien.
\newblock The SmallGroups Library --- A Library of Groups of Small Order.
\newblock (2005) {An accepted and refereed GAP package, available also in MAGMA}

\bibitem{Com25} 
\newblock F. Comellas. 
\newblock New results on the degree-diameter problem for undirected graphs.  {\it  Electron. J. Graph Theory Appl.} {\bf 13 (1)} (2025), 211--215. 
% https://doi.org/10.5614/ejgta.2025.13.1.14

\bibitem{Com26} 
\newblock F. Comellas. 
\newblock Table of Large Degree/Diameter Graphs, Mendeley Data, V11 (2026), doi: 10.17632/d75dzbjd4k.11 
%   https://data.mendeley.com/datasets/d75dzbjd4k/11


\bibitem{CoMa13}
\newblock M. Conder, J. Ma.
\newblock Arc-transitive abelian regular covers of cubic graphs. 
\newblock {\it J. Algebra}. {\bf 387} (2013) 215--242.
% https://doi.org/10.1016/j.jalgebra.2013.02.035

\bibitem{Co11} 
\newblock M. Conder. 
\newblock Description of optimal Cayley graphs found by Marston Conder.
\newblock Page at Combinatoricswiki.org.  {\small \url{http://combinatoricswiki.org/wiki/Description_of_optimal_Cayley_graphs_found_by_Marston_Conder}.} (2011, updated 2019, accessed May30,2026)

\bibitem{Co26} 
\newblock M. Conder. 
\newblock Largest Cayley graphs with prescribed degree and diameter.
\newblock {Submitted, July 2026}. 

%\bibitem{D91}
%\newblock  M. J.  Dinneen.
%\newblock {\it Algebraic Methods for Efficient Network Constructions}.
%\newblock Master Thesis, Department of Computer Science,   University of Victoria, Victoria, B.C., Canada, 1991.

%\bibitem{Di98}
%\newblock  M. J.  Dinneen.
%\newblock {\it Group-theoretic methods for designing networks}.
%\newblock Centre for Discrete Mathematics and
%Theoretical Computer Science. Research Report Series CDMTCS-082, May 1998.
%University of Auckland. {\small \url{https://researchspace.auckland.ac.nz/server/api/core/bitstreams/173eba24-1d89-4c7d-9fe6-dfd34340d403/content}}


\bibitem{DH94}
\newblock M. J. Dinneen, P. Hafner.
\newblock New results for the degree/diameter problem.
\newblock{\em Nertworks}, {\bf 24} (1994) 359--367. (arXiv version 
{\small\url{https://arxiv.org/pdf/math/9504214.pdf}}

\bibitem{Hae80}
\newblock W. H. Haemers.
\newblock {\em Eigenvalue techniques in design and graph theory},
\newblock Mathematical Centre Tracts 121, Mathematisch Zentrum, Amsterdam,  1980.

\bibitem{Lo09} 
\newblock E. Loz,
\newblock {\em The degree-diameter and cage problems : a study in structural graph theory.}
\newblock  {\it Ph.D. Thesis. The University of Auckland  (2009)} \url{https://hdl.handle.net/2292/52008}. 

\bibitem{LoSi06} 
\newblock E. Loz, J. \v{S}ir\'{a}\v{n}.
\newblock {\em New record graphs in the degree-diameter problem.}
\newblock  {\it  Australas. J. Combin.}  {\bf 41} (2006) 63--80. 

\bibitem{MiSi13}
\newblock M. Miller,  J. \v{S}ir\'{a}\v{n}.
\newblock {Moore graphs and beyond:  A survey of the degree/diameter problem},
\newblock {\em Electron. J. Combin.}, {\bf DS14} (2013) 92 pp.

\bibitem{Mo06}
\newblock S.G. Molodsov.
\newblock Largest graphs of diameter 2 and maximum degree 6,
\newblock {\it General Theory of Information Transfer and Combinatorics, Springer, Berlin, Heidelberg (2006) 853--857. doi: 10.1007/11889342\_54}
\newblock {\small \url{https://link.springer.com/chapter/10.1007/11889342_54#citeas}}

\bibitem{Po13}
\newblock P. Poto\v{c}nik, P. Spiga, G. Verret,
\newblock {Cubic vertex-transitive graphs on up to 1280 vertices},
\newblock {\em J. Symb. Comput.}, {\bf 30} (2013) 465-477.
\newblock {\small \url{https://graphsym.net/}}


\bibitem{Ro13}
\newblock A. Rodr\'{\i}guez de los Santos.
\newblock {\it B\'usquedas masivas de grafos de gran orden con grado y di\'ametro acotados}.
\newblock Tesis de maestr\'{\i}a. Universidad de la Rep\'ublica (Uruguay). Facultad de Ingenier\'{\i}a., 2013. 
\newblock  {\small\url{https://graphsym.net/}}

\bibitem{Sa97}
\newblock M. Sampels.
\newblock Large networks with small diameter
\newblock {\it Proc. 23rd Int. Workshop on Graph-Theoretic Concepts in Computer Science (WG '97), LNCS 1335, pp. 288-302. Springer-Verlag, 1997 }
\newblock {\small\url{https://link.springer.com/content/pdf/10.1007/BFb0024505?pdf=chapter+toc}} 

\bibitem{Tw97} 
\newblock L. Twele. % Lutz Twele
\newblock {\em Effiziente Implementierung des Todd-Coxeter Algorithmus im Hinblick auf Grad/Durchmesser-Optimierung von knotentransitiven Graphen,} \newblock  {\it Diplomarbeit Universitat Oldenburg, 1997}

\bibitem{Wo96}
\newblock O. Wohlmuth. % Otto Wohlmuth.
\newblock A new dense group graph discovered by an evolutionary approach. 
\newblock {\it In Paralleles und Verteiltes Reehnen, Beitr\"age zum 4. Workshop \"uber wissensehaftliches Rechnen. Shaker Verlag, 1996. }


%%%%%%%%%%%%
\end{thebibliography}
\end{document}